\documentclass{article}

\usepackage{microtype}
\usepackage{booktabs} 

\usepackage{hyperref}

\usepackage[accepted]{icmla}

\usepackage{amsmath}
\usepackage{amssymb}
\usepackage{mathtools}
\usepackage{amsthm}

\usepackage[capitalize,noabbrev]{cleveref}

\usepackage[textsize=tiny]{todonotes}

\usepackage{enumitem,amsfonts,amssymb,mathtools,amsthm,bm,pifont,bbm,listings,xcolor,threeparttable,xcolor}
\usepackage{ulem,algorithmic}
\usepackage{graphicx}
\usepackage{subcaption}

\newtheorem{theorem}{Theorem}[section]
\newtheorem{proposition}[theorem]{Proposition}
\newtheorem{lemma}[theorem]{Lemma}

\theoremstyle{definition}
\newtheorem{definition}[theorem]{Definition}
\newtheorem{assumption}[theorem]{Assumption}
\theoremstyle{remark}

\def\a{\mathbf{a}}\def\b{\mathbf{b}}\def\c{\mathbf{c}}\def\g{\mathbf{g}}\def\p{\mathbf{p}}\def\r{\mathbf{r}}\def\v{\mathbf{v}}\def\x{\mathbf{x}}\def\y{\mathbf{y}}\def\z{\mathbf{z}}

\def\A{\mathbf{A}}\def\B{\mathbf{B}}\def\C{\mathbf{C}}\def\E{\mathbf{E}}\def\G{\mathbf{G}}\def\H{\mathbf{H}}\def\I{\mathbf{I}}

\def\M{\mathbf{M}}\def\N{\mathbf{N}}\def\P{\mathbf{P}}\def\Q{\mathbf{Q}}\def\U{\mathbf{U}}

\def\trans{^\mathsf{T}}  \def\zero{\mathbf{0}}    
\def\randn{\mathrm{randn}}
 \def\trans {^\mathsf{T}} 

\newcommand{\bfit}[1]{\textit{\textbf{#1}}}
\newcommand{\beq}{\begin{eqnarray}}
\newcommand{\eeq}{\end{eqnarray}}
\newcommand{\beqq}{\begin{equation}}
\newcommand{\eeqq}{\end{equation}}
\newcommand{\la}{\langle}
\newcommand{\ra}{\rangle}
\newcommand{\noi}{\noindent}
\newcommand{\nn}{\nonumber}
\def\noi{\noindent}
\def\nn{\nonumber}
\def\la{\langle}
\def\ra{\rangle}

\def\dist{{\rm{dist}}}
\newcommand{\step}[1]{\text{\ding{#1}}}

\def\varepsilons{\boldsymbol{\varepsilon}}

\def\lambdas{\boldsymbol{\lambda}}

\def\argmin{\operatorname*{arg\,min}}

\def\B{\mathrm{B}}
\def\N{\mathrm{B^c}}
\def\P{\mathbb{P}}

\def\B{\texttt{\textup{B}}}
\def\Bc{\B^c}
\def\UB{\mathrm{U}_{\B}}

\def\UBc{\mathrm{U}_{\Bc}}
\def\UBt{\mathrm{U}_{\B^t}}
\def\UBtt{\mathrm{U}_{{\B}^{t+1}}}

\usepackage{accents}
\usepackage{tikz}
\def\halfcheckmark{\tikz\draw[scale=0.32,fill=black](0,.35) -- (.25,0) -- (1,.7) -- (.25,.15) -- cycle (0.75,0.2) -- (0.77,0.2)  -- (0.3,0.6) -- cycle;}

\def\E{\mathbb{E}}

\def\Vup{\overline{\textup{V}}}
\def\Vlow{\underline{\textup{V}}}
\def\Aup{\overline{\textup{A}}}
\def\Alow{\underline{\textup{A}}}
\def\Sup{\overline{\textup{H}}}
\def\Slow{\underline{\textup{H}}}

\icmltitlerunning{Nonsmooth Sparsity Constrained Optimization}

\begin{document}

\twocolumn[

\icmltitle{Smoothing Proximal Gradient Methods for Nonsmooth Sparsity Constrained Optimization: Optimality Conditions and Global Convergence}



\icmlsetsymbol{equal}{}

\begin{icmlauthorlist}
\icmlauthor{Ganzhao Yuan}{yyy}
\end{icmlauthorlist}

\icmlaffiliation{yyy}{Peng Cheng Laboratory, China}

\icmlcorrespondingauthor{Ganzhao Yuan}{yuangzh@pcl.ac.cn}

\icmlkeywords{Machine Learning, ICML}

\vskip 0.3in
]



\printAffiliationsAndNotice{} 

\begin{abstract}

Nonsmooth sparsity constrained optimization encompasses a broad spectrum of applications in machine learning. This problem is generally non-convex and NP-hard. Existing solutions to this problem exhibit several notable limitations, including their inability to address general nonsmooth problems, tendency to yield weaker optimality conditions, and lack of comprehensive convergence analysis. This paper considers Smoothing Proximal Gradient Methods (SPGM) as solutions to nonsmooth sparsity constrained optimization problems. Two specific variants of SPGM are explored: one based on Iterative Hard Thresholding (SPGM-IHT) and the other on Block Coordinate Decomposition (SPGM-BCD). It is shown that the SPGM-BCD algorithm finds stronger stationary points compared to previous methods. Additionally, novel theories for analyzing the convergence rates of both SPGM-IHT and SPGM-BCD algorithms are developed. Our theoretical bounds, capitalizing on the intrinsic sparsity of the optimization problem, are on par with the best-known error bounds available to date. Finally, numerical experiments reveal that SPGM-IHT performs comparably to current IHT-style methods, while SPGM-BCD consistently surpasses them.


\end{abstract}

\section{Introduction}

This paper mainly focuses on the following nonsmooth sparsity constrained optimization problem (`$\triangleq$' means define):
\beq \label{eq:main}
\min_{\x \in \mathbb{R}^n}\,F(\x) \triangleq f(\x) + h(\A\x - \b),\,s.t.\,\|\x\|_0\leq s.
\eeq
\noi Here, $\A\in\mathbb{R}^{m\times n}$, $\b\in\mathbb{R}^{m}$, $s\in[n]$ is a positive integer, $f(\x):\mathbb{R}^{n} \mapsto \mathbb{R}$ is a smooth convex function, and $h(\y): \mathbb{R}^{m} \mapsto \mathbb{R}$ is a convex but not necessarily smooth function. For any vector $\c \in \mathbb{R}^{m}$ and any positive constant $\mu \in \mathbb{R}$, we assume that the following proximal operator of $h(\cdot)$ can be computed efficiently:
\beq \label{eq:proximal:operator}
\P_\mu(\c)  \triangleq \arg \min_{\y} h(\y) + \frac{1}{2\mu} \| \c - \y\|_2^2.
\eeq

Problem (\ref{eq:main}) captures a diverse range of applications in machine learning. To illustrate, nonsmooth functions including $h(\x) \triangleq \|\A\x - \b\|_{1}$, $h(\x) \triangleq \|\A\trans(\A\x - \b)\|_{\infty}$, and $h(\x) \triangleq \|\max(0,\A\x-\b)\|_1$ have been used in robust regression, Digzig selector computation, and support vector machines, respectively \cite{yuan2020jmlr}. Furthermore, Problem (\ref{eq:main}) covers a multitude of significant applications, such as sparse logistic regression \cite{bahmani2013greedy}, sparse censored regression \cite{bian2020smoothing}, impulse noise removal \cite{YuanG19}, sparse isotonic regression \cite{chen18d}, and sparse quantile regression \cite{bian2020smoothing}, as specific instances.

\begin{table*}[ht]
\centering
\caption{Comparison of existing nonsmooth sparsity constrained optimization methods. $T$ denotes the iteration counter, $L_F$ is the Lipschitz constant for $F(\x)$, and $\bar{\x}$ is the global optimal solution satisfying $\|\bar{\x}\|_0\leq s$.}
\scalebox{0.9}{\begin{tabular}{|c|c|c|c|}
\hline
  & General Nonsmooth &Optimality Conditions & Convergence Rate$^{a}$   \\
  \hline
\textbf{PDM} \cite{LuZ13} & \ding{52} & Lipschitz Stationary & Not Available  \\
 \hline
 \textbf{DIHT} \cite{yuan2020jmlr}  & \ding{56}   & Lipschitz Stationary & $\mathcal{O}(\frac{1}{\sqrt{T}})^{b}$ \\
\hline
\textbf{PSGD} \cite{liu2019one}&  \halfcheckmark $^{c}$ & Lipschitz Stationary & $\mathcal{O}(\frac{1}{T} ) + L_F\|\bar{\x}\| $   \\
\hline
\textbf{SPGM-IHT} [ours] & \ding{52} & Lipschitz Stationary & $\mathcal{O}(\frac{1}{T} ) + \mathcal{O}(\frac{\ln(T)}{T} ) + L_F\|\bar{\x}\|^{d}$ \\
\hline
\textbf{SPGM-BCD} [ours] &\ding{52} & Block-$k$ Stationary & $\mathcal{O}(\frac{1}{T} ) + \mathcal{O}(\frac{\ln(T)}{T} ) + 3 L_F \|\bar{\x}\|^{e}$ \\
\hline
\end{tabular}}
\smallskip 
\scalebox{0.85}{\begin{tabular}{@{}p{\linewidth}@{}}
\footnotesize Note $a$: The rate describes the decrease in objective values towards the optimum $F(\bar{\x})$, not in the distance to the optimal solution $\bar{\x}$. \\
\footnotesize Note $b$: \textbf{DIHT} establishes the convergence rate solely for the primal-dual gap, without addressing the primal convergence rate. \\
\footnotesize Note $c$: \textbf{PSGD} is less versatile, unable to solve nonsmooth problems when $h(\x)$ lacks Lipschitz continuity.  \\
\footnotesize Note $d$: The irreducible estimation error term $L_F\|\bar{\x}\|$ precisely aligns with the \textbf{PSGD} bound. Refer to Theorems \ref{theorem:IHT:stepsize:constant}, \ref{theorem:IHT:stepsize:diminishing}.\\
\footnotesize Note $e$: The irreducible estimation error term $3 L_F\|\bar{\x}\|$ is three times that of the \textbf{PSGD} bound. Refer to Theorems \ref{theorem:DEC:stepsize:constant},\ref{theorem:DEC:stepsize:diminishing}.
\end{tabular}}
\label{tab:intro}
\vspace{-10pt}
\end{table*}

Solving Problem (\ref{eq:main}) presents a challenge primarily due to the combinatorial nature of the cardinality constraint. A conventional approach involves replacing the non-convex $\ell_0$ norm with its convex relaxations, such as the $\ell_1$ norm \cite{candes2005decoding} and top-$k$ norm relaxation. However, studies have revealed that non-convex approximation techniques, such as the Schatten $\ell_p$ norm \cite{xu2012l, ZengLX16} and reweighted $\ell_1$ norm \cite{candes2008enhancing}, often yield superior accuracy compared to their convex counterparts \cite{zhang10a, YuanG19}. Furthermore, alternative strategies like multi-stage convex relaxation techniques have been introduced \cite{zhang10a, bi2014exact}, aiming to refine solutions obtained through convex methods. Recent efforts have primarily focused on directly minimizing the non-convex formulation in (\ref{eq:main}). Greedy pursuit methods \cite{bahmani2013greedy,tropp2007signal} selectively choose a variable coordinate to update, leading to optimality guarantees in certain scenarios. Iterative Hard Thresholding (IHT) methods \cite{BahmaniRajBoufounos2013, NguyenNW17} maintain sparsity by iteratively zeroing out small magnitude elements in a gradient descent fashion. Convergence rates and parameter estimation errors for IHT-style methods have been rigorously established under restricted smoothness and strong convexity conditions \cite{YuanLZ17, jain2014iterative}. The work of \cite{beck2013sparsity, Beck2016,beck2016minimization,beck2019optimization} introduced a novel optimality criterion based on coordinate-wise optimality for sparsity constrained optimization. It has been proven that this condition is stronger than the IHT-based optimality condition. Additionally, a new block coordinate optimality condition \cite{yuankdd2020, yuancvpr2019} was introduced for general sparse optimization, which is more powerful than the coordinate-wise optimality condition, encompassing it as a special case.

Another challenge in solving Problem (\ref{eq:main}) arises from the nonsmooth nature of the objective function. One widely adopted approach to address this issue is the Alternating Direction Method of Multipliers (ADMM) \cite{HeY12}. ADMM introduces dual variables to address linear constraints, iteratively optimizing primal variables with other primal and dual variables kept static, and employs a gradient ascent strategy to update the dual variables. However, it has been noted in \cite{LuZ13} that ADMM often yields unsatisfactory solution quality. This observation has motivated the exploration of Penalty Decomposition Methods (PDM) for solving generally nonlinear sparsity constrained optimization problems \cite{LuZ13}. Additionally, Projective Subgradient Descent (PSGD) methods have been proposed for solving nonsmooth one-bit compressed sensing problems \cite{liu2019one}, operating by iteratively projecting the intermediate solution onto the nonconvex sparsity constraint after each sub-gradient descent update. Furthermore, Dual Iterative Hard Thresholding (DIHT) \cite{yuan2020jmlr} applies projective subgradient methods to the dual of sparsity constraint optimization problems, offering proven guarantees on primal-dual gap convergence and sparsity recovery. Their duality theory establishes sufficient and necessary conditions for solving the original non-convex problem equivalently or approximately through a concave dual approach.

In summary, existing methods for solving Problem (\ref{eq:main}) exhibit three main limitations. \bfit{(i)} Inability to handle general nonsmooth problems. Block decomposition \cite{yuankdd2020} and dual IHT \cite{yuan2020jmlr} methods are limited to smooth sparsity constrained problems, while PSGD methods are restricted to objectives that are Lipschitz continuous. These methods struggle with general non-Lipschitz problems, which are better addressed by penalty decomposition \cite{LuZ13} or smoothing proximal gradient methods \cite{bian2020smoothing,chen2012smoothing}. \bfit{(ii)} Tendency to yield weaker optimality conditions. Predominantly relying on IHT, current methods often result in suboptimal optimality guarantees and subpar practical accuracy \cite{beck2013sparsity,yuankdd2020,yuan2023coordinate}. \bfit{(iii)} Lack of comprehensive convergence analysis. Despite the integration of IHT-style methods into penalty decomposition methods \cite{LuZ13}, a thorough convergence analysis is lacking. Additionally, the duality theory in \cite{yuan2020jmlr} is constrained by its assumption of smooth objective functions, as evident in Theorem 15 and Theorem 17 in \cite{yuan2020jmlr}.

To address these limitations, this paper introduces Smoothing Proximal Gradient Methods (SPGM) for nonsmooth sparsity constrained optimization, featuring two SPGM variants: \textbf{SPGM-IHT} and \textbf{SPGM-BCD}. These methods, rooted in smoothing techniques, tackle a wide range of nonsmooth problems, with \textbf{SPGM-BCD} ensuring superior optimality conditions. We also establish the convergence rate of both methods towards the global optimum. Our theoretical bounds, which leverage the inherent sparsity of the optimization problem, match the best-known error bounds \cite{liu2019one} currently available (details in Table \ref{tab:intro}).

\textbf{Contributions.} The contributions of this paper are threefold. \bfit{(i)} Algorithmically, we explore Smoothing Proximal Grdient Methods (SPGM) for solving Problem (\ref{eq:main}), including SPGM based on Iterative Hard Thresholding (\textbf{SBCD-IHT}) and SPGM based on Block Coordinate Decomposition (\textbf{SPGM-BCD}) (See Section \ref{sect:SPGM}). We offer smooth and optimality analysis for the smoothing reformulation problem, demonstrating that \textbf{SPGM-BCD} attains stronger stationary points compared to existing solutions (see Section \ref{sect:smooth:analysis}). \bfit{(ii)} Theoretically, we develop novel theories to analyze the convergence rate of both \textbf{SPGM-IHT} and \textbf{SPGM-BCD} (See Section \ref{sect:rate}). \bfit{(iii)} Empirically, we have conducted experiments on two nonsmooth sparsity constrained optimization tasks to show the superiority of our methods (See Section \ref{sect:exp}).



\textbf{Notations.} All vectors are column vectors, with superscript $\trans$ indicating transpose. For a vector $\x\in \mathbb{R}^n$, $\x_i$ represents its $i$-th component for any $i\in [n] \triangleq \{1,2,...,n\}$. The Euclidean inner product between vectors $\x$ and $\x'$ is expressed as $\la \x,\x' \ra$ or $(\x')\trans\x$. The identity matrix in $\mathbb{R}^{n\times n}$ is denoted by $\I_n$. $\|\A\|$ represents the spectral norm of $\A$. The notations $\C \succeq \mathbf{0}$ and $\C \succ \mathbf{0}$ indicate positive semidefiniteness and definiteness of $\C$, respectively. For any $\C$ with $\C \succeq \mathbf{0}$, we define $\|\x\|_{\C}\triangleq \sqrt{\x\trans\C\x}$ as a generalized vector norm, and denote $\lambdas_{\max}(\C)$ and $\lambdas_{\min}(\C)$ as respectively the largest and smallest eigenvalue of $\C$. If $\beta$ is a constant, $\beta^t$ refers to its $t$-th power, while if $\beta$ is an optimization variable, $\beta^t$ signifies the value in the $t$-th iteration. The subdifferential of the function $h:\mathbb{R}^{n} \mapsto (-\infty,+\infty]$ at $\x$, defined as $\partial h(\x) \triangleq \{\g: h(\y) \geq h(\x) + \la\g, \y - \x \ra \}$, includes all subgradients of $h(\x)$. The squared distance between sets $\Xi$ and $\Xi'$ is defined as $\dist^2(\Xi,\Xi') \triangleq \inf_{\v \in \Xi, \v' \in \Xi'}\|\v - \v'\|_2^2$.

For a set $\B \in \mathbb{N}^k$ containing $k$ unique integers selected from $\{1, 2, ..., n\}$, we define $\Bc \triangleq \{1,2,...,n\}\setminus \B$, and denote $\C_{\B\B}$ as the sub-matrix of $\C$ indexed by $\B$. $C_n^k$ counts the combinations to select $k$ items from $n$ without repetition. $\Omega_n^k \triangleq \{\mathcal{\B}_{(1)},\mathcal{\B}_{(2)},...,\mathcal{\B}_{(C_n^k)}\}$ represents the set of all index vector combinations for this selection, with each $\mathcal{\B}_{(i)} \in \mathbb{N}^k$.

\section{Smoothing Proximal Gradient Methods} \label{sect:SPGM}

This section explores Smoothing Proximal Gradient Methods (SPGM) for Problem (\ref{eq:main}), detailing two versions: SPGM-IHT, using Iterative Hard Thresholding \cite{blumensath2008gradient,BLUMENSATH2009265}, and SPGM-BCD, employing Block Coordinate Decomposition \cite{yuankdd2020,yuancvpr2019}.

In the sequel of this paper, we impose the following assumptions on Problem (\ref{eq:main}).

\quad\begin{assumption} \label{ass:1}
The functions $f(\cdot)$ and $h(\cdot)$ are Lipschitz continuous with some constants $L_f$ and $L_h$, satisfying $\| \nabla f(\x)\| \leq L_f$ for all $\|\x\|_0 \leq s$ and $\| \partial h(\y)\|\leq L_h$ for all $\y\in\mathbb{R}^m$. Consequently, $F(\x)$ is Lipschitz continuous with constant $L_F \triangleq L_f + \|\A\| L_h$.
\end{assumption}

\quad\begin{assumption}\label{ass:2}
The function $f(\cdot)$ is restricted $V_{s}$-strongly convex and restricted $M_{s}$-smooth, satisfying:
\beq
\tfrac{V_{s}}{2}\|\x-\x'\|_2^2 \leq  \mathcal{Q}(\x,\x')  \leq  \tfrac{M_s}{2}\|\x-\x'\|_2^2\nn
\eeq
\noi for all $\|\x\|_0 \leq s$ and $\|\x'\|_0 \leq s$, where $\mathcal{Q}(\x,\x')\triangleq f(\x') - f(\x) - \la \x'-\x,\nabla  f(\x)\ra$. Additionally, a symmetric matrix $\tilde{\M}\in\mathbb{R}^{n\times n}$ exists, fulfilling $\zero \prec V_{s} \I_n \preceq\tilde{\M} \preceq M_s \I_n$ and
\beq
\mathcal{Q}(\x,\x')  \leq  \tfrac{1}{2}\|\x-\x'\|_{\tilde{\M}}^2 \label{eq:inequality:H}
\eeq
\noi for all $\|\x\|_0 \leq s$ and $\|\x'\|_0 \leq s$.
\end{assumption}

\quad\begin{assumption}\label{ass:3}
\noi A constant $A_{s}>0$ exists, ensuring $\|\A(\x-\x')\|  \leq A_{s}\|\x-\x'\|$ for all $\x\in\mathbb{R}^n,\x'\in\mathbb{R}^n$ with $\|\x\|_0 \leq s$, $\|\x'\|_0 \leq s$.
\end{assumption}

\noi \textbf{Remarks}. \bfit{(i)} Assumptions \ref{ass:1}, \ref{ass:2}, and \ref{ass:3} are broadly applicable, meeting conditions of various applications like robust regression and support vector machines (see \cite{YuanLZ17}). \bfit{(ii)} Assumption \ref{ass:3} is less stringent than $\|\A(\x-\x')\|  \leq\|\A\|\|\x-\x'\|$. \bfit{(iii)} Common choices for nonsmooth $h(\y)$ include $\{\|\y\|_1$, $\|\max(0,\y)\|_1$, $\|\y\|_{\infty}\}$, with their corresponding $L_h$ values being $\{\sqrt{m},\sqrt{m},1\}$, respectively. \bfit{(iv)} When $f(\x)$  takes the form of a quadratic function with $f(\x)\triangleq \frac{1}{2}\x\trans \hat{\Q} \x + \x\trans\hat{\p}$ for some $\hat{\Q}\in\mathbb{R}^{n\times n}$ and $\hat{\p}\in\mathbb{R}^n$, Inequality (\ref{eq:inequality:H}) holds with $\mathcal{Q}(\x,\x')  =  \tfrac{1}{2}\|\x-\x'\|_{\tilde{\M}}^2$, where $\tilde{\M}=\hat{\Q}$.

Introducing a new variable $\y\in\mathbb{R}^m$, we reframe Problem (\ref{eq:main}) as: $\min_{\x,\y}\,f(\x)  +  h(\y),\,s.t.\,\A\x- \b = \y,\,\|\x\|_0\leq s$. In \textbf{SPGM}, a smoothing parameter $\mu \rightarrow 0$ is incorporated to penalize the squared error in the linear constraints, leading to the subsequent minimization problem:
\beq \label{eq:penalty:J}
\min_{\x,\y}\,\mathcal{J}(\x,\y;\mu) \triangleq \mathcal{R}(\x,\y;\mu) +   h(\y) + \delta(\x),\,\nn \\
\text{where}\,\mathcal{R}(\x,\y;\mu) \triangleq f(\x) + \tfrac{1}{2\mu} \| \A\x - \b - \y\|_2^2,
\eeq
\noi and $\delta(\x) \triangleq { {\tiny \left\{\begin{array}{ll}0, &\hbox{$\|\x\|_0 \leq s$} \\
\infty, &\hbox{else} \end{array}\right.}}$. In each iteration, we employ proximal point strategies to alternatively minimize \textit{w.r.t.} $\x$ and $\y$ \cite{tseng2009coordinate}. Notably, SPGM is closely related to alternating minimization methods, block coordinate descent methods \cite{xu2013block}, and penalty decomposition methods \cite{LuZ13} in the literature.

\noi $\blacktriangleright$ $\x$-subproblem. Keeping parameters $\y^t$ and $\mu^t$ constant at their current values, we minimize $\mathcal{J}(\x,\y^t;\mu^t)$ \textit{w.r.t.} $\x$, resulting in the following optimization problem:
\beq
\min_{\x}\,\mathcal{R}(\x,\y^t;\mu^t),\,s.t.\,\|\x\|_0\leq s.\nn
\eeq
\noi The function $\mathcal{R}(\x,\y^t;\mu^t)$ is differentiable in $\x$, with its gradient at $\x^t$ given by:
\beq
\nabla_{\x} \mathcal{R}(\x^t,\y^t;\mu^t)  =  \nabla f(\x^t) + \tfrac{1}{\mu^t}\A\trans (\A\x^t - \b - \y^t)  \triangleq \mathbf{r}^t.\nn
\eeq
\noi To solve the $\x$-subproblem, we consider state-of-the-art sparse optimization methods, including the IHT strategy \cite{YuanLZ17,yuan2020jmlr,jain2014iterative,lu2014iterative} and the BCD strategy \cite{yuankdd2020}.

\begin{algorithm}[!t]
\caption{ {\bf The  Smoothing Proximal Gradient Methods based on IHT or BCD strategy} }
\begin{algorithmic}[0]

\STATE Input: working set size $k\in[n]$, the proximal point parameters $\theta>0$, $\theta_1>0$, $\theta_2>0$, initial feasible solution $\x^1$, an initial parameter $\mu^{1}$. 

\FOR{$t = 1$ to $T$}
\STATE (S1) Solve the following $\x$-subproblem using IHT or BCD strategy.

\quad $\blacktriangleright$ Option I (\textbf{IHT}): Solve the following problem globally \cite{blumensath2008gradient}:
\beq \label{eq:subprob:x:1}
\begin{split}
&\x^{t+1} \in  \argmin_{\|\x\|_0\leq s} \, \dot{\mathcal{M}}(\x,\x^t,\y^t;\mu^t) \triangleq R^t  \\
&+  \tfrac{ H^t }{2} \|\x-\x^t\|_2^2 +  \la \x-\x^t,\r^t \ra,
\end{split}
\eeq
\noi where $H^t\triangleq A_s^2 / \mu^t + M_s + \theta \in \mathbb{R}$, and $R^t\triangleq\mathcal{R}(\x^t,\y^t;\mu^t)$.

\quad$\blacktriangleright$ Option II (\textbf{BCD}): Use a random or/and a greedy method to find a working set $\B^t$ of size $k$ for the $t$-th iteration. Let $\B=\B^t$ and $\B^c \triangleq \{1,...,n\}\setminus \B$. Solve the following problem globally \cite{yuankdd2020}:
\beq \label{eq:subprob:x:2}
\begin{split}
&\x^{t+1} \in \argmin_{\|\x\|_0\leq s,\x_{\B^c} = \x^t_{\B^c}} \, \ddot{\mathcal{M}}(\x,\x^t,\y^t;\mu^t) \triangleq R^t   \\
& ~~~~~+ \tfrac{ 1 }{2} \|\x-\x^t\|_{ \H^t}^2+ \la \x-\x^t,\r^t \ra,  \\
\end{split}
\eeq
\noi where $\H^t\triangleq (\A\trans \A + \theta_1 \I_n) /\mu^t +  \tilde{\M} + \theta_2 \I_n\in \mathbb{R}^{n\times n}$, and $R^t\triangleq\mathcal{R}(\x^t,\y^t;\mu^t)$.


\STATE (S2) Solve the following $\y$-subproblem:
 \beq \label{eq:subprob:y}
 \y^{t+1} &=& \argmin_{\y} \mathcal{J}(\x^{t+1},\y;\mu^{t}) \nn\\
 &=& \P_{\mu^t}(\A\x^{t+1} - \b)
 \eeq
\STATE (S3) Choose a new parameter $\mu^{t+1}$ with $\mu^{t+1} \leq \mu^{t}$.
\ENDFOR
\end{algorithmic}
\label{algo:main}
\end{algorithm}

\hspace{7pt} \boxed{\text{IHT strategy}.} We observe it always holds that:
\beq \label{eq:major:1}
\mathcal{R}(\x,\y^t;\mu^t) \leq \dot{\mathcal{M}}(\x,\x^t,\y^t;\mu^t)
\eeq
\noi for all $\|\x\|_0\leq s$, where $\dot{\mathcal{M}}(\x,\x^t,\y^t;\mu^t)$ is defined in Equation (\ref{eq:subprob:x:1}), and $\theta>0$ is a constant. The IHT strategy aims to minimize the majorization function $\dot{\mathcal{M}}(\x,\x^t,\y^t;\mu^t)$, while adhering to the sparsity constraint. This approach simultaneously reduces the objective function and identifies the active variables, as indicated by the update in (\ref{eq:subprob:x:1}). We note that (\ref{eq:subprob:x:1}) is equivalent to the following problem:
\beq \label{eq:opt:xt1}
\x^{t+1} \in \arg \min_{\|\x\|_0 \leq s}\, \tfrac{1}{2}\| \x-\x^t_+  \|_2^2  =  \Pi_s(\x^t_+),
\eeq
\noi where $\x^t_+ \triangleq \x^t - \r^t / H^t$, $H^t\triangleq A_s^2 / \mu^t + M_s + \theta$, and $\Pi_s(\x)$ is an operator that sets all but the largest (in magnitude) $s$ elements of $\x$ to zero.

\hspace{7pt} \boxed{\text{BCD strategy}.} We notice the following inequality consistently holds:
\beq \label{eq:major:2}
\mathcal{R}(\x,\y^t;\mu^t) \leq \ddot{\mathcal{M}}(\x,\x^t,\y^t;\mu^t),
\eeq
\noi for all $\|\x\|_0\leq s$, where $\ddot{\mathcal{M}}(\x,\x^t,\y^t;\mu^t)$ is defined in Equation (\ref{eq:subprob:x:2}), and $\{\theta_1,\theta_2\}$ are given positive constants. The \textbf{BCD} strategy aims to minimize the majorization function $\ddot{\mathcal{M}}(\x,\x^t,\y^t;\mu^t)$ using a block coordinate fashion. It employs either a random method or a greedy method to select a subset of coordinates of size $k$ as the working set $\B$, where $k\geq 2$. It then conducts a global combinatorial search over this working set, based on the quadratic majorization function, as indicated by the update in (\ref{eq:subprob:x:2}). Problem (\ref{eq:subprob:x:2}) can be equivalently rewritten as: $\x_{\B}^{t+1}  \in \arg \min_{\z_{\B}}  \ddot{\mathcal{M}}(\UB \z_{\B} + \UBc \x^t_{\Bc},\x^t,\y^t;\mu^t) + \delta(\UB \z_{\B} + \UBc \x^t_{\Bc})$, where $\Bc \triangleq \{1,...,n\} \setminus \B$, $\UB \in \mathbb{R}^{n\times k}$, $\UBc \in \mathbb{R}^{n\times (n-k)}$, and
\beq
[\UB]_{ji}=\left\{
             \begin{array}{ll}
               1, & \hbox{$\B_i=j$;} \\
               0, & \hbox{else.}
             \end{array}
           \right.,\,[\UBc]_{ji} = \left\{
             \begin{array}{ll}
               1, & \hbox{$\Bc_i=j$;} \\
               0, & \hbox{else.}
             \end{array}
           \right..\nn
\eeq
\noi We have: $\x=(\U_{\N}\U_{\N}\trans+\U_{\B}\U_{\B}\trans)\x=\U_{\B}\x_{\B} + \U_{\N}\x_{\N}$, and $\x_{\B}=\U_{\B}\trans\x$. Thus, Problem (\ref{eq:subprob:x:2}) reduces to the following problem:
\beq \label{eq:BCD:2}
&\x_{\B}^{t+1}  \in \arg \min_{\z_{\B} \in\mathbb{R}^k}\, \frac{1}{2}(\z_{\B} - \x_{\B}^t)\trans [\H^t_{\B\B}](\z_{\B} - \x_{\B}^t),\nn\\
& + \la \z_{\B} - \x_{\B}^t,\, \r^t_{\B} \ra,\,s.t.\,\|\z_{\B}\|_0 + \|\x^t_{\B^c}\|_0\leq s,
\eeq
\noi where $\H^t \triangleq (\A\trans \A + \theta_1 \I_n) /\mu^t +  \tilde{\M} + \theta_2 \I_n$. Problem (\ref{eq:BCD:2}) involves $k$ unknown decision variables, and can be tackled by solving a set of $2^k$ linear equations. The BCD strategy combines the efficacy of combinatorial search methods with the efficiency of coordinate descent methods, allowing it to efficiently identify stronger stationary points than the IHT strategy when minimizing smooth functions under sparsity constraints \cite{yuankdd2020,yuancvpr2019}.

\noi $\blacktriangleright$ $\y$-subproblem. With the parameters $\x^{t+1}$ and $\mu^t$ fixed at their current estimates, we encounter an optimization problem \textit{w.r.t.} $\y$ as in Equation (\ref{eq:subprob:y}), which is equivalent to the computation of the proximal operator as described in Equation (\ref{eq:proximal:operator}).

We summarize the SPGM algorithm in Algorithm \ref{algo:main}.

\section{Smooth and Optimality Analysis}
\label{sect:smooth:analysis}

This section provides smooth and optimality analysis for the smoothing function as in Problem (\ref{eq:penalty:J}).

\subsection{Smooth Analysis}

Problem (\ref{eq:penalty:J}) becomes equivalent to the original optimization problem in (\ref{eq:main}) as $\mu \rightarrow 0$. This equivalence is expressed as:
\beq
[\min_{\x}F(\x) + \delta(\x)] \equiv [\min_{\x,\y} \lim_{\mu \rightarrow 0} \mathcal{J}(\x,\y;\mu) ]. \nn
\eeq
\noi Thus, we conduct a smooth analysis for Problem (\ref{eq:penalty:J}). By eliminating $\y$, Problem (\ref{eq:penalty:J}) simplifies to:
\beq \label{eq:def:GGG}
\begin{split}
 &\min_{\x}\,\mathcal{G}(\x;\mu) = f(\x) +  h(\P_{\mu}(\A\x - \b)) \nn\\
 &~~+ \tfrac{1}{2\mu} \| \A\x - \b - \P_{\mu}(\A\x - \b)\|_2^2,\,s.t.\,\|\x\|_0\leq s.
\end{split}
\eeq
\noi $\mathcal{G}(\x;\mu)$ is smooth \textit{w.r.t.} $\x$ and its gradient is given by:
\beq
\nabla_{\x}\mathcal{G}(\x;\mu) = \nabla f(\x) +  \tfrac{1}{\mu} \A\trans (\A\x - \b - \P_{\mu}(\A\x - \b)).\nn
\eeq

\noi We have the following useful lemmas \footnote{All proofs can be found in the Appendix.}.

\begin{lemma} \label{lemma:lip:variable:mu}
\noi (Proof in Appendix \ref{app:lemma:lip:variable:mu}) Fix $\x$ with $\|\x\|_0\leq s$. The function $\psi(\mu) \triangleq \mathcal{G}(\x;\mu)$ is decreasing and $(\tfrac{1}{2}L_h^2)$-Lipschitz continuous for all $\mu>0$. In other words, for all $0< \mu_1 < \mu_2$, we have: $0 \leq \frac{\psi(\mu_1) - \psi(\mu_2) }{\mu_2 - \mu_1}  \leq \tfrac{1}{2}L_h^2$.

\end{lemma}

\begin{lemma} \label{lemma:lip:variable:x}
\noi (Proof in Appendix \ref{app:lemma:lip:variable:x}) Fix $\mu>0$. For all $\x$ with $\|\x\|_0\leq s$, we have:

\begin{enumerate}[label=\textbf{(\alph*)}, leftmargin=22pt, itemsep=1pt, topsep=1pt, parsep=0pt, partopsep=0pt]

\item It holds that: $F(\x)-\tfrac{\mu}{2}L_h^2 \leq \mathcal{J}(\x,\P_{\mu}(\A\x - \b);\mu) = \mathcal{G}(\x;\mu)  \leq F(\x)$.

\item It holds that: $\|\partial F(\x)\| \leq L_F$, $\|\nabla \mathcal{G}(\x;\mu)\| \leq L_F$, where $L_F \triangleq L_f + L_h \|\A\|$.

\item $\mathcal{G}(\x;\mu)$ is restricted $V_{s}$-strongly convex and restricted $(M_{s} + A_s \|\A\|/\mu )$-smooth.

\end{enumerate}

\end{lemma}

\noi \textbf{Remarks}. \bfit{(i)} Lemmas \ref{lemma:lip:variable:mu} and \ref{lemma:lip:variable:x} can be derived using Assumptions \ref{ass:1}, \ref{ass:2}, and \ref{ass:3}, along with the optimality of the proximal operator $\P_{\mu}(\c)$ for any $\c$. \bfit{(ii)} The inequalities in Lemma \ref{lemma:lip:variable:mu} and \textbf{Part (a)} of Lemma \ref{lemma:lip:variable:x} are closely linked to smooth approximation functions as discussed in \cite{chen2012smoothing} and the Moreau-Yosida approximation \cite{bauschke2011} in the literature. These properties play a crucial role in the development of smoothing methods for nonsmooth optimization. \bfit{(iii)} Given that $\mathcal{G}(\x^{t};\mu^{t-1})$ serves as a smooth approximation function for $F(\x^t)$, we can assess the convergence rate of $F(\x^t)$ by estimating the convergence rate of $\mathcal{G}(\x^{t};\mu^{t-1})$.


\subsection{Optimality Analysis}

To provide optimality analysis for SPGM, we begin by introducing some fundamental definitions.

\begin{definition}
(Basic Stationary Point) A solution $\check{\x}$ is a basic stationary point if the following condition is met: $F(\check{\x}) = \min_{\x}\,F(\x),\, s.t.\,[\x]_{\mathrm{J}^c }=\mathbf{0}$. Here, $\mathrm{J}^c\triangleq \{1,...,n\} \setminus \mathrm{J}$, where $\mathrm{J}$ represents the known support set of the solution $\check{\x}$ with $|\mathrm{J}|\leq s$.

\end{definition}

\noi \textbf{Remarks}. The basic stationary point implies that the solution attains global optimality when the support set is restricted \cite{beck2013sparsity}.


\begin{definition} \label{def:point:L}
(Lipschitz Stationary Point) Fix $\mu > 0$ as a sufficiently small constant. A solution $(\dot{\x},\dot{\y})$ is a Lipschitz stationary point if the following condition holds:
\beq
\dot{\y} & \in  & \arg \min_{\y}\,\mathcal{J}(\dot{\x},\y;\mu) ,  \nn\\
\dot{\x}  &\in & \arg \min_{\x}\, \dot{\mathcal{M}}(\x,\dot{\x},\dot{\y};\mu) + \delta(\x) , \nn
\eeq
\noi where $\dot{\mathcal{M}}(\x,\x^t,\y^t;\mu^t)$ is defined in Equation (\ref{eq:subprob:x:1}).

\end{definition}
\noi \textbf{Remarks}. The Lipschitz stationary point states that if we minimize the smoothing function $\mathcal{J}(\dot{\x},\y;\mu)$ over $\y$ and the majorization function $\dot{\mathcal{M}}(\x,\dot{\x},\dot{\y};\mu)$ over $\x$, the quality of the solution $(\dot{\x},\dot{\y})$ cannot be further improved.

\begin{definition} \label{def:point:K}
(Block-$k$ Stationary Point) Fix $\mu>0$ as a sufficiently small constant. We denote $\B^c \triangleq \{1,...,n\}\setminus \B$. A solution $(\ddot{\x},\ddot{\y})$ is a block-$k$ stationary point if the following condition is met:
\beq
\ddot{\y}  & \in   & \arg \min_{\y}\,\mathcal{J}(\ddot{\x},\y;\mu) ,  \nn\\
\ddot{\x}_{\B} & \in  & \argmin_{\z_{\B},\|\z_{\B}\|_0 + \|\ddot{\x}_{\B^c}\|_0\leq s}\, \ddot{\mathcal{M}}(\U_{\B}\z_{\B} + \U_{\B^c}\ddot{\x}_{\B^c},\ddot{\x},\ddot{\y};\mu)  \nn
\eeq
\noi for all $\B\in \Omega_n^k$. Here, $\Omega_n^k\triangleq\{\mathcal{\B}_{(i)}\}_{i=1}^{C_n^k}$ denotes all the combinations of the index vector choosing $k$ items from $n$ without repetition, and $\ddot{\mathcal{M}}(\x,\x^t,\y^t;\mu^t)$ is defined in Equation (\ref{eq:subprob:x:2}).

\end{definition}

\noi \textbf{Remarks}. \bfit{(i)} Block-$k$ stationary point capture more intrinsic structures of the nonconvex problem than Lipschitz stationary points, and it holds that $\ddot{\mathcal{M}}(\x,\x^t,\y^t;\mu^t)\leq \dot{\mathcal{M}}(\x,\x^t,\y^t;\mu^t)$ for all $\x$. \bfit{(ii)} Deterministically finding a block-$k$ stationary point requires evaluating $C_n^k$ subproblems, which can be time-consuming. However, using a random strategy to select the working set $\B$ from the $C_n^k$ combinations allows for an expected block-$k$ stationary point.

The following proposition states the relation between different types of the stationary point above.

\begin{proposition} \label{proposition:hierarchy}
\textbf{Optimality Hierarchy \cite{yuankdd2020}.} We denote the sets $\{\check{\x}\}$ (basic stationary points), $\{\dot{\x}\}$ (Lipschitz stationary points), $\{\ddot{\x}_{[k]}\}$ (block-$k$ stationary points), and $\{\bar{\x}\}$ (global optimal points). The following relation holds for all $2\leq k\leq n-1$:
\beq
\{\bar{\x}\} \equiv \{\ddot{\x}_{[n]}\} \subseteq\{\ddot{\x}_{[k+1]}\} \subseteq \{\ddot{\x}_{[k]}\} \subseteq  \{\dot{\x}\}\subseteq \{\check{\x}\} . \nn
\eeq

\end{proposition}

We establish the optimality hierarchy among the optimality conditions by directly applying the results of Proposition 1 in \cite{yuankdd2020}, which addresses the minimization of smooth functions under sparsity constraints.

\section{Convergence Analysis}
\label{sect:rate}

In this section, we develop novel theories to analyze the convergence rate of \textbf{SPGM-IHT} and \textbf{SPGM-BCD}.

In our analysis, we consider two strategies for updating $\mu^t$ for all $t = 1, 2, \ldots, \infty$.

\begin{itemize}
  \item $\mu^t = \bar{\mu}$, where $\bar{\mu}>0$ is a sufficiently small constant.
  \item $\mu^t = \tfrac{\eta}{t+t_0}$, where $\eta>0$ and $t_0\geq 1$ are constants.
\end{itemize}


\noi We notice the following relation between $\nabla_{\x} \mathcal{G}(\x^t;{\mu^{t-1}})$ and $\nabla_{\x}\mathcal{R}(\x^t,\y^t;{\mu^{t}})$:
\beq \label{eq:relation}
\begin{split}
\g^t\triangleq&\, \nabla \mathcal{G}(\x^t;{\mu^{t-1}})  = \nabla_{\x}\mathcal{R}(\x^t,\y^{t+1};{\mu^{t-1}})  \nn \\
=&\, \nabla f(\x^{t}) +  \tfrac{1}{\mu^{t-1}} \A\trans (\A\x^{t} - \b - \y^{t} )\nn\\
=&\,\underbrace{\nabla_{\x}\mathcal{R}(\x^t,\y^t;{\mu^{t}})}_{\triangleq\,\r^t} + \underbrace{(\tfrac{1}{\mu^{t-1}} - \tfrac{1}{\mu^t})   \A\trans (\A\x^{t} - \b - \y^t)}_{ \triangleq \,\varepsilons^t  }
\end{split}
\eeq

\noi We derive the following results for both \textbf{SPGM-IHT} and \textbf{SPGM-BCD}.

\begin{lemma} \label{lemma:first:bound}
\noi (Proof in Appendix \ref{app:lemma:first:bound}) For all $t = 1,2,...,\infty$, we have:

\begin{enumerate}[label=\textbf{(\alph*)}, leftmargin=22pt, itemsep=1pt, topsep=1pt, parsep=0pt, partopsep=0pt]

\item $\| \A\x^{t}-\y^{t} - \b   \| \leq L_h \mu^{t-1}$.

\item $\|\y^{t+1}-\y^{t}\| \leq  \|\A\|  \| \x^{t+1}  -\x^{t}\|  + 2L_h \mu^{t-1}$.

\item $ \|\r^t\| \leq \begin{scriptsize}\left\{
  \begin{array}{ll}
L_F\triangleq L_f + L_h \|\A\|, & \hbox{$\mu^t=\bar{\mu}$;} \\
 L'_F \triangleq L_f + \frac{t_0+1}{t_0}L_h \|\A\|   , & \hbox{$\mu^t =\eta / (t+t_0)$}
  \end{array}
\right.\end{scriptsize}$.


\item $\|\varepsilons^t\|  \leq(\tfrac{\mu^{t-1}}{\mu^t} - 1) \|\A\|  L_h$.


\item $\mathcal{J}(\x^{t},\y^{t};\mu^{t})-\mathcal{J}(\x^{t},\y^{t};\mu^{t-1}) \leq \Psi^t$, where $\Psi^t\triangleq \tfrac{L_h^2}{2}   (\tfrac{(\mu^{t-1})^2}{\mu^{t}} - {\mu^{t-1}} )$.

\item $ [\sum_{t=1}^{\infty}\Psi^t] \leq  \begin{scriptsize}\left\{
  \begin{array}{ll}
0, & \hbox{$\mu^t=\bar{\mu}$;} \\
 \eta L_h^2   , & \hbox{$\mu^t =\eta / (t+t_0)$}
  \end{array}
\right.\end{scriptsize}$.

\end{enumerate}

\end{lemma}

\noi \textbf{Remarks}. \bfit{(i)} Given our choices of $\mu^t$ and the fact that $\| \A\x^{t}-\y^{t} - \b   \| \leq L_h \mu^{t-1}$, it follows that: $\|\A\x^{t}-\y^{t} - \b\| \rightarrow0$. \bfit{(ii)} We notice that $\|\y^{t+1} - \y^{t}\| \rightarrow 0$ when $\|\x^{t+1} - \x^{t}\| \rightarrow 0$. \bfit{(iii)} $\|\r^t\|$ is bounded, and $L_F' \rightarrow L_F$ as $t_0 \rightarrow \infty$. \bfit{(iv)} We observe $\mu^{t-1}\rightarrow0$ and $\tfrac{\mu^{t-1}}{\mu^t}\rightarrow 1$ as $t\rightarrow \infty$, resulting in $\|\varepsilons^t\|\rightarrow0$. \bfit{(v)} The inequalities in Lemma \ref{lemma:first:bound} are independent of the choice of strategies for solving the $\x$-subproblem, and they hold deterministically.

The following lemma is useful in our subsequent analysis.

\begin{lemma} \label{lemma:bound:negative:l2}
\noi (Proof in Appendix \ref{app:lemma:bound:negative:l2}) Let $\bar{\x}$ be any \textit{global optimal solution} of Problem (\ref{eq:main}). We have:

\begin{enumerate}[label=\textbf{(\alph*)}, leftmargin=22pt, itemsep=1pt, topsep=1pt, parsep=0pt, partopsep=0pt]

\item It holds that: $\la \r^t, \bar{\x} - \x^t \ra\leq \Upsilon^t  -\frac{V_s}{2}\|\x^t -\bar{\x}\|_2^2$, where $\Upsilon^t\triangleq  F(\bar{\x}) - F(\x^t) + \tfrac{\mu^{t-1} L_h^2 }{2}  +  \tfrac{2 L_F L_h}{V_s} \| \A\|(\tfrac{\mu^{t-1}}{\mu^{t}} - 1)$.

\item If $\mu^t = \bar{\mu}$, we have: $\Upsilon^t \leq \tfrac{1}{2} \bar{\mu}L_h^2 + F(\bar{\x}) - [\min_{i=1}^t F(\x^i)]$.

\item If $\mu^t = \frac{\eta}{t+t_0}$, we have: $\sum_{i=1}^t \Upsilon^i \leq  C_{\Upsilon} (\ln(t) + 1)  - t  [\min_{i=1}^t F(\x^i) ]  +  t  F(\bar{\x})$, where $C_{\Upsilon}\triangleq \frac{\eta L_h^2 }{2}  +  \tfrac{2 L_F L_h}{V_s}  \| \A\|$.

\end{enumerate}

\end{lemma}

\noi \textbf{Remarks}. Noticing that $F(\bar{\x}) - F(\x^t)\leq 0$, we have $\Upsilon^t\rightarrow 0$ as $t\rightarrow \infty$, and it holds that $\la \r^t,\x^t-\bar{\x} \ra\leq -\tfrac{V_s}{2}\|\x^t -\bar{\x}\|_2^2$ in the limit.

\subsection{Convergence Rate for SPGM-IHT} \label{sect:spgm:iht}

In this subsection, we assume that IHT strategy is used for solving the $\x$-subproblem.

We denote any limit point of SPGM-IHT as $(\dot{\x},\dot{\y})$ and present the following useful definition.

\begin{definition}
(Approximate Lipschitz Stationary Point) Given any constant $\epsilon>0$. Fix $\mu>0$ to be a sufficiently small constant. A solution $(\dot{\x},\dot{\y})$ is a $\epsilon$-approximate Lipschitz stationary point if: $\dist^2(\dot{\y}, \arg \min_{\y}\mathcal{J}(\dot{\x},\y;\mu))+\dist^2(\dot{\x},\arg \min_{\x}\delta(\x)+\dot{\mathcal{M}}(\x,\dot{\x},\dot{\y};\mu))\leq \epsilon$, where $\dot{\mathcal{M}}(\cdot,\cdot,\cdot;\cdot)$ is defined in Equation (\ref{eq:subprob:x:1}).

\end{definition}

The following theorem establishes the convergence of \textbf{SPGM-IHT}.

\begin{theorem} \label{theorem:convergence:iht}
\noi (Proof in Appendix \ref{app:theorem:convergence:iht}) \bfit{Convergence to Lipschitz Stationary Solutions}. We let $\mathcal{J}^{t+1} \triangleq\mathcal{J}(\x^{t+1},\y^{t+1};\mu^{t})$. We define $\Psi^t$ as in Lemma \ref{lemma:first:bound}. We have:

\begin{enumerate}[label=\textbf{(\alph*)}, leftmargin=22pt, itemsep=1pt, topsep=1pt, parsep=0pt, partopsep=0pt]

\item $\tfrac{\theta}{2}\|\x^{t+1}-\x^{t}\|_2^2  + \tfrac{1}{2\mu^0}\|\y^{t+1}-\y^t\|_2^2 \leq  \Psi^t + \mathcal{J}^t  -\mathcal{J}^{t+1}$.

\item $\sum_{t=1}^T [\tfrac{\theta}{2}\|\x^{t+1}-\x^{t}\|_2^2 + \tfrac{1}{2\mu^0}\|\y^{t+1}-\y^{t}\|_2^2] \leq \mathcal{J}^1 - \mathcal{J}^{T+1} +  \eta L_h^2 \triangleq C <\infty$.

\item Algorithm \ref{algo:main} finds an $\epsilon$-approximate Lipschitz stationary point of Problem (\ref{eq:main}) in at most $T$ iterations, where $T\leq \lceil  \frac{2C}{\epsilon \min(\theta,(\mu^0)^{-1}) }\rceil = \mathcal{O}(\epsilon^{-1})$.

\end{enumerate}
\end{theorem}

\noi \textbf{Remarks}. The introduction of parameter $\theta>0$ is important since it guarantees sufficient decrease condition and global convergence of Algorithm \ref{algo:main}.

In what follows, we present enhanced convergence results for \textbf{SBCD-IHT}, leading to the attainment of the global optimal solution $\bar{\x}$. We use the following quantities to measure the distance between $\x^t$ and $\bar{\x}$:
\beq
\Delta_{\x}^t \triangleq \|\x^{t} - \bar{\x}\|_2^2,\,\text{and}\,\Delta_{F}^t \triangleq [\min_{i=1}^t F(\x^i) ]  - F(\bar{\x}) .\nn
\eeq

We first have the following useful lemma.

\begin{lemma} \label{lemma:iht:bound}
\noi (Proof in Appendix \ref{app:lemma:iht:bound}) We define $H^t \triangleq  A_s^2 /\mu^t + M_s+\theta$. We have:

\begin{enumerate}[label=\textbf{(\alph*)}, leftmargin=22pt, itemsep=1pt, topsep=1pt, parsep=0pt, partopsep=0pt]

\item $\|\x^{t+1} - \x^t\| \leq \tfrac{2}{ H^t}\|\r^t\|$.

\item $\la \x^t_+, \x^t_+\ra \geq \la \x^t_+,\x^{t+1} \ra$, where $\x^t_+ \triangleq \x^t - \r^t /H^t$.

\item $\tfrac{1}{2}H^t\Delta_{\x}^{t+1} \leq \tfrac{1}{2}(H^t-V_s)\Delta_{\x}^t + \tfrac{3\|\r^{t}\|_2^2}{H^t} + \Upsilon^t + \|\r^{t}\| \|\bar{\x} \|$.

\end{enumerate}

\end{lemma}

\noi \textbf{Remarks}. As $\mu^t\rightarrow 0$, we have $H^t \rightarrow + \infty$, leading to $\|\x^{t+1} - \x^t\| \rightarrow  0$.

The following theorems establish the convergence
of \textbf{SPGM-IHT} to the global optimal solution $\bar{\x}$.

\begin{theorem}\label{theorem:IHT:stepsize:constant}
\noi (Proof in Appendix \ref{app:theorem:IHT:stepsize:constant}) \bfit{Convergence to the Global Optimal Solutions for Constant Stepsizes}. Assume constant stepsizes are used with $\mu^{t} = \bar{\mu}$. We define $\gamma \triangleq 1-{V_s}/{H}$, and $H \triangleq  A_s^2  /\bar{\mu} + M_s+\theta$. We have the following recursive inequality: $\tfrac{1}{2}\Delta_{\x}^{t+1}  \leq   \tfrac{\gamma}{2}\Delta_{\x}^{t} + \frac{L_F}{H} \|\bar{\x}\|+  \tfrac{3L_F^2}{(H)^2} + (\tfrac{1}{2}\bar{\mu}L_h^2 - \Delta_{F}^t)/H$. Furthermore, it holds that:
\beqq
\Delta_{F}^t   \leq    K_1 \gamma^t   +  D_1  \bar{\mu}   +    L_F \|\bar{\x}\|, \label{eq:IHT:constant:F}
\eeqq
\beqq
\Delta_{\x}^{t+1}  \leq \left( K_1 \gamma^t  +  D_1 \bar{\mu} +   L_F \|\bar{\x}\| \right) \tfrac{2}{V_s}, \label{eq:IHT:constant:x}
\eeqq
\noi where $K_1\triangleq \frac{V_s}{2}\Delta_{\x}^1$, and $D_1 \triangleq   \tfrac{3L_F^2 }{A_s^2} +  \tfrac{1}{2}L_h^2$.

\end{theorem}

\begin{theorem}\label{theorem:IHT:stepsize:diminishing}
\noi (Proof in Appendix \ref{app:theorem:IHT:stepsize:diminishing}) \bfit{Convergence to Global Optimal Solutions for Diminishing Stepsizes}. Assume diminishing stepsizes are used with $\mu^t = \tfrac{\eta}{t+t_0}$, where $\eta={A_s^2}/{V_s}$. We let $L'_F \triangleq L_f + \frac{t_0+1}{t_0}L_h \|\A\|$, and $H^t \triangleq  A_s^2 /\mu^t + M_s+\theta$. We define $\Upsilon^t$ as in Lemma \ref{lemma:bound:negative:l2}. We have the following recursive inequality: $\tfrac{1}{2}(H^{t+1}  - V_s) \Delta_{\x}^{t+1}  \leq  \tfrac{1}{2}(H^t - V_s)\Delta_{\x}^t +
{3   (L_F')^2  } / (V_s \cdot t) +  \Upsilon^t  +   L_F' \|\bar{\x} \|$. Furthermore, it holds that:
\beqq
\Delta_{F}^t  \leq   \tfrac{K_2 }{t}     +  \tfrac{D_2 ( \ln(t)+1) }{t }      +   L_F' \|\bar{\x} \|    , \label{eq:IHT:diminishing:F}
\eeqq
\beqq
\Delta_{\x}^{t+1} \leq (  \tfrac{K_2}{t+1} +   \tfrac{D_2 (\ln(t)+1) }{t+1}   +    L_F'  \|\bar{\x} \|  )\tfrac{2}{V_s}, \label{eq:IHT:diminishing:x}
\eeqq
\noi where $K_2 \triangleq \tfrac{H^1}{2}\Delta_{\x}^1$, and $D_2 \triangleq \tfrac{3 (L_F')^2}{V_s} +  C_{\Upsilon}$, with $C_{\Upsilon}$ defined in Lemma \ref{lemma:bound:negative:l2}.

\end{theorem}

\noi \textbf{Remarks}. \bfit{(i)} As $t_0\rightarrow +\infty$, we have $L_F'=L_f + \frac{t_0+1}{t_0} L_h\|\A\|\rightarrow  L_F$. \bfit{(ii)} The irreducible estimation error terms $L_F\|\bar{\x}\|$ and $L'_F \|\bar{\x}\|$, as specified in (\ref{eq:IHT:constant:F}) and in (\ref{eq:IHT:diminishing:F}) respectively, match the best-known error bounds for this nonconvex NP-hard problem detailed in \cite{liu2019one}, specifically in \textbf{Corollary III.4} and \textbf{Corollary III.8}. \bfit{(iii)} Given $\|\bar{\x}\|\leq s$, we obtain: $\|\bar{\x}\| \leq \sqrt{s} \|\bar{\x}\|_{\infty}$, rendering the irreducible error terms small. Hence, our theoretical bounds can exploit the inherent sparsity structure of the problem. \bfit{(iv)} The irreducible error terms in Inequalities (\ref{eq:IHT:constant:F}), (\ref{eq:IHT:constant:x}), (\ref{eq:IHT:diminishing:F}), and (\ref{eq:IHT:diminishing:x}) depend on $\|\bar{\x}\|$, the Lipschitz constant $L_F$ of $F(\x)$, and the strong convexity parameter $V_s$ of $f(\x)$, indicating the difficulty of solving this NP-hard problem. \textbf{SPGM-IHT} is more likely to converge to the global optimum when $\|\bar{\x}\|$ and $L_F$ are small while $V_s$ is large.

\subsection{Convergence Rate for \textbf{SPGM-BCD}}\label{sect:spgm:bcd}

In this subsection, we assume that BCD strategy is used for solving the $\x$-subproblem.

We assume that the working set $\B$ is selected randomly and uniformly from $\Omega_n^k \triangleq \{\mathcal{\B}_{1},\mathcal{\B}_{2},...,\mathcal{\B}_{C_n^k}\}$. \textbf{SPGM-BCD} generates
a random output $\x^t$ with $t = 1,2, ...$, based on the observed realization of the random variable $\xi^{t-1} \triangleq \{\B^1,\B^2,..., \B^{t-1}\}$. The expectation of a random variable is denoted by $\E_{\xi^t}[\cdot ]$. The following lemma is useful in this context.

\begin{lemma} \label{lemma:bcd:simple:basic}
\noi (Proof in Appendix \ref{app:lemma:bcd:simple:basic}) For any $\x\in\mathbb{R}^n$ and $\z\in\mathbb{R}^n$, we have $\frac{1}{ C_n^k}\sum_{\B \in \Omega_n^k}\,\x\trans (\UB \UB\trans ) \z  = Z_n^k \la \x,\z\ra$, and $\E_{\B} \|\x_{\B}\|_2^2 = Z_n^k \|\x\|_2^2$, where $Z_n^k \triangleq \frac{k}{n}$.

\end{lemma}

We denote any limit point of \textbf{SPGM-BCD} as $(\ddot{\x},\ddot{\y})$ and offer the following useful definition.

\begin{definition}
(Approximate block-$k$ Stationary Point) Given any constant $\epsilon>0$. Fix $\mu>0$ to be a sufficiently small constant. A solution $(\ddot{\x},\ddot{\y})$ is an $\epsilon$-approximate block-$k$ stationary point if: $\dist^2(\ddot{\y}, \arg \min_{\y}\mathcal{J}(\ddot{\x},\y;\mu))+\tfrac{1}{C_n^k}\sum_{\B\in \Omega_n^k} \dist^2(\ddot{\x}_{\B}, \arg \min_{\x_{\B}}\delta(\UB\x_{\B}+\UBc \ddot{\x}_{\Bc})+\ddot{\mathcal{M}}(\UB\x_{\B}+\UBc \ddot{\x}_{\Bc},\ddot{\x},\ddot{\y};\mu))\leq\epsilon$, where $\ddot{\mathcal{M}}(\cdot,\cdot,\cdot;\cdot)$ is defined in Equation (\ref{eq:subprob:x:2}).

\end{definition}

The following theorem establishes the convergence of \textbf{SPGM-BCD}.

\begin{theorem} \label{theorem:convergence:bcd}

\noi (Proof in Appendix \ref{app:theorem:convergence:bcd}) \bfit{Convergence to Block-$k$ Stationary Solutions}. We let $\theta \triangleq \frac{\theta_1}{\mu^1}+\theta_2$, and $\mathcal{J}^{t+1} \triangleq \E_{\xi^{t}}[\mathcal{J}(\x^{t+1},\y^{t+1};\mu^{t})]$. We define $\Psi$ as in Lemma \ref{lemma:first:bound}. We have:

\begin{enumerate}[label=\textbf{(\alph*)}, leftmargin=22pt, itemsep=1pt, topsep=1pt, parsep=0pt, partopsep=0pt]

\item $\E_{\xi^t}[\tfrac{\theta}{2}\|\x^{t+1}-\x^{t}\|_2^2+\tfrac{1}{2\mu^1}\|\y^{t+1}-\y^t\|_2^2 ]   \leq  \Psi^t + \mathcal{J}^t  -\mathcal{J}^{t+1}$.

\item $\E_{\xi^T}[\sum_{t=1}^T \tfrac{\theta}{2}\|\x^{t+1}-\x^{t}\|_2^2 +  \tfrac{1}{2\mu^1}\|\y^{t+1}-\y^{t}\|_2^2 ] \leq \mathcal{J}^1 - \mathcal{J}^{T+1} +  \eta L_h^2 \triangleq C <\infty$.

\item Algorithm \ref{algo:main} finds an $\epsilon$-approximate block-$k$ stationary point of Problem (\ref{eq:main}) in at most $T$ iterations in the sense of expectation, where $T\leq \lceil  \frac{2C}{\epsilon \min(\theta,(\mu^0)^{-1}) }\rceil = \mathcal{O}(\epsilon^{-1})$.

\end{enumerate}

\end{theorem}

\noi \textbf{Remarks}. Theorem \ref{theorem:convergence:bcd} resembles Theorem \ref{theorem:convergence:iht}, with the key distinction being that \textbf{SPGM-IHT} deterministically converges to a Lipschitz stationary point, whereas \textbf{SPGM-BCD} converges to a block-$k$ stationary point in expectation.

In what follows, we present enhanced convergence results for SBCD-BCD, leading to the attainment of the global optimal solution $\bar{\x}$. For notation convenience, we define
\begin{align}
&\Delta_{\x}^t \triangleq \E_{\xi^{t}} [\|\x^t - \bar{\x}\|_2^2],\,\Delta_{F}^t \triangleq \E_{\xi^{t}} [(\min_{i=1}^t F(\x^i) )  -   F(\bar{\x})], \nn\\
&\Vup\triangleq \max_{\B\in\Omega_n^k} \lambdas_{\max}(\tilde{\M}_{\B\B}),\,\Vlow\triangleq \min_{\B\in\Omega_n^k}\lambdas_{\min}(\tilde{\M}_{\B\B}),  \nn\\
&\Aup\triangleq\max_{\B\in\Omega_n^k}\lambdas_{\max}(\A\trans\A)_{\B\B}),\,\Alow\triangleq \min_{\B\in\Omega_n^k}\lambdas_{\min}((\A\trans\A)_{\B\B}),  \nn\\
&\Sup^t \triangleq \tfrac{\Aup+\theta_1}{\mu^t} + \Vup+\theta_2,\,\Slow^t \triangleq \tfrac{\Alow+\theta_1}{\mu^t} + \Vlow+\theta_2, \kappa^t\triangleq \tfrac{\Sup^t}{\Slow^t},\label{eq:def:H}\nn
\end{align}
\noi where $\Vlow\geq V_s$, and $\Alow$ can be zero.

We first have the following two useful lemmas.

\begin{lemma} \label{lemma:condition:constant}

(Proof in Appendix \ref{app:lemma:condition:constant})  Given any constant $\epsilon>0$. If $\theta_1$ and $\theta_2$ are sufficiently large such that $\theta_1 \geq \mathcal{T}_1(\epsilon)\triangleq \frac{\Aup - \Alow (1+\epsilon)}{\epsilon}$ and $\theta_2 \geq \mathcal{T}_2(\epsilon)\triangleq \frac{\Vup - \Vlow (1+\epsilon)}{\epsilon}$, we have: $\kappa^t \triangleq \frac{\Sup^t}{\Slow^t} \leq 1 + \epsilon$.

\end{lemma}

\begin{lemma} \label{lemma:bcd:bound}

\noi (Proof in Appendix \ref{app:lemma:bcd:bound}) We let $\H^t\triangleq (\A\trans \A + \theta_1 \I_n)/ {\mu^t} + \tilde{\M} + \theta_2 \I_n$. We define $Z_n^k$ as in Lemma \ref{lemma:bcd:simple:basic}. For all $t\geq 1$, we have:

\begin{enumerate}[label=\textbf{(\alph*)}, leftmargin=22pt, itemsep=1pt, topsep=1pt, parsep=0pt, partopsep=0pt]

\item $\E_{\xi^t}[\|\x^{t+1} - \x^{t}\|] \leq  \frac{2}{\Slow^t} \E_{\xi^t}[\|\r^t\|]$.

\item $\E_{\xi^t}[ \la [[\H^t]_{\B\B}]( \x_{\B}^{t+1} - \x_{\B}^{t}), \x_{\B}^{t+1} \ra ] = -Z_n^k  \la \r^t, \x^{t+1}\ra$.

\item $\E_{\xi^t}[\frac{1}{2}\|\x^{t+1}-\bar{\x}\|_{\H_{\ast}^t}^2-\frac{1}{2}\| \x^{t} - \bar{\x}\|_{\H_{\ast}^t}^2]  \leq Z_n^k  [ 2/\Slow^t  \|\r^t\|_2^2  +  \Upsilon^t + ( 1+2\kappa^t)\|\r^t\| \|\bar{\x}\| - \frac{V_s}{2}\| \x^{t} - \bar{\x} \|_2^2]$, where $\H_{\ast}^t \triangleq \UBt \UBt\trans \H^t \UBt  \UBt\trans$.

\end{enumerate}

 \end{lemma}

\noi \textbf{Remarks}. When $\mu^t \rightarrow 0$, we have $\Slow^t \rightarrow +\infty$, leading to $\E_{\xi^t}[\|\x^{t+1} - \x^t\|] \rightarrow  0$.


The following theorems establish the convergence
of \textbf{SPGM-BCD} to the global optimal solution $\bar{\x}$.

\begin{theorem} \label{theorem:DEC:stepsize:constant}

(Proof in Appendix \ref{app:theorem:DEC:stepsize:constant}) \bfit{Convergence to Global Optimal Solutions for Constant Stepsizes}. Assume constant stepsizes are used with $\mu^{t} = \bar{\mu}$. Given any constant $\epsilon>0$. Assume that $\theta_1\geq \mathcal{T}_1(\epsilon)$ and $\theta_2\geq \mathcal{T}_2(\epsilon)$, where $\mathcal{T}_1(\cdot)$ and $\mathcal{T}_2(\cdot)$ are define in Lemma \ref{lemma:condition:constant}. We let $\Sup \triangleq \frac{\Aup+\theta_1}{\bar{\mu}} + \Vup+\theta_2,\,\Slow \triangleq \frac{\Alow+\theta_1}{\bar{\mu}} + \Vlow+\theta_2$, and $\gamma \triangleq 1 - \frac{V_s}{\Sup} \in (0,1)$. We define $\H^t_{\ast}$ as in Lemma \ref{lemma:bcd:bound}, and $\Upsilon^t$ as in Lemma \ref{lemma:bound:negative:l2}. We have the following recursive inequality: $\E_{\xi^{t+1}}[\tfrac{1}{2}\|\x^{t+1} -  \bar{\x}\|_{\H_{\ast}^{t+1}}^2] \leq \gamma  \E_{\xi^{t}}[\tfrac{1}{2}\|\x^t -  \bar{\x}\|_{\H_{\ast}^t}^2]+Z_n^k [\Upsilon^t + \frac{2(L_F)^2 }{ \Slow }+(3+2\epsilon) L_F \|\bar{\x}\|]$. Furthermore, it holds that:
\beqq
\Delta_{F}^{t} \leq     K_3 \gamma^t  + D_3 \bar{\mu} + (3 + 2 \epsilon) L_F \|\bar{\x}\|    , \label{eq:BCD:constant:F}
\eeqq
\beqq
\Delta_{\x}^{t+1}  \leq ( K_3 \gamma^t +  D_3\bar{\mu}    + ( 3 +2 \epsilon) L_F \|\bar{\x}\| )  \tfrac{2}{V_s} , \label{eq:BCD:constant:x}
\eeqq
\noi where $K_3 \triangleq \tfrac{V_s (1+\epsilon) }{2} \Delta_{\x}^1$, and $D_3\triangleq \frac{2(L_F)^2 }{  \theta_1 + \Aup  } + \frac{L_h^2}{2}$.

\end{theorem}

\begin{theorem} \label{theorem:DEC:stepsize:diminishing}

(Proof in Appendix \ref{app:theorem:DEC:stepsize:diminishing}) \bfit{Convergence to Global Optimal Solutions for Diminishing Stepsizes}. Assume diminishing stepsizes are used with $\mu^t = \tfrac{\eta}{t+t_0}$, where $\eta = \frac{\Aup + \theta_1}{V_s}$. Given any constant $\epsilon>0$. Assume that $\theta_1\geq \mathcal{T}_1(\epsilon)$ and $\theta_2\geq \mathcal{T}_2(\epsilon)$, where $\mathcal{T}_1(\cdot)$ and $\mathcal{T}_2(\cdot)$ are define in Lemma \ref{lemma:condition:constant}. We let $\Sup^t \triangleq \frac{\Aup+\theta_1}{\mu^t} + \Vup+\theta_2$, $\Slow^t\triangleq \frac{\Alow+\theta_1}{\mu^t} + \Vlow+\theta_2$, and $L'_F \triangleq L_f + \frac{t_0+1}{t_0}L_h \|\A\|$. We define $\H^t_{\ast}$ as in Lemma \ref{lemma:bcd:bound}, and $\Upsilon^t$ as in Lemma \ref{lemma:bound:negative:l2}. We have the following recursive inequality: $\Phi^{t+1} - \Phi^t  \leq   Z_n^k  \Upsilon^t +  \frac{2  Z_n^k  (L'_F)^2 }{ \Slow^t }     +  ( 3+2\epsilon )     Z_n^k  L'_F \|\bar{\x}\|$, where $\Phi^t \triangleq \E_{\xi^t} [\frac{1}{2}\| \x^{t} - \bar{\x}\|_{\H_{\ast}^t}^2] - Z_n^k  \frac{V_s}{2}\|\x^t-\bar{\x}\|_2^2$. Furthermore, it holds that:
\beqq
\Delta_{F}^t \leq  \tfrac{K_4}{t} +  \tfrac{ D_4 (1 + \ln(t) ) }{t}      + ( 3 +    2\epsilon )   L_F' \|\bar{\x}\|,  \label{eq:BCD:diminishing:F}
\eeqq
\beqq
\Delta_{\x}^{t+1} \leq(\tfrac{K_4}{t+1} +  \tfrac{ D_4 (\ln(t)+1 ) }{t+1}+( 3 +    2\epsilon )   L_F' \|\bar{\x}\|)\tfrac{2}{V_s},  \label{eq:BCD:diminishing:x}
\eeqq
\noi where $K_4\triangleq \frac{\Sup^1}{2} \Delta_{\x}^{1}$, and $D_4 \triangleq    \frac{2 (L_F')^2 }{V_s} + C_{\Upsilon}$, with $C_{\Upsilon}$ defined in Lemma \ref{lemma:bound:negative:l2}.

\end{theorem}

\noi \textbf{Remarks}. \bfit{(i)} The convergence rates in Theorems \ref{theorem:DEC:stepsize:constant} and \ref{theorem:IHT:stepsize:constant} are similar, as are those in Theorems \ref{theorem:DEC:stepsize:diminishing} and \ref{theorem:IHT:stepsize:diminishing}. However, analyzing \textbf{SPGM-BCD} is more intricate than \textbf{SPGM-IHT} due to the utilization of a general Hessian matrix $\H^t$ and a stochastic mechanism of \textbf{SPGM-BCD}, in contrast to the utilization of a scaled identity matrix $H^t \I_n$ and a deterministic mechanism of \textbf{SPGM-IHT}. Consequently, their strategies differ significantly. \bfit{(ii)} As $\epsilon \rightarrow 0$ and $t_0\rightarrow+\infty$, the irreducible estimation error terms for $\Delta_{F}^t$ in (\ref{eq:BCD:constant:F}) and (\ref{eq:BCD:diminishing:F}) simplify to $3 L_F\|\bar{\x}\|$, which is three times the bound of \textbf{PSGD} in \cite{liu2019one}. Our bounds leverage the inherent sparsity of the problem.

\section{Experiments}
\label{sect:exp}

This section evaluates the effectiveness of \textbf{SPGM-IHT} and \textbf{SPGM-BCD}, comparing them with five state-of-the-art nonsmooth sparsity constrained optimization algorithms: \bfit{(i)} Projective Subgradient Descent (PSGD) \cite{liu2019one}. \bfit{(ii)} Alternating Direction Method of Mutipliers based on IHT (ADMM-IHT) \cite{HeY12}. \bfit{(iii)} Dual Iterative Hard Thresholding(DIHT) \cite{yuan2020jmlr}. \bfit{(iv)} Convex $\ell_1$ Approximation Method (CVX-$\ell_1$) \cite{candes2005decoding}. \bfit{(v)} Nonconvex $\ell_p$ Approximation Method (NCVX-$\ell_p$) \cite{xu2012l}.

Our experiments reveal that \textbf{SPGM-IHT} is on par with existing IHT-style methods, and \textbf{SPGM-BCD} consistently delivers the best performance. This outcome is expected as \textbf{SPGM-IHT} is an IHT-style method itself, while \textbf{SPGM-BCD} excels in identifying stronger stationary points compared to other approaches. Due to space constraints, detailed experiment results are provided in the \textbf{Appendix}.

\section{Conclusions}

This paper explores Smoothing Proximal Gradient Method (\textbf{SPGM}) for solving nonsmooth sparsity constrained optimization problems. We discuss two specific variants of SPGM: one based on Iterative Hard Thresholding (\textbf{SPGM-IHT}) and the other on Block Coordinate Decomposition (\textbf{SPGM-BCD}). We provide both smooth and optimality analyses for the smoothing functions, demonstrating that \textbf{SPGM-BCD} discovers stronger stationary points of the nonsmooth nonconvex problem. We offer theoretical insights into the convergence rates of the \textbf{SPGM-IHT} and \textbf{SPGM-BCD} algorithms. Our bounds depend on the Lipschitz constant of the objective function, the strong convexity parameter of its smooth component, and the $\ell_2$ norm of the global optimal point. Leveraging the inherent sparsity of the optimization problem, our bounds align with the most competitive error estimates in the field. Finally, numerical experiments demonstrate that \textbf{SPGM-IHT} performs on par with existing IHT-style methods, while \textbf{SPGM-BCD} consistently delivers state-of-the-art numerical performance.


\normalem
\bibliographystyle{icml2024}
\bibliography{my}

\clearpage
\onecolumn
\appendix
{\huge \textbf{Appendix}}

The appendix is organized as follows.

Appendix \ref{app:sect:useful} contains some useful lemmas.

Appendix \ref{app:sect:smooth:analysis} includes the proofs for Section \ref{sect:smooth:analysis}.

Appendix \ref{app:sect:rate} presents the proofs for Section \ref{sect:rate}.

Appendix \ref{app:sect:exp} provides the experimental results.

\section{Some Useful Lemmas}  \label{app:sect:useful}


We present some useful lemmas that will be used subsequently.


\begin{lemma}\label{lemma:relation}
(\textbf{Pythagoras Relation}) For any symmetric matrix $\H \in \mathbb{R}^{n\times n}$ with $\H\trans = \H$ and any vectors $\a\in\mathbb{R}^n$, $\b\in \mathbb{R}^n$, $\c\in\mathbb{R}^n$, we have:
\beq
\frac{1}{2}\|\a-\b\|_{\H}^2  - \frac{1}{2} \|\c-\b\|_{\H}^2 =  \frac{1}{2}\|\a-\c\|_{\H}^2 - \la \a-\c,{\H}(\b-\c)\ra.\nn
\eeq

\end{lemma}

\begin{lemma} \label{lemma:two:non:negative:sequences}
Assume $\gamma\in(0,1)$. Denote $\gamma^t$ as the $t$-th power of $\gamma$. Let $\{\Phi^t\}_{t=1}^{\infty}$ and $\{\Lambda^t\}_{t=1}^{\infty}$ be any two non-negative sequences. We have:
\beq
\left(\Phi^{t+1} \leq \gamma \Phi^{t} + \Lambda^t\right) \Rightarrow\left( \Phi^{t+1}  \leq \gamma^t \Phi^{1} +  \tfrac{ \max_{i=1}^t(\Lambda^i)  }{1-\gamma} \right)\nn.
\eeq

\begin{proof}
    Using basic induction, we have the following results:
    \beq
    &&t=1,\,\Phi^{2} \leq \gamma \Phi^{1} + \Lambda^1 \nn\\
    &&t=2,\,\Phi^{3} \leq \gamma \Phi^{2} + \Lambda^2 \leq  \gamma (\gamma \Phi^{1} + \Lambda^2) + \Lambda^1 = \gamma^2 \Phi^{1} +   (\Lambda^2 + \gamma \Lambda^1)     \nn\\
    &&t=3,\,\Phi^{4} \leq \gamma \Phi^{3} +\Lambda^3 \leq  \gamma (\gamma^2 \Phi^{1} +   (\Lambda^2 + \gamma \Lambda^1)) + \Lambda^3 = \gamma^3 \Phi^{1} +   (\Lambda^3+ \gamma \Lambda^2  + \gamma^2 \Lambda^1)     \nn\\
    && ... \nn
    \eeq
\noi Therefore, we obtain:
\beq
\Phi^{T+1} \leq \gamma^T \Phi^{1} +  \sum_{i=1}^T \Lambda^i\gamma^{T-i}  \overset{\step{172}}{\leq}  \gamma^T \Phi^{1} + (\max_{i=1}^T \Lambda^i) \cdot (\sum_{i=1}^T \gamma^{T-i}) \overset{\step{173}}{\leq} \gamma^T \Phi^{1} +  \frac{(\max_{i=1}^T \Lambda^i)}{ 1-\gamma},\nn
\eeq
\noi where step \step{172} uses the Cauchy-Schwarz Inequality; step \step{173} uses the fact that:
\beq
\sum_{i=1}^t\gamma^{t-i}= 1 + \gamma^1 + \gamma^2 + ... +  \gamma^{t-1} = \frac{1-\gamma^t}{1-\gamma} < \frac{1}{1-\gamma}.\nn
\eeq

\end{proof}

\end{lemma}

\section{Proofs for Section \ref{sect:smooth:analysis}}
\label{app:sect:smooth:analysis}

\subsection{Proof of Lemma \ref{lemma:lip:variable:mu}} \label{app:lemma:lip:variable:mu}

\begin{proof} Without loss of generality, we assume $\mu_1 < \mu_2$. For all $\x$ with $\|\x\|_0\leq s$, we define:
\beq \label{eq:GG}
\psi(\mu) \triangleq \G(\x;\mu)   = f(\x)+ h(\P_{\mu}(\c)) + \tfrac{1}{2\mu} \|\c - \P_{\mu}(\c) \|_2^2 \,\text{with}\,\c \triangleq \A\x-\b.
\eeq
\noi Using the definition of $\P_\mu(\c)$ as shown in (\ref{eq:proximal:operator}), we have for any given $\mu_1$ and $\mu_2$:
\beq
\P_{\mu_1}(\c) = \arg \min_{\y}\,h(\y) + \tfrac{1}{2\mu_1} \|\y - \c\|_2^2 ,\,\text{and}\, \P_{\mu_2}(\c) = \arg \min_{\y}\,h(\y) + \tfrac{1}{2\mu_2} \|\y - \c\|_2^2 .\nn
\eeq
\noi By the optimality of $\P_{\mu_1}(\c)$ and $\P_{\mu_2}(\c)$, we obtain:
\beq
\c - \P_{\mu_1}(\c)  \in \mu_1 \partial h(\P_{\mu_1}(\c))  ,\,\text{and}\, \c - \P_{\mu_2}(\c)  \in \mu_2  \partial h(\P_{\mu_2}(\c)) .\label{eq:optimality:mu12}
\eeq
\noi \textbf{(a)} We now prove that $\psi(\mu)$ is a decreasing function. For any $\p_1 \in \partial h(\P_{\mu_1}(\c))$ and $\p_2\in \partial h(\P_{\mu_2}(\c))$, we have:
\beq
\psi(\mu_2) - \psi(\mu_1)&\overset{\step{172}}{=}& h(\P_{\mu_2}(\c)) + \tfrac{1}{2\mu_2} \|\c -  \P_{\mu_2}(\c) \|_2^2  - h(\P_{\mu_1}(\c)) - \tfrac{1}{2\mu_1} \|\c - \P_{\mu_1}(\c) \|_2^2 \nn\\
&\overset{\step{173}}{\leq}& \la \P_{\mu_2}(\c)-\P_{\mu_1}(\c),\p_2\ra  + \tfrac{1}{2\mu_2}\| \c -  \P_{\mu_2}(\c) \|_2^2 - \tfrac{1}{2\mu_1}\| \c -  \P_{\mu_1}(\c)\|_2^2 \nn\\
&\overset{\step{174}}{=}& \la \mu_1 \p_1 -\mu_2 \p_2,\p_2\ra  + \tfrac{\mu_2}{2}\|  \p_2\|_2^2 - \tfrac{\mu_1}{2}\| \p_1\|_2^2 \nn\\
&=& \la \mu_1 \p_1,\p_2\ra  - \tfrac{\mu_2}{2}\|  \p_2\|_2^2 - \tfrac{\mu_1}{2}\| \p_1\|_2^2 \nn\\
&\overset{\step{175}}{\leq}& \la \mu_1 \p_1 ,\p_2\ra  - \tfrac{\mu_1}{2}\|  \p_2\|_2^2 - \tfrac{\mu_1}{2}\| \p_1\|_2^2 \nn\\
&\overset{}{=}&   -\tfrac{\mu_1}{2} \|\p_1 - \p_2\|_2^2 \leq 0, \nn
\eeq
\noi where step \step{172} uses the definition of $\psi(\mu)$ in (\ref{eq:GG}); step \step{173} uses the convexity of $h(\cdot)$;  step \step{174} uses the optimality of $\P_{\mu_1}(\c)$ and $\P_{\mu_2}(\c)$ in (\ref{eq:optimality:mu12}); step \step{175} uses $\mu_1 < \mu_2$.

\noi \textbf{(b)} We now prove that $\psi(\mu)$ is $(\tfrac{1}{2}L_h^2)$-Lipschitz. For any $\p_1 \in \partial h(\P_{\mu_1}(\c))$ and $\p_2\in \partial h(\P_{\mu_2}(\c))$, we have:
\beq
 \psi(\mu_1) - \psi(\mu_2)&\overset{\step{172}}{=}&  h(\P_{\mu_1}(\c)) + \tfrac{1}{2\mu_1} \|\c -  \P_{\mu_1}(\c) \|_2^2-h(\P_{\mu_2}(\c)) - \tfrac{1}{2\mu_2} \|\c - \P_{\mu_2}(\c) \|_2^2   \nn\\
&\overset{\step{173}}{\leq}&  \la \P_{\mu_1}(\c) - \P_{\mu_2}(\c) ,  \p_1\ra + \tfrac{1}{2\mu_1} \|\c -  \P_{\mu_1}(\c) \|_2^2 - \tfrac{1}{2\mu_2} \|\c - \P_{\mu_2}(\c) \|_2^2   \nn\\
&\overset{\step{174}}{=}&  \la [\c - \P_{\mu_2}(\c) ] - [\c-\P_{\mu_1}(\c)], \p_1\ra + \tfrac{1}{2\mu_1} \|\c -  \P_{\mu_1}(\c) \|_2^2 - \tfrac{1}{2\mu_2} \|\c - \P_{\mu_2}(\c) \|_2^2  \nn\\
&\overset{\step{175}}{=}&  \la \mu_2 \p_2- \mu_1 \p_1,  \p_1\ra + \tfrac{\mu_1}{2} \| \p_1 \|_2^2 - \tfrac{\mu_2}{2} \|\p_2\|_2^2  \nn\\
&=&  - \tfrac{\mu_2}{2} \| \p_2\|_2^2  + \mu_2\la \p_1, \p_2\ra - \tfrac{\mu_1}{2} \|\p_1\|_2^2 \nn\\
&\overset{\step{176}}{\leq}&  \tfrac{\mu_2}{2}\|\p_1\|_2^2  - \tfrac{\mu_1}{2} \|\p_1\|_2^2   \nn\\
&\overset{\step{177}}{\leq}&   \tfrac{ \mu_2 - \mu_1   }{2} \cdot L_h^2,\nn
\eeq
\noi where step \step{172} uses the definition of $\psi(\mu)$ in (\ref{eq:GG}); step \step{173} uses the convexity of $h(\cdot)$;  step \step{174} uses the fact that $\P_{\mu_1}(\c) - \P_{\mu_2}(\c)= [\c - \P_{\mu_2}(\c)] - [\c-\P_{\mu_1}(\c)]$; step \step{175} uses the optimality of $\P_{\mu_1}(\c)$ and $\P_{\mu_2}(\c)$ in (\ref{eq:optimality:mu12}); step \step{176} uses the inequality that: $-\tfrac{\mu}{2} \| \p_2\|_2^2  +  \mu \la \p_1, \p_2\ra\leq \tfrac{\mu}{2} \|\p_1\|_2^2$ for all $\mu>0$ and for all $\p_1\in \mathbb{R}^m$ and $\p_2\in \mathbb{R}^m$; step \step{177} uses $\|\p_1\|\leq L_h$. Dividing both sides by $(\mu_2-\mu_1)$, we conclude that $\psi(\mu)$ is $(\tfrac{1}{2}L_h^2)$-Lipschitz.

\end{proof}

\subsection{Proof of Lemma \ref{lemma:lip:variable:x}}

\label{app:lemma:lip:variable:x}

\begin{proof} \noi We fix $\mu>0$ to be a constant. For any given $\x\in \mathbb{R}^n$ and $\x'\in\mathbb{R}^n$ with $\|\x\|_0\leq s$ and $\|\x'\|_0\leq s$, we define
\beq \label{eq:def:yy}
\c \triangleq  \A\x-\b,\,\text{and}\,\c' \triangleq \A\x'-\b.
\eeq
\noi Using the definition of $\P_\mu(\cdot)$ as shown in (\ref{eq:proximal:operator}), we have:
\beq
\P_{\mu}(\c)  = \arg\min_{\y} h(\y) + \tfrac{1}{2\mu} \|\y - \c\|_2^2,\,\text{and}\,\P_{\mu}(\c') = \arg\min_{\y} h(\y) + \tfrac{1}{2\mu} \|\y - \c' \|_2^2. \nn
\eeq
\noi By the optimality condition of $\P_{\mu}(\c)$ and $\P_{\mu}(\c')$, we have:
\beq
 \c - \P_{\mu}(\c)   \in  \mu  \partial h( \P_{\mu}(\c) ) ,\,\text{and}\,  \c' - \P_{\mu}(\c') \in  \mu \partial h( \P_{\mu}(\c) ) .\label{eq:opt:PcPc}
\eeq
\noi The function $\mathcal{G}(\x;\mu)$ defined in (\ref{eq:def:GGG}) is differentiable and its gradient at $\x$ and $\x'$ can be respectively computed as:
\beq
\nabla_{\x} \mathcal{G}(\x;\mu) = \nabla f(\x) +  \tfrac{1}{\mu} \A\trans (\c - \P_{\mu}(\c))  \,\text{and}\,\nabla_{\x} \mathcal{G}(\x';\mu) = \nabla f(\x') + \tfrac{1}{\mu} \A\trans (\c' - \P_{\mu}(\c')). \label{eq:grad:G:mu}
\eeq
\noi \textbf{(a)} We notice that $F(\x) = \lim_{\bar{\mu} \rightarrow 0} \mathcal{G}(\x;\bar{\mu}) $ and $\mathcal{G}(\x;\mu)$ is a decreasing function \textit{w.r.t.} $\mu$. The inequality $F(\x) \geq \mathcal{G}(\x;\mu)$ clearly holds. We now prove that $F(\x) - \tfrac{\mu  }{2}   L_h^2   \leq \mathcal{G}(\x;\mu)$. For any $\x$ with $\|\x\|_0\leq s$ and $\p\in \partial h( \P_{\mu}(\c) )$, we obtain:
\beq
 F(\x) -\mathcal{G}(\x;\mu) &{\overset{\step{172}}{=}}&  [f(\x) + h(\A\x-\b)] - [f(\x) + h(\P_{\mu}(\A\x - \b)) + \tfrac{1}{2\mu} \| \A\x - \b - \P_{\mu}(\A\x - \b)\|_2^2 ]\nn\\
&{\overset{\step{173}}{=}}&  h(\c) - h(\P_{\mu}(\c)) - \tfrac{1}{2\mu} \| \c - \P_{\mu}(\c) \|_2^2   \nn\\
&{\overset{\step{174}}{ \leq}}&    \la \c -  \P_{\mu}(\c) , \partial h(   \P_{\mu}(\c)   ) \ra  - \tfrac{1}{2\mu} \|\mu  \partial h( \P_{\mu}(\c) ) \|_2^2   \nn\\
&{\overset{}{\leq}}& \la \c -  \P_{\mu}(\c) ,\p \ra - \tfrac{\mu}{2} \| \p \|_2^2 \nn\\
&{\overset{\step{175}}{\leq}}& \tfrac{1}{2\mu} \|  \c -  \P_{\mu}(\c) \|_2^2\nn\\
&{\overset{\step{176}}{=}}& \tfrac{1}{2\mu} \|  \mu \p \|_2^2\nn\\
&{\overset{\step{177}}{\leq}}& \tfrac{\mu}{2} L_h^2 ,\nn
\eeq
\noi where step \step{172} uses the definition of $F(\x) \triangleq f(\x) + h(\A\x - \b)$ in (\ref{eq:main}) and the definition of $\mathcal{G}(\x;\mu)$ in (\ref{eq:def:GGG}); step \step{173} uses $\A\x-\b=\c$; step \step{174} uses the convexity of $h(\cdot)$ and the optimality of $\P_{\mu}(\c)$ as shown in (\ref{eq:opt:PcPc}); step \step{175} uses the inequality $-\tfrac{\mu}{2}\| \p\|_2^2 + \la\v, \p\ra\leq \tfrac{1}{2\mu} \|\v\|_2^2$ for all $\v$ and $\mu>0$; step \step{176} uses (\ref{eq:opt:PcPc}); step \step{177} uses $\|\p\|_2\leq L_h$.

\noi \textbf{(b)} We now prove that $F(\x)$ is $(L_f + L_h \|\A\|)$-Lipschitz. We have:
\beq
 \| \partial F(\x) \| &\overset{\step{172}}{=}& \|\nabla f(\x) + \A\trans \partial h(\A\x-\b)\|\nn\\
&\overset{\step{173}}{\leq}&  \|\nabla f(\x)\|+ \|\A\|\|\partial h(\A\x-\b)\| \nn\\
&\overset{\step{174}}{\leq}& L_f + L_h \|\A\|,\nn
\eeq
\noi where step \step{172} uses $\partial F(\x)= \nabla f(\x) + \A\trans \partial h(\A\x-\b)$; step \step{173} uses the norm inequality; step \step{174} uses the fact that $h(\cdot)$ is $L_h$-Lipschitz and $f(\cdot)$ is $L_f$-Lipschitz.

We now prove that $\mathcal{G}(\x,\mu)$ is $(L_f + L_h \|\A\|)$-Lipschitz. We obtain:
\beq
\|\nabla_{\x} \mathcal{G}(\x;\mu)\| &\overset{}{=}& \|\nabla f(\x) + \tfrac{1}{\mu} \A\trans (\c -  \P_{\mu}(\c))\| \nn\\
&\leq& \|\nabla f(\x)\| + \|\A\|\cdot \|\tfrac{1}{\mu}  (\c  -  \P_{\mu}(\c))\|\nn\\
&\overset{\step{172}}{=}& \|\nabla f(\x)\| + \|\A\|\cdot\| \partial h( \P_{\mu}(\c) ) \| \nn\\ &\leq& L_f +  L_h \|\A\|, \nn
\eeq
 \noi where step \step{172} uses the optimality condition of $\P_{\mu}(\c)$ as shown in (\ref{eq:opt:PcPc}) that $   \c - \P_{\mu}(\c) \in  \mu  \partial h( \P_{\mu}(\c) ) $.

\noi \textbf{(c)} Noticing $f(\x)$ is restricted $V_s$-strongly convex, we directly conclude that $\mathcal{G}(\x,\mu)$ is also restricted $V_s$-strongly convex. We now prove that the function $\mathcal{G}(\x,\mu)$ is restricted $(M_s+\frac{A_s \|\A\| }{\mu})$-smooth. For any $\x\in\mathbb{R}^n$, $\x'\in\mathbb{R}^n$, $\p_1 \in \partial h( \P_{\mu}(\c) )$, and $\p_2 \in \partial h( \P_{\mu}(\c') )$, we derive:
\beq \label{eq:lip:bound}
&&\| [\A\x' -  \A\x]  + [\P_{\mu}(\c)- \P_{\mu}(\c')]\|_2^2 \nn \\
&\overset{\step{172}}{=} & \| \A\x'  -  \A\x\|_2^2 + \|\P_{\mu}(\c)-\P_{\mu}(\c')\|_2^2 + 2\la \P_{\mu}(\c) -  \P_{\mu}(\c'),(\A \x' -\b) - (\A\x-\b)\ra \nn\\
&\overset{\step{173}}{=} &  A^2_s  \|\x'-\x\|_2^2 + \| \P_{\mu}(\c) - \P_{\mu}(\c')\|_2^2 + 2\la \P_{\mu}(\c)-\P_{\mu}(\c'),[\mu \p_1 + \P_{\mu}(\c')] - [\mu \p_1 +  \P_{\mu}(\c) ]\ra \nn\\
&\overset{}{=} &  A^2_s  \|\x'-\x\|_2^2 - \| \P_{\mu}(\c)-\P_{\mu}(\c')\|_2^2 + 2\la \P_{\mu}(\c)-\P_{\mu}(\c'),\mu \p_2-\mu \p_1 \ra \nn\\
&\overset{\step{174}}{\leq}&  A^2_s  \|\x'-\x\|_2^2 + 0 +0 ,
\eeq
\noi where step \step{172} uses the Pythagoras relation; step \step{173} uses Assumption \ref{ass:3} and the optimality conditions in (\ref{eq:opt:PcPc}); step \step{174} uses the convexity of $h(\cdot)$ that $\la \y'-\y,  \partial h(\y') - \mu \partial h(\y)\ra \geq 0$ for all $\y$.

\noi Finally, we have the following inequalities:
\beq
&&\|\nabla_{\x} \mathcal{G}(\x';\mu) - \nabla_{\x} \mathcal{G}(\x;\mu) \| \nn\\
&\overset{\step{172}}{=} &   \| [\nabla f(\x')  + \tfrac{1}{\mu} \A\trans (\c' - \P_{\mu}(\c')) ] - [ \nabla f(\x) + \tfrac{1}{\mu} \A\trans (\c - \P_{\mu}(\c))  ] \| \nn\\
&\overset{\step{173}}{\leq} &   \| \nabla f(\x') - \nabla f(\x)\| + \|\tfrac{1}{\mu} \A\trans (\A\x' - \b - \P_{\mu}(\c')) - \tfrac{1}{\mu} \A\trans (\A\x - \b - \P_{\mu}(\c))] \|  \nn\\
&\overset{\step{174}}{\leq } &   M_s \| \x - \x'\| + \tfrac{1}{\mu} \cdot\|\A\| \cdot \|[\A\x' - \A\x] + [ \P_{\mu}(\c) - \P_{\mu}(\c')] \|  \nn\\
&\overset{\step{175}}{\leq} &   M_s \| \x - \x'\| + \tfrac{1}{\mu} \cdot \|\A\| \cdot A_s \cdot\|\x-\x'\|, \nn
\eeq
\noi where step \step{172} uses the definition of $\nabla_{\x} \mathcal{G}(\x;\mu)$ in (\ref{eq:grad:G:mu}); step \step{173} uses the norm inequality; step \step{174} uses the fact that $f(\x)$ is restricted $M_s$-smooth as shown in Assumption \ref{ass:2} and norm inequality; step \step{175} uses Inequality (\ref{eq:lip:bound}).

\end{proof}

\subsection{Proof of Lemma \ref{lemma:first:bound}}

\label{app:lemma:first:bound}

\begin{proof} \textbf{(a)} We now bound $\| \y^{t+1} + \b - \A\x^{t+1} \|$ using these inequalities:
\beq \label{eq:eee:t1}
\| \y^{t+1} + \b - \A\x^{t+1} \| \overset{\step{172}}{=} \mu^{t} \|  \partial h(\y^{t+1})\| \overset{\step{173}}{\leq}  L_h \mu^{t},
\eeq
\noi where step \step{172} uses the optimality condition of $\y^{t+1}$ with $\y^{t+1}=\arg \min_{\y}\,h(\y) + \tfrac{1}{2\mu^{t}} \|\A\x^{t+1} - \b- \y\|_2^2$, which yields:
\beq \label{eq:opt:yyyy}
 \A\x^{t+1}- \b -\y^{t+1}   \in \mu^{t} \partial  h(\y^{t+1});
\eeq
\noi step \step{173} uses Assumption \ref{ass:1}.

\noi \textbf{(b)} We now bound $\|\y^{t+1}-\y^{t}\|$ using these inequalities:
\beq
\|\y^{t+1}-\y^{t}\| &\overset{\step{172}}{=}& \| \A(\x^{t+1} - \x^{t}) + \mu^{t-1} \partial h(\y^{t}) - \mu^t \partial h(\y^{t+1})\| \nn\\
 &\overset{\step{173}}{\leq}& \|\A\| \| \x^{t+1}  -\x^{t}\| + \|\mu^t \partial h(\y^{t+1}) \| + \| \mu^{t-1} \partial h(\y^{t})\| \nn\\
 &\overset{\step{174}}{\leq}&\|\A\| \| \x^{t+1}  -\x^{t}\| + 2\mu^{t-1} L_h ,\nn
\eeq
\noi where step \step{172} uses (\ref{eq:opt:yyyy}); step \step{173} uses the triangle inequality and norm inequality; step \step{174} uses $\|\partial h(\y)\| \leq L_h$ and $\mu^{t}\leq\mu^{t-1}$.

\noi \textbf{(c)} We now bound $\| \nabla_{\x} \mathcal{R}(\x^{t},\y^{t};\mu^{t}) \|$ using these inequalities:
\beq
\| \nabla_{\x} \mathcal{R}(\x^{t},\y^{t};\mu^{t}) \| &\overset{\step{172}}{=}& \| \nabla f(\x^{t})  + \tfrac{1}{\mu^{t}} \A\trans (\A\x^{t} - \y^{t} - \b)\|\nn\\
&\overset{}{\leq}&  \| \nabla f(\x^{t})\|  + \tfrac{1}{\mu^{t}}\|\A\| \| \A\x^{t} - \y^{t} - \b\|\nn\\
&\overset{\step{173}}{\leq}& L_f + \tfrac{ \mu^{t-1} }{\mu^{t}} L_h \|\A\|, \nn
\eeq
\noi where step \step{172} uses the definition of $\nabla_{\x} \mathcal{R}(\x^{t},\y^{t};\mu^{t})$ as in (\ref{eq:penalty:J}); step \step{173} uses \textbf{Part (a)} of this lemma.

\noi If we choose $\mu^t = \bar{\mu}$, we have: $\tfrac{\mu^{t-1}}{\mu^{t}} = 1,\,\text{and}\,\| \nabla_{\x}\mathcal{R}(\x^{t},\y^{t};\mu^{t}) \|\leq L_f +  L_h \|\A\|$. \\
\noi If we choose $\mu^{t} = \tfrac{\eta}{t+t_0}$, we have $\tfrac{\mu^{t-1}}{\mu^{t}}  = \frac{t + t_0}{ t -1 + t_0} \leq \max_{i=1}^{\infty} \frac{i+t_0}{i-1+t_0} \leq \frac{t_0+1}{t_0},\,\text{and}\,\|\nabla_{\x}\mathcal{R}(\x^{t},\y^{t};\mu^{t}) \|\leq L_f +   \frac{t_0+1}{t_0}L_h\|\A\|$.

\noi \textbf{(d)} We now bound the term $\|\varepsilons^t\|$ using these inequalities:
 \beq
 \|\varepsilons^t\|  &= & \|  (\tfrac{1}{\mu^{t-1}} - \tfrac{1}{\mu^{t}}) \cdot  \A\trans (\A\x^{t} - \b - \y^t) \| \nn\\
 &\overset{\step{172}}{\leq}& (\tfrac{1}{\mu^t} - \tfrac{1}{\mu^{t-1}}) \|\A\| \|\A\x^{t} - \b - \y^t\| \nn\\
 &\overset{\step{173}}{=}& (\tfrac{1}{\mu^t} - \tfrac{1}{\mu^{t-1}}) \|\A\|\cdot L_h \mu^{t-1},\nn
  \eeq
\noi where step \step{172} uses the norm inequality; step \step{173} uses (\ref{eq:eee:t1}).

\noi \textbf{(e)} We now bound $\mathcal{J}(\x^{t},\y^{t};\mu^{t}) -   \mathcal{J}(\x^{t},\y^{t};\mu^{t-1})$ using these inequalities:
\beq
 \mathcal{J}(\x^{t},\y^{t};\mu^{t}) -   \mathcal{J}(\x^{t},\y^{t};\mu^{t-1}) &\overset{\step{172}}{=}&  \tfrac{1}{2}\| \A\x^t - \b - \y^t\|_2^2 \cdot (\tfrac{1}{\mu^t} - \tfrac{1}{\mu^{t-1}} ) \nn\\
 &\overset{\step{173}}{\leq}&  \tfrac{1}{2}(\mu^{t-1} L_h)^2 (\tfrac{1}{\mu^t} - \tfrac{1}{\mu^{t-1}} )  =   \tfrac{1}{2}L_h^2 (\tfrac{  (\mu^{t-1})^2  }{\mu^t} -\mu^{t-1} ), \nn
\eeq
\noi where step \step{172} uses the definition of $\mathcal{J}(\x^{t},\y^{t};\mu^{t})\triangleq f(\x^t) + \tfrac{1}{2\mu^t} \| \A\x^t - \b - \y^t\|_2^2$; step \step{173} uses \textbf{Part (a)} of Lemma \ref{lemma:first:bound}. We now prove that $\Psi'\triangleq\sum_{t=0}^{\infty}\Psi^t = \tfrac{1}{2}  L_h^2 \sum_{t=0}^{\infty}  (\tfrac{(\mu^{t-1})^2}{\mu^{t}} - {\mu^{t-1}} ) $ is upper bounded by $\Lambda \leq \begin{scriptsize}\left\{
  \begin{array}{ll}
0, & \hbox{$\mu^t=\bar{\mu}$;} \\
\eta L_h^2   , & \hbox{$\mu^t =\eta (t+t_0)^{-1}$}
  \end{array}
\right.\end{scriptsize}$.

\noi \textbf{(f)} We now bound $[\sum_{t=1}^{\infty}\Psi^t]$. We discuss two cases for $\mu^t$.

\noi \quad\quad \textbf{Case 1)}. When $\mu^t = \bar{\mu}$, we have:
\beq
[\sum_{t=1}^{\infty}\Psi^t] \triangleq \tfrac{1}{2}  L_h^2\sum_{t=0}^{\infty} (\tfrac{(\mu^{t-1})^2}{\mu^{t}} - {\mu^{t-1}} )
= \tfrac{1}{2}  L_h^2\sum_{t=0}^{\infty} (\bar{\mu} - \bar{\mu} ) = 0 .\nn
\eeq

\noi \quad \quad\textbf{Case 2)}. When $\mu^t = \tfrac{\eta}{t+t_0}$, we have:
\beq
[\sum_{t=1}^{\infty}\Psi^t] &\overset{\step{172}}{=}& (\tfrac{1}{2}  L_h^2  )\cdot \sum_{t=1}^{\infty} (\tfrac{(\mu^{t-1})^2}{\mu^{t}} - {\mu^{t-1}} )   \overset{\step{173}}{=}  (\tfrac{1}{2}  L_h^2  )\cdot(\sum_{t=1}^{\infty}\frac{\eta}{(t+t_0-1)^2}) \nn\\
&\overset{\step{174}}{\leq}& (\tfrac{1}{2}  L_h^2   ) \cdot (\sum_{t=1}^{\infty}\frac{\eta}{t^2} )  \overset{\step{175}}{<}  (\tfrac{1}{2}  L_h^2  ) \cdot 2 \eta,\nn
\eeq
\noi where step \step{172} uses the definition of $\Psi^t\triangleq \tfrac{L_h^2}{2}   (\tfrac{(\mu^{t-1})^2}{\mu^{t}} - {\mu^{t-1}} )$; step \step{173} uses $\mu^t = \frac{\eta}{t+t_0}$; step \step{174} uses $t_0\geq 1$; step \step{175} uses $\sum_{t=1}^{\infty} \frac{1}{t^2} = \frac{\pi^2}{6}< 2$.


\end{proof}

\subsection{Proof of Lemma \ref{lemma:bound:negative:l2}}
\label{app:lemma:bound:negative:l2}

\begin{proof}

\textbf{(a)} We first now bound the term $\|\x^t - \bar{\x}\|$ using these inequalities:
\beq
\tfrac{V_s}{2}\|\x^t - \bar{\x}\|_2^2 {\overset{\step{172}}{\leq}}  F(\bar{\x})-F(\x^t) - \la \bar{\x} - \x^t,\partial F(\x^t) \ra {\overset{\step{173}}{\leq}}  0 + L_F\|\bar{\x}-\x^t\|,\nn
\eeq
\noi where step \step{172} uses the restricted strong convexity of $F(\cdot)$; step \step{173} uses $F(\bar{\x})\leq F(\x^t)$ and $\|\partial F(\x)\| \leq L_F$. Dividing both sides by $(\tfrac{V_s}{2}\|\bar{\x}-\x^t\|)$, we have:
\beq\label{eq:bound:xt1xs}
\|\x^t - \bar{\x}\| \leq \tfrac{2 L_F}{V_s}.
\eeq
\noi We now now bound the term $\frac{V_s }{2}\|\x^t - \bar{\x}\|_2^2$ using these inequalities:
\beq
&& \tfrac{V_s}{2}\|\x^t - \bar{\x}\|_2^2    \nn \\
 &{\overset{\step{172}}{\leq}}  & \la \nabla_{\x} \mathcal{G}(\x^t;{\mu^{t-1}}),\,\x^t -\bar{\x}\ra + \mathcal{G}(\bar{\x};{\mu^{t-1}}) - \mathcal{G}(\x^t;{\mu^{t-1}})     \nn\\
&{\overset{\step{173}}{\leq }}&    \la \r^t +\varepsilons^t ,\,\x^t-\bar{\x} \ra +  [ F(\bar{\x}) - F(\x^t) + \tfrac{{\mu^{t-1}}}{2}L_h^2]   \nn  \\
&{\overset{\step{174}}{\leq }}&      \la \r^t ,\x^t-\bar{\x} \ra  + \|\varepsilons^t\| \cdot \|\bar{\x}-\x^t\| +  [ F(\bar{\x}) - F(\x^t) + \tfrac{{\mu^{t-1}}}{2}L_h^2] \nn \\
&{\overset{\step{175}}{\leq }}&      \la \r^t ,\x^t-\bar{\x} \ra  + \underbrace{ (\tfrac{1}{\mu^{t}} - \tfrac{1}{\mu^{t-1}})L_h \mu^{t-1} \| \A\|  \cdot \tfrac{2L_F}{V_s} +  [ F(\bar{\x}) - F(\x^t)  + \tfrac{{\mu^{t-1}}}{2}L_h^2]}_{\triangleq \Upsilon^t} \nn
\eeq
\noi where step \step{172} uses the restricted strong convexity of $\mathcal{G}(\x;\mu^{t-1})$; step \step{173} uses the the relation between $\nabla_{\x} \mathcal{G}(\x^t;{\mu^{t-1}})$ and $\nabla_{\x}\mathcal{R}(\x^t,\y^t;{\mu^{t}})$ and \textbf{Part (a)} in Lemma (\ref{lemma:lip:variable:x}) that
\beq
\mathcal{G}(\bar{\x};{\mu^{t-1}}) \leq F(\bar{\x}),\, \mathcal{G}(\x^t;{\mu^{t-1}})  \geq  F(\x^t ) - \frac{\mu^{t-1}}{2} L_h^2;\nn
\eeq
\noi step \step{174} uses the norm inequality; step \step{175} uses Inequality (\ref{eq:bound:xt1xs}) and the inequality in \textbf{Part (d)} of Lemma \ref{lemma:first:bound} that $\|\varepsilons^t\| \leq  (\tfrac{1}{\mu^{t}} - \tfrac{1}{\mu^{t-1}})L_h \mu^{t-1} \| \A\|$.

\noi \textbf{(b)} When $\mu^t = \bar{\mu}$, we have the following results:
\beq
\Upsilon^t &\triangleq & F(\bar{\x}) - F(\x^t) + \tfrac{1}{2} \mu^{t-1}L_h^2 +  \tfrac{2 L_F L_h}{V_s}  \| \A\|(\tfrac{\mu^{t-1}}{\mu^{t}} - 1) \nn\\
&\overset{\step{172}}{=} & F(\bar{\x}) - F(\x^t) + \tfrac{1}{2} \bar{\mu}L_h^2.\nn\\
&\overset{\step{173}}{\leq } &  F(\bar{\x}) -  [\min_{i=1}^t F(\x^i) ] + \tfrac{1}{2} \bar{\mu}L_h^2 ,\nn
\eeq
\noi where step \step{172} uses $\frac{\mu^{t+1}}{\mu^{t}}=1$; step \step{173} uses $[\min_{i=1}^t F(\x^i) ]\leq F(\x^t)$.

\noi \textbf{(c)} When $\mu^t = \frac{\eta}{t+t_0}$, we have the following results:
\beq
 \sum_{t=1}^T  \Upsilon^t  &\overset{\step{172}}{=} & \sum_{t=1}^T \left(F(\bar{\x}) - F(\x^t) + \tfrac{1}{2} \mu^{t-1}L_h^2 +  \tfrac{2 L_F L_h}{V_s}  \| \A\|(\tfrac{\mu^{t-1}}{\mu^{t}} - 1) \right)\nn\\
&\overset{\step{173}}{\leq} &  T F(\bar{\x}) -  T  [\min_{t=1}^T F(\x^t) ]  +  \sum_{t=1}^T \left( \frac{L_h^2 }{2} \frac{\eta}{t+t_0-1} +  \tfrac{2 L_F L_h}{V_s}  \| \A\| \frac{1}{t + t_0 - 1} \right)\nn\\
&\overset{\step{174}}{\leq} & T  F(\bar{\x}) - T  [\min_{t=1}^T F(\x^t)] +  \left[1 + \ln(T)\right] \cdot \left( \frac{\eta L_h^2 }{2}  +  \tfrac{2 L_F L_h}{V_s}  \| \A\| \right),\nn
\eeq
\noi where step \step{172} uses the definition of $\Upsilon^t$ as shown in Lemma \ref{lemma:bound:negative:l2}; step \step{173} uses $\max_{t=1}^T  [- F(\x^t) ]  = - [\min_{t=1}^T F(\x^t)]$, the fact that: $\frac{\mu^{t-1}}{\mu^t} - 1 = \frac{t+t_0}{t+t_0-1} - 1   = \frac{1}{t + t_0 - 1}$; step \step{174} uses $t_0\geq 1$, and the fact that:
\beq
\sum_{t=1}^T \frac{1}{t+t_0 - 1 } \leq \sum_{t=1}^T \frac{1}{t}  \leq  1 + \ln(T).\nn
\eeq
\noi Using the definition of $C_{\Upsilon} \triangleq \frac{\eta L_h^2 }{2}  +  \tfrac{2 L_F L_h}{V_s}  \| \A\|$, we finish the proof of this lemma.

\end{proof}

\section{Proofs for Section \ref{sect:rate}}
\label{app:sect:rate}

\subsection{Proof of Theorem \ref{theorem:convergence:iht}}
\label{app:theorem:convergence:iht}

\begin{proof}

We denote $\r^t \triangleq \nabla_{\x} \mathcal{R}(\x^t,\y^t;\mu^t)$ and $H^t=  A_s^2/\mu^t + M_s + \theta$.

\textbf{(a)} We focus on the $\x$-subproblem. We have from Problem (\ref{eq:subprob:x:1}) that:
\beq
\la \r^t,\,\x^{t+1}-\x^{t}\ra   +\tfrac{H^t}{2}\|\x^{t+1}-\x^{t}\|_2^2 \leq \la \r^t,\,\x^{t}-\x^{t}\ra   +\tfrac{H^t}{2}\|\x^{t}-\x^{t}\|_2^2 = 0.\nn
\eeq
\noi Since $\mathcal{R}(\x^t,\y^t;\mu^t)$ is restricted $(A_s^2/\mu^t + M_s)$-smooth \textit{w.r.t.} $\x$, we have:
\beq
\mathcal{R}(\x^{t+1},\y^t;\mu^t) \leq \mathcal{R}(\x^t,\y^t;\mu^t)+\la \r^t, \x^{t+1}- \x^{t} \ra + \tfrac{A_s^2/\mu^t + M_s}{2}\|\x^t-\x^{t+1}\|_2^2.\nn
\eeq
\noi We observe that the following equality holds:
\beq
\mathcal{R}(\x^{t+1},\y^t;\mu^t) -\mathcal{R}(\x^t,\y^t;\mu^t) = \mathcal{J}(\x^{t+1},\y^{t};\mu^t) - \mathcal{J}(\x^{t},\y^{t};\mu^t)\nn
\eeq
\noi Summing up these three inequalities, we have:
\beq
\mathcal{J}(\x^{t+1},\y^{t};\mu^t) - \mathcal{J}(\x^{t},\y^{t};\mu^t) \leq -\tfrac{\theta}{2}\|\x^{t+1}-\x^{t}\|_2^2. \label{eq:opt:x}
\eeq
\noi We now focus on the $\y$-subproblem. We derive the following inequalities for all $\y\in\mathbb{R}^m$:
\beq
&&\mathcal{J}(\x^{t+1},\y^{t+1};\mu^{t}) - \mathcal{J}(\x^{t+1},\y;\mu^{t})\nn\\
&\overset{\step{172}}{\leq}& - \frac{1}{2\mu^t} \|\y^{t+1}-\y\|_2^2 -  \la \y-\y^{t+1}, \partial_{\y} \mathcal{J}(\x^{t+1},\y^{t+1};\mu^{t})   \ra  \nn\\
& \overset{\step{173}}{=} &- \frac{1}{2\mu^t} \|\y^{t+1}-\y\|_2^2   \nn\\
& \overset{\step{174}}{\leq} &- \frac{1}{2\mu^1} \|\y^{t+1}-\y\|_2^2 ,  \nn
\eeq
\noi where step \step{172} uses the fact that $\mathcal{J}(\x^{t+1},\y;\mu^{t})$ is $\frac{1}{\mu}$-strongly convex \textit{w.r.t.} $\y$; step \step{173} uses the optimality of $\y^{t+1}$ that $0\in\partial_{\y} \mathcal{J}(\x^{t+1},\y^{t+1};\mu^{t})$; step \step{174} uses the fact that the sequence $\{\mu^t\}_{t=1}^{\infty}$ is non-increasing. Letting $\y = \y^t$, we obtain:
\beq \label{eq:opt:y}
\mathcal{J}(\x^{t+1},\y^{t+1};\mu^{t}) - \mathcal{J}(\x^{t+1},\y^{t};\mu^{t})    \leq  - \frac{1}{2\mu^1}\|\y^{t+1}-\y^t\|_2^2
\eeq
\noi Using the continuity of $\mathcal{J}(\x^{t},\y^{t};\mu)$ \textit{w.r.t.} $\mu$ as shown in \textbf{Part (e)} of Lemma \ref{lemma:first:bound}, we obtain:
\beq\label{eq:connut:mu}
0\leq \mathcal{J}(\x^{t},\y^{t};\mu^{t}) -   \mathcal{J}(\x^{t},\y^{t};\mu^{t-1}) \leq \tfrac{L_h^2}{2}   (\tfrac{(\mu^{t-1})^2}{\mu^{t}} - {\mu^{t-1}} ) \triangleq  \Psi^t.
\eeq
\noi Summing up Inequalities (\ref{eq:opt:x}), (\ref{eq:opt:y}), and (\ref{eq:connut:mu}) together, we have:
\beq \label{eq:sss:ddd:ccc:1}
&& \frac{1}{2\mu^1}\|\y^{t+1}-\y^t\|_2^2 + \tfrac{\theta}{2}\|\x^{t+1}-\x^{t}\|_2^2  \nn\\
&\leq&    \mathcal{J}(\x^{t},\y^{t};\mu^{t-1}) - \mathcal{J}(\x^{t+1},\y^{t+1};\mu^{t}) + \tfrac{1}{2}  L_h^2 (\tfrac{(\mu^{t-1})^2}{\mu^{t}} - {\mu^{t-1}} ) \nn \\
&=&\mathcal{J}^{t} -\mathcal{J}^{t+1} + \Psi^t,
\eeq
 \noi where $\mathcal{J}^{t+1} \triangleq\mathcal{J}(\x^{t+1},\y^{t+1};\mu^{t})$.

\noi \textbf{(b)} Summing up the inequality in (\ref{eq:sss:ddd:ccc:1}) over $t=1,2,...,T$, we have:
\beq
&& \textstyle \sum_{t=1}^T \tfrac{1}{2\mu^1}\|\y^{t+1}-\y^{t}\|_2^2 + \sum_{t=1}^T \tfrac{\theta}{2}\|\x^{t+1}-\x^{t}\|_2^2\nn\\
& \leq &   \textstyle \mathcal{J}^1 - \mathcal{J}^{T+1} + [\sum_{t=1}^T \Psi^t]  \overset{\step{172}}{\leq}   \mathcal{J}^1 - \mathcal{J}^{T+1} + \eta L_h^2 \triangleq C < + \infty,\nn
\eeq
\noi where step \step{172} uses $\sum_{t=1}^T \Psi^t < \sum_{t=1}^{\infty} \Psi^t \leq \eta L_h^2 $ which is shown in \textbf{Part (f)} of Lemma \ref{lemma:first:bound}.

\noi \textbf{(c)} As a result, there exists an index $\bar{t}$ with $1 \leq \bar{t} \leq T$ such that: $\tfrac{1}{2\mu^1}\|\y^{\bar{t}+1}-\y^{\bar{t}}\|_2^2 + \tfrac{\theta}{2}\|\x^{\bar{t}+1}-\x^{\bar{t}}\|_2^2 \leq \frac{C}{T}$, leading to:
\beq
\|\y^{\bar{t}+1}-\y^{\bar{t}}\|_2^2 + \|\x^{\bar{t}+1}-\x^{\bar{t}}\|_2^2 \leq \frac{2C}{T\cdot \min(\theta,(\mu^1)^{-1}) }. \label{eq:lll:1}
\eeq
\noi Letting $\Gamma_x(\x,\y;\mu) \triangleq \dist^2(\x,\arg \min_{\x'}\mathcal{M}(\x',\x,\y;\mu))$ and $\Gamma_y(\x,\y;\mu) \triangleq \dist^2(\y, \arg \min_{\y'}\mathcal{J}(\x,\y';\mu))$, we have:
\beq
\|\x^{\bar{t}+1}-\x^{\bar{t}}\|_2^2 + \|\y^{\bar{t}+1}-\y^{\bar{t}}\|_2^2 \geq
\Gamma_x(\x^{\bar{t}},\y^{\bar{t}};\mu) + \Gamma_y(\x^{\bar{t}},\y^{\bar{t}};\mu) \label{eq:lll:2}
\eeq
\noi for all $\bar{t}\geq 1$ and some sufficiently small $\mu=\mu^{\bar{t}}>0$. Combining Inequality (\ref{eq:lll:1}) and Inequality (\ref{eq:lll:2}), we have:
\beq
\Gamma_x(\x^{\bar{t}},\y^{\bar{t}};\mu^{\bar{t}}) + \Gamma_y(\x^{\bar{t}},\y^{\bar{t}};\mu^{\bar{t}}) \leq \frac{2C}{T\cdot \min(\theta,(\mu^1)^{-1}) }\nn
\eeq

\noi Therefore, we conclude that Algorithm \ref{algo:main} finds an $\epsilon$-approximate Lipschitz stationary point of Problem (\ref{eq:main}) in at most $T$ iterations, where $T\leq \lceil  \frac{2C}{\epsilon \min(\theta,(\mu^1)^{-1}) }\rceil$.

\end{proof}

\subsection{Proof of Lemma \ref{lemma:iht:bound}}
\label{app:lemma:iht:bound}

\begin{proof}

We define $H^t\triangleq   A_s^2 /\mu^t + M_s+\theta$, $\x^t_+  \triangleq \x^t -  \r^t/H^t \in \mathbb{R}^n$, $\mathrm{J} \triangleq \{i\,|\,\x^{t+1}_i \neq 0\}$, and $\mathrm{J}^c\triangleq \{i\,|\,\x^{t+1}_i = 0\}$.

\noi \textbf{(a)} Due to the optimality of $\x^{t+1}$ in (\ref{eq:subprob:x:1}) that: $\x^{t+1} = \arg \min_{\|\x\|_0 \leq s}\, \tfrac{1}{2}\| \x - \x^t_+\|_2^2$, we have $\| \x^{t+1} - \x^t_+  \| \leq \| \x - \x^t_+\|$ for all $\|\x\|_0\leq s$. Given that $\|\x^{t}\|_0\leq s$, we let $\x=\x^{t}$, resulting in:
\beq \label{eq:opt:zt}
\|\x^{t+1} - \x^t_+ \| \leq \| \x^t - \x^t_+ \|.
\eeq
We derive the following inequalities:
\beq
\begin{split}
\|\x^{t+1} - \x^{t}\| \overset{\step{172}}{\leq}\|\x^{t+1} - \x^t_+\| + \|\x^t_+ - \x^{t}\|  \overset{\step{173}}{\leq} \tfrac{1}{H^t}\| \r^{t} \| + \tfrac{1}{H^t}\| \r^{t}\| = \tfrac{2}{H^t}\| \r^{t} \| , \nn
\end{split}
\eeq
\noi where step \step{172} uses the triangle inequality; step \step{173} uses (\ref{eq:opt:zt}).

\noi \textbf{(b)} We have the following inequalities:
\beq \label{eq:xt:xt1:xt:xt1}
\la \x^t_+ , \x^t_+\ra & = &\la [\x^t_+]_{\mathrm{J}} , [\x^t_+]_{\mathrm{J}}\ra + \| [\x^t_+]_{\mathrm{J}^c}\|_2^2  \geq   \la [\x^t_+]_{\mathrm{J}} , [\x^t_+]_{\mathrm{J}}\ra+ 0  \nn\\
&\overset{\step{172}}{=} & \la [\x^t_+]_{\mathrm{J}} , [\x^{t+1}]_{\mathrm{J}}\ra   \nn\\
&\overset{\step{173}}{=} & \la [\x^t_+]_{\mathrm{J}} , [\x^{t+1}]_{\mathrm{J}}\ra + \la [\x^t_+]_{\mathrm{J}^c} , [\x^{t+1}]_{\mathrm{J}^c}\ra   \nn\\
&\overset{\step{174}}{=} & \la  \x^t_+ ,\x^{t+1}\ra,
\eeq
\noi where step \step{172} uses the fact $[\x^t_+]_{\mathrm{J}} = [\x^{t+1}]_{\mathrm{J}}$; step \step{173} uses $[\x^{t+1}]_{\mathrm{J}^c}=\zero$; step \step{174} uses $\mathrm{J} \cup \mathrm{J}^c=\{1,2,...,n\}$.

\noi \textbf{(c)} We derive the following inequalities:
\beq
&& \tfrac{H^t}{2}\|\x^{t+1} - \bar{\x}\|_2^2 - \tfrac{H^t}{2}\|\x^t -\bar{\x}\|_2^2 - \tfrac{H^t}{2}\|\x^{t+1} - \x^{t} \|_2^2  \nn\\
&\overset{\step{172}}{=} &    H^t  \la   \x^t - \x^{t+1}, \bar{\x}-\x^t\ra   \nn\\
&\overset{\step{173}}{= } &  H^t \la   \x^t_+ - \x^{t+1}  , \bar{\x}-\x^t\ra +   H^t \la  \r^t/ H^t, \bar{\x}-\x^t\ra \nn \\
&\overset{\step{174}}{\leq } &   H^t\la \x^t_+ - \x^{t+1}  , \bar{\x}-\x^t\ra +  \Upsilon^t  -\tfrac{V_s}{2}\|\x^t -\bar{\x}\|_2^2\label{eq:lemma:conv:1}
\eeq
\noi where step \step{172} uses the Pythagoras relation that $\|\a-\b\|_2^2  -  \|\c-\b\|_2^2=\|\a-\c\|_2^2 - 2 \la \a-\c,\b-\c\ra$ for all $\a,\,\b,\,\c$; step \step{173} uses $\x^t  = \x^t_+ +  \r^t/H^t$; step \step{174} uses \textbf{Part (b)} of this lemma.

\noi We now bound the first term  of the right-hand side in Inequality (\ref{eq:lemma:conv:1}) using the following inequalities:
\beq\label{eq:lemma:conv:2}
H^t\la   \x^t_+ - \x^{t+1}  , \bar{\x}-\x^t\ra  &\overset{\step{172}}{=}  & H^t\la   \x^t_+ - \x^{t+1}  , \bar{\x}-\x^t - \x^t_+ + \x^t_+\ra  \nn \\
&\overset{\step{173}}{\leq} & H^t\la   \x^t_+ - \x^{t+1}  , \bar{\x}-\x^t  + \x^t_+\ra   \nn \\
&\overset{}{=}  & H^t\la   \x^t_+ - \x^{t+1}  , \bar{\x}\ra + H^t\la   \x^t_+ - \x^{t+1}  , \x^t_+ -\x^t  \ra    \nn\\
&\overset{\step{174}}{\leq} & H^t ( \|\x^t_+ - \x^{t+1}\|\|\bar{\x}\| + \|\x^t_+ - \x^{t+1}\|\|\x^t_+-\x^t\| )    \nn\\
&\overset{\step{175}}{\leq}  & H^t\|\x^t_+ - \x^{t}\| ( \|\bar{\x}\| + \|\x^t_+-\x^t\| )    \nn\\
&\overset{\step{176}}{=}  & \|\r^t\| ( \|\bar{\x}\| +  \tfrac{\|\r^t\|}{H^t} ),
\eeq
\noi where step \step{172} uses the fact that $\bar{\x}-\x^{t} = (\bar{\x} -\x^t_+) + (\x^t_+-\x^{t})$; step \step{173} uses Inequality (\ref{eq:xt:xt1:xt:xt1}); step \step{174} uses the Cauchy-Schwarz Inequality, step \step{175} uses (\ref{eq:opt:zt}); step \step{176} uses $\|\x^t_+ - \x^t\| = \tfrac{1}{H^t}\|\r^t\|$.

Finally, we have from (\ref{eq:lemma:conv:1}) and (\ref{eq:lemma:conv:2}):
\beq
&&\tfrac{H^t}{2}\|\x^{t+1} - \bar{\x}\|_2^2 - (\tfrac{H^t}{2} - \tfrac{V_s}{2}) \|\x^t -\bar{\x}\|_2^2  \nn\\
&\overset{}{\leq}  & \tfrac{H^t}{2}\|\x^{t+1} - \x^{t}\|_2^2 +  \tfrac{1}{H^t}\|\r^{t}\|_2^2 + \Upsilon^t + \|\r^{t}\| \|\bar{\x} \|    \nn \\
&\overset{\step{172}}{\leq}  & \tfrac{2+1}{H^t}\|\r^{t}\|_2^2  + \Upsilon^t +  \|\r^{t}\| \|\bar{\x} \| , \nn
\eeq
\noi where step \step{172} uses $\|\x^{t+1} - \x^{t}\| \leq \frac{2}{H^t}\|\r^t\|$ as shown in \textbf{Part (a)} of this lemma.


\end{proof}

\subsection{Proof of Theorem \ref{theorem:IHT:stepsize:constant}}
\label{app:theorem:IHT:stepsize:constant}

\begin{proof}

Assume constant stepsizes are used with $\mu^t = \bar{\mu}$ for all $t\geq 1$.

\noi We define $H \triangleq A_s^2 /\bar{\mu} + M_s +\theta$, $\gamma \triangleq 1 - \frac{V_s}{H}$, $\Delta_{\x}^t \triangleq \|\x^{t} - \bar{\x}\|_2^2$, and $\Delta_{F}^t \triangleq [\min_{i=1}^t F(\x^i)]  -   F(\bar{\x})$.

\noi First, using \textbf{Part (c)} in Lemma \ref{lemma:first:bound}, we have: $\|\r^t\|\leq L_F$.

\noi Second, using \textbf{Part (b)} Lemma \ref{lemma:bound:negative:l2}, we have the upper bound of $\Upsilon^t$ that $\forall t,\,\Upsilon^t \leq \tfrac{\bar{\mu}}{2} L_h^2 - \Delta_{F}^t$.

\noi Third, it holds that $H^t=H$ for all $t\geq 1$.

\noi \textbf{(a)} Using the inequality in \textbf{Part (b)} in Lemma \ref{lemma:iht:bound}, we have the following recursive formulation:
\beq \label{eq:iht:hhhhhh}
\Delta_{\x}^{t+1}  &\leq &   (1 -  \tfrac{V_s}{H})\Delta_{\x}^{t} +  \tfrac{6}{H^2}\|\r^{t}\|_2^2  + \tfrac{2}{H}\Upsilon^t + \tfrac{2}{H} \|\r^t \|\cdot \|\bar{\x}\|  \nn \\
& \overset{\step{172}}{\leq}  &   \gamma \Delta_{\x}^{t} +  \tfrac{6 L_F^2}{H^2} + \tfrac{   \bar{\mu}L_h^2 - 2\Delta_{F}^t }{H} + \tfrac{2L_F \|\bar{\x}\|}{H},\nn
\eeq
\noi where step \step{172} uses the definition of $\gamma \triangleq 1 - \frac{V_s}{H}$, $\|\r^t\|\leq L_F$, and $\Upsilon^t\leq \tfrac{1}{2} \bar{\mu}L_h^2 - \Delta_{F}^t$.

\noi \textbf{(b)} Let $T\geq 1$ be any integer. Applying Lemma \ref{lemma:two:non:negative:sequences} with $\Phi^t=\Delta_{\x}^t$ and $\Lambda^t =    \tfrac{6L_F^2}{H^2} + \frac{ \bar{\mu}L_h^2 - 2\Delta_{F}^t }{H} + \frac{2 L_F \|\bar{\x}\|}{H}$, we have:
\beq \label{eq:final:IHT:mu:1}
\Delta_{\x}^{T+1}     &\leq &  \gamma^T \Delta_{\x}^{1} + \tfrac{1}{1-\gamma}\cdot \max_{t=1}^T \left(    \tfrac{6L_F^2}{H^2} + \tfrac{ \bar{\mu}L_h^2 - 2\Delta_{F}^t }{H} + \tfrac{2 L_F \|\bar{\x}\| }{H}  \right) \nn\\
& \overset{\step{172}}{=}  &   \gamma^T \Delta_{\x}^{1} +    \tfrac{6 L_F^2}{ V_s H} +  \tfrac{\bar{\mu}L_h^2}{ V_s }  - \tfrac{2\Delta_{F}^T}{V_s}  + \tfrac{2 L_F \|\bar{\x}\| }{  V_s } \nn\\
& \overset{\step{173}}{\leq }  &  \gamma^T \Delta_{\x}^{1}  +  \tfrac{6L_F^2 }{ V_s A_s^2} \bar{\mu} +  \tfrac{L_h^2}{  V_s }  \bar{\mu} - \tfrac{ 2\Delta_{F}^T }{V_s}  +  \tfrac{2 L_F \|\bar{\x}\| }{ V_s }  \nn\\
& \overset{\step{174}}{=}  &   \tfrac{2}{V_s} \left( K_1 \gamma^T + D_1 \bar{\mu}  - \Delta_{F}^T +  L_F \|\bar{\x}\|\right),
\eeq
\noi where step \step{172} uses $\max_{t=1}^T  [-\Delta_{F}^t] = -\Delta_{F}^T$ since $\Delta_{F}^1 \geq \Delta_{F}^2\geq ... \geq \Delta_{F}^T \geq 0$ and $\gamma \triangleq 1 - \frac{V_s}{H}$; step \step{173} uses $H \triangleq A_s^2 /\bar{\mu} + M_s +\theta \geq A_s^2 /\bar{\mu}$; step \step{174} uses the definitions of $K_1\triangleq \frac{1}{2}V_s\Delta_{\x}^1$ and $D_1 \triangleq   \tfrac{3L_F^2 }{A_s^2} +  \tfrac{1}{2}L_h^2$.

\noi We now focus on (\ref{eq:final:IHT:mu:1}). Using the fact that $\Delta_{\x}^{T+1}\geq 0$, we obtain: $\Delta_{F}^T  \leq K_1  \gamma^t    +  D_1  \bar{\mu} +   L_F \|\bar{\x}\|$.

Using the fact that $\Delta_{F}^T \geq 0$, we obtain: $\Delta_{\x}^{T+1}   \leq \left( K_1 \gamma^T      + D_1 \bar{\mu} +   L_F \|\bar{\x}\| \right) \frac{2}{V_s}$.

\end{proof}

\subsection{Proof of Theorem \ref{theorem:IHT:stepsize:diminishing}}
\label{app:theorem:IHT:stepsize:diminishing}

\begin{proof}

Assume diminishing stepsizes are used with $\mu^t = \tfrac{\eta}{t+t_0}$ for all $t\geq 1$, where $\eta = \tfrac{A_s^2}{V_s}$.

We define $H^t \triangleq A_s^2/\mu^t + M_s + \theta$, $\Delta_{\x}^t \triangleq \|\x^{t} - \bar{\x}\|_2^2$, and $\Delta_{F}^t \triangleq [\min_{i=1}^t F(\x^i)]  -   F(\bar{\x})$.

First, using \textbf{Part (c)} in Lemma \ref{lemma:first:bound}, we have: $\|\r^t\|\leq L'_F$.

Second, using Lemma \ref{lemma:bound:negative:l2}, we have: $\sum_{t=1}^T \Upsilon^t \leq  C_{\Upsilon}\left(1 + \ln(T)\right) -  T  \Delta_{F}^T$ for any $T\geq 1$.

Third, using the definition of $H^t$ and the choice of $\eta=\frac{A_s^2}{V_s}$, we have:
\beq
H^t = V_s \eta/\mu^t + M_s + \theta = V_s(t+t_0) + M_s + \theta, \label{eq:St1:H}
\eeq

\noi \textbf{(a)} We have the following inequalities:
\beq\label{eq:IHT:diminishing:to:sum}
 \tfrac{1}{2}(H^{t+1}  - V_s) \Delta_{\x}^{t+1} &\overset{\step{172}}{\leq} &  \tfrac{1}{2} H^{t}  \Delta_{\x}^{t+1}\nn\\
&\overset{\step{173}}{\leq} & \tfrac{1}{2}(H^t- V_s)\Delta_{\x}^t  + \tfrac{3\|\r^{t}\|_2^2}{H^t} +  \Upsilon^t +  \|\r^{t}\| \|\bar{\x} \| \nn\\
&\overset{\step{174}}{\leq} & \tfrac{1}{2}(H^t - V_s)\Delta_{\x}^t +
\tfrac{3   (L_F')^2  }{ V_s  \cdot t } +  \Upsilon^t  +   L_F' \|\bar{\x} \|,
\eeq
\noi where step \step{172} uses $H^{t+1} =  H^{t} +  V_s$, which can be implied by Equation (\ref{eq:St1:H}); step \step{173} uses \textbf{Part (b)} in Lemma \ref{lemma:iht:bound}; step \step{174} uses $\|\r^t\| \leq L_F'$ and $H^t \geq A_s^2 / \mu^t = A_s^2 (t+t_0) / \eta \geq A_s^2 t / \eta = V_s t$.

\noi \textbf{(b)} Let $T\geq 1$ be any integer. Summing Inequality (\ref{eq:IHT:diminishing:to:sum}) over $t=1,2,...,T$, we have:
\beq \label{eq:final:IHT:mu:last}
0  &\leq& \textstyle  -  \tfrac{1}{2}(H^{T+1} - V_s)\Delta_{\x}^{T+1} + \tfrac{1}{2}(H^1-V_s)\Delta_{\x}^{1}   +   \tfrac{3   (L_F')^2}{V_s}  \sum_{t=1}^T \tfrac{1}{t}   + \sum_{t=1}^T \Upsilon^t + T L_F' \|\bar{\x} \|  \nn\\
&\overset{\step{172}}{\leq} & \textstyle - \tfrac{V_s}{2} (T+1)\Delta_{\x}^{T+1} + \tfrac{H^1}{2} \Delta_{\x}^1    +  [ \tfrac{3 (L_F')^2}{ V_s } + C_{\Upsilon}]  ( \ln(T)+1)   -  T\Delta_{F}^T + T L_F' \|\bar{\x} \|  \nn\\
&\overset{\step{173}}{\leq} & -    \tfrac{V_s}{2} (T+1) \Delta_{\x}^{T+1} + K_2       +   D_2    (\ln(T)+1)      -  T\Delta_{F}^T  +  T L_F' \|\bar{\x} \|,
\eeq
\noi where step \step{172} uses $H^{T+1} - V_s =  V_s(T+1+t_0) + M_s + \theta  - V_s \geq  V_s (T+1) + M_s + \theta \geq  V_s (T+1)$, $-\tfrac{1}{2}V_s \Delta_{\x}^1\leq 0$, the fact that $\sum_{t=1}^T  \frac{1 }{  t } \leq \ln(T)+1$, and the upper bound $\sum_{t=1}^T\Upsilon^t \leq  C_{\Upsilon}\left(1 + \ln(T)\right) -  T  \Delta_{F}^T$; step \step{173} uses the definition of $K_2 \triangleq \tfrac{H^1}{2}\Delta_{\x}^1$, and the definition of $D_2 \triangleq \tfrac{3(L_F')^2}{ V_s } +  C_{\Upsilon}$.

\noi We now focus on (\ref{eq:final:IHT:mu:last}). Using the fact that $\frac{V_s}{2}(T+1)\Delta_{\x}^{T+1}\geq 0$, we obtain: $T\Delta_F^T \leq K_2     +  T L_F' \|\bar{\x} \|   +   D_2    (\ln(T)+1)$, leading to $\Delta_F^T \leq \frac{K_2}{T}     +  L_F' \|\bar{\x} \|   +   D_2 \frac{\ln(T)+1}{T}$.

Using the fact that $T\Delta_{F}^T \geq 0$, we obtain: $\tfrac{V_s}{2} (T+1) \Delta_{\x}^{T+1} \leq K_2     +  T L_F' \|\bar{\x} \|   +   D_2    (\ln(T)+1)$, leading to $\Delta_{\x}^{T+1} \leq (\frac{K_2}{T+1}   +   D_2 \frac{\ln(T)+1}{T+1}  +  L_F' \|\bar{\x} \|   ) \frac{2}{V_s}$.

\end{proof}

\subsection{Proofs for Lemma \ref{lemma:bcd:simple:basic}}

\label{app:lemma:bcd:simple:basic}

 \begin{proof}

We denote $\Omega_n^k\triangleq\{\mathcal{\B}_{(i)}\}_{i=1}^{C_n^k}$ as all the possible combinations of the index vectors choosing $k$ items from $n$ with $\mathcal{\B}_{i} \in \mathbb{N}^{k},\,\forall i$. For any vector $\x\in\mathbb{R}^n$, we have:
 \beq
 \textstyle
\sum_{\B \in \Omega_n^k}\,\x\trans (\UB \UB\trans ) \z  \overset{\step{172}}{=} \sum_{\B \in \Omega_n^k}\,\la \x_{\B} ,\z_{\B}\ra  \overset{\step{173}}{=}  C_n^k \frac{k}{n} \la \x,\z\ra,\nn
 \eeq
\noi where step \step{172} uses $\UB\trans \x=\x_{\B}$ and $\UB\trans \z=\z_{\B}$; step \step{173} uses the basic induction that every entry $(\x_i\cdot\z_i)$ is present within the term $(\sum_{\B \in \Omega_n^k}\,\la \x_{\B} ,\z_{\B}\ra)$ for a total of $(C_n^k \cdot\frac{k}{n})$ times for all $i\in[n]$.

\noi Given $\B$ is chosen from $\Omega_n^k$ randomly and uniformly, we have: $\E_{\B} [\|\x_{\B}\|_2^2] = \frac{1}{C_n^k} \sum_{i=1}^{C_n^k}\|\x_{\B_{i}}\|_2^2 = \frac{k}{n}\|\x\|_2^2$.

 \end{proof}

\subsection{Proof of Theorem \ref{theorem:convergence:bcd}}
\label{app:theorem:convergence:bcd}

\begin{proof}

We denote $\r^t \triangleq \nabla_{\x} \mathcal{R}(\x^t,\y^t;\mu^t)$, $\H^t= (\A\trans \A + \theta_1 \I_n) /\mu^t +  \tilde{\M} + \theta_2 \I_n$, and $\theta = \frac{\theta_1}{\mu^1}+\theta_2$.

\textbf{(a)} We focus on the $\x$-subproblem. We have from Problem (\ref{eq:subprob:x:2}) that:
\beq
\E_{\xi^t}[\la \r^t,\x^{t+1}-\x^{t}\ra+\tfrac{1}{2}\|\x^{t+1}-\x^{t}\|_{\H^t}^2 ] \leq \E_{\xi^t}[\la \r^t,\x^{t}-\x^{t}\ra   +\tfrac{1}{2}\|\x^{t}-\x^{t}\|_{\H^t}^2] = 0.\nn
\eeq

\noi Using Assumption \ref{ass:2} and the inherent structure of the function $\mathcal{R}(\mathbf{x}^t,\mathbf{y}^t;\mu^t)$, we have:
\beq
\mathcal{R}(\x^{t+1},\y^t;\mu^t) \leq \mathcal{R}(\x^t,\y^t;\mu^t)+\la \r^t, \x^{t+1}- \x^{t} \ra + \tfrac{1}{2}\|\x^t-\x^{t+1}\|_{[ \A\trans \A /\mu^t + \tilde{\M} ]}^2.\nn
\eeq
\noi We observe that the following equality holds:
\beq
\mathcal{R}(\x^{t+1},\y^t;\mu^t) -\mathcal{R}(\x^t,\y^t;\mu^t) = \mathcal{J}(\x^{t+1},\y^{t};\mu^t) - \mathcal{J}(\x^{t},\y^{t};\mu^t).\nn
\eeq
\noi Summing up these three inequalities, we obtain:
\beq \label{eq:opt:x:2}
\E_{\xi^t}[\mathcal{J}(\x^{t+1},\y^{t};\mu^t) - \mathcal{J}(\x^{t},\y^{t};\mu^t)] &\leq& \E_{\xi^t}[-(\tfrac{\theta_1}{\mu^t}+\theta_2) \cdot \frac{1}{2}\|\x^{t+1}-\x^{t}\|_2^2]\nn\\
&\overset{\step{172}}{\leq} & \E_{\xi^t}[-\theta\cdot \frac{1}{2}\|\x^{t+1}-\x^{t}\|_2^2],
\eeq
\noi where step \step{172} uses $\frac{\theta_1}{\mu^t} + \theta_2 \leq \frac{\theta_1}{\mu^1} + \theta_2 \triangleq \theta$ as the sequence $\{\mu^t\}_{t=1}^{\infty}$ is non-increasing.

\noi We now focus on the $\y$-subproblem. Similar to the proof for Theorem \ref{theorem:convergence:iht}, we have:
\beq \label{eq:opt:y:2}
 \mathcal{J}(\x^{t+1},\y^{t+1};\mu^{t}) - \mathcal{J}(\x^{t+1},\y^{t};\mu^{t})     \leq  - \frac{1}{2\mu^1}\|\y^{t+1}-\y^t\|_2^2.
\eeq
\noi Using the continuity of $\mathcal{J}(\x^{t},\y^{t};\mu)$ \textit{w.r.t.} $\mu$ as detailed in \textbf{Part (e)} of Lemma \ref{lemma:first:bound}, we obtain:
\beq\label{eq:connut:mu:2}
0\leq \E_{\xi^t}[\mathcal{J}(\x^{t},\y^{t};\mu^{t}) -   \mathcal{J}(\x^{t},\y^{t};\mu^{t-1})] \leq \tfrac{L_h^2}{2}   (\tfrac{(\mu^{t-1})^2}{\mu^{t}} - {\mu^{t-1}} ) \triangleq  \Psi^t.
\eeq

\noi Summing up Inequalities (\ref{eq:opt:x:2}), (\ref{eq:opt:y:2}), and (\ref{eq:connut:mu:2}) together, we have:
\beq \label{eq:sss:ddd:ccc}
&& \E_{\xi^T}[\frac{1}{2\mu^1}\|\y^{t+1}-\y^t\|_2^2 +  \tfrac{\theta}{2}\|\x^{t+1}-\x^{t}\|_2^2]  \nn\\
&\leq&  \E_{\xi^{t-1}}[\mathcal{J}(\x^{t},\y^{t};\mu^{t-1}) ] - \E_{\xi^{t}}[ \mathcal{J}(\x^{t+1},\y^{t+1};\mu^{t}) ] + \tfrac{1}{2}  L_h^2 (\tfrac{(\mu^{t-1})^2}{\mu^{t}} - {\mu^{t-1}} ) \nn \\
&=&\mathcal{J}^t - \mathcal{J}^{t+1} + \Psi^t,
 \eeq
\noi where $\mathcal{J}^{t+1} \triangleq \E_{\xi^t}[\mathcal{J}(\x^{t+1},\y^{t+1};\mu^{t})]$.

\noi \textbf{(b)} Summing up the inequality in (\ref{eq:sss:ddd:ccc}) over $t=1,2,...,T$, we have:
\beq
&& \textstyle \sum_{t=1}^T \tfrac{1}{2\mu^1}\|\y^{t+1}-\y^{t}\|_2^2 + \sum_{t=1}^T \tfrac{\theta}{2}\|\x^{t+1}-\x^{t}\|_2^2 \nn\\
& \leq & \textstyle   \mathcal{J}^1 - \mathcal{J}^{T+1} + [\sum_{t=1}^T \Psi^t ]   \overset{\step{172}}{\leq}  \textstyle \mathcal{J}^1 - \mathcal{J}^{T+1} + \eta L_h^2  =  C < + \infty,\nn
\eeq
\noi where step \step{172} uses $\sum_{t=1}^T \Psi^t < \sum_{t=1}^{\infty} \Psi^t \leq \eta L_h^2 $, as demonstrated in \textbf{Part (f)} of Lemma \ref{lemma:first:bound}.

\noi \textbf{(c)} As a result, there exists an index $\bar{t}$ with $1 \leq \bar{t} \leq T$ such that: $\tfrac{1}{2\mu^1}\|\y^{\bar{t}+1}-\y^{\bar{t}}\|_2^2 + \tfrac{\theta}{2}\|\x^{\bar{t}+1}-\x^{\bar{t}}\|_2^2 \leq \frac{C}{T}$, leading to:
\beq
\|\y^{\bar{t}+1}-\y^{\bar{t}}\|_2^2 +  \|\x^{\bar{t}+1}-\x^{\bar{t}}\|_2^2 \leq \frac{2C}{T \cdot \min(\theta, (\mu^1)^{-1}) }. \label{eq:lll:bcd:1}
\eeq
\noi We define $\Gamma_x(\x,\y;\mu) \triangleq \frac{1}{C_n^k} \sum_{\B \in \Omega_n^k } \dist^2(\x_{\B},\arg\min_{ \z_{\B}}\,\delta(\U_{\B}\z_{\B} + \U_{\B^c}\x_{\B^c})+\ddot{\mathcal{M}}(\U_{\B}\z_{\B} + \U_{\B^c}\x_{\B^c},\x,\y;\mu)$ and $\Gamma_y(\x,\y;\mu) \triangleq \dist^2(\y, \arg \min_{\y'}\mathcal{J}(\x,\y';\mu))$. It is important to note that $\x^{\bar{t}+1}$ and $\x^{\bar{t}}$ differ in at most $k$ coordinates. We have:
\beq
\|\x^{\bar{t}+1}-\x^{\bar{t}}\|_2^2 + \|\y^{\bar{t}+1}-\y^{\bar{t}}\|_2^2 \geq
\Gamma_x(\x^{\bar{t}},\y^{\bar{t}};\mu) + \Gamma_y(\x^{\bar{t}},\y^{\bar{t}};\mu) \label{eq:lll:bcd:2}
\eeq
\noi for all $\bar{t}\geq 1$ and some sufficiently small $\mu=\mu^{\bar{t}}>0$. Combining Inequality (\ref{eq:lll:bcd:1}) and Inequality (\ref{eq:lll:bcd:2}), we have:
\beq
\Gamma_x(\x^{\bar{t}},\y^{\bar{t}};\mu^{\bar{t}}) + \Gamma_y(\x^{\bar{t}},\y^{\bar{t}};\mu^{\bar{t}}) \leq \frac{2C}{T\cdot \min(\theta,(\mu^1)^{-1})}. \nn
\eeq
\noi Therefore, we conclude that Algorithm \ref{algo:main} finds an $\epsilon$-approximate block-$k$ stationary point of Problem (\ref{eq:main}) in at most $T$ iterations in the sense of expectation, where $T\leq \lceil  \frac{2C}{\epsilon \min(\theta,(\mu^1)^{-1}) }\rceil = \mathcal{O}(\epsilon^{-1})$.

\end{proof}

\subsection{Proof of Lemma \ref{lemma:bcd:bound}}
\label{app:lemma:bcd:bound}

 \begin{proof}

We denote $\B=\B^t$. We define $\H^t\triangleq (\A\trans \A + \theta_1 \I_n)/ {\mu^t} + \tilde{\M} + \theta_2 \I_n \in \mathbb{R}^{n\times n}$, and $\H_{\ast}^t \triangleq \UBt \UB \trans \H^t \UB \UBt\trans \in \mathbb{R}^{n\times n}$. $\Sup^t \triangleq \Vup+\theta+\frac{\theta_1+\Aup}{\mu^t}$, $\Slow^t \triangleq\Vlow+\theta+\frac{\theta_1+\Alow}{\mu^t}$.

\noi \textbf{(a)} Problem (\ref{eq:subprob:x:2}) in Algorithm \ref{algo:main} is equivalent to solving the following optimization problem:
\beq
\x_{\B}^{t+1}  \in \arg \min_{\z_{\B}\in\mathbb{R}^k}\,\mathcal{W}(\z_{\B})\,s.t.\,\|\z_{\B}\|_0 + \|\x^t_{\B^c}\|_0\leq s, \label{eq:dec:subp}
\eeq
\noi where $\mathcal{W}(\z_{\B}) \triangleq \la \z_{\B} - \x_{\B}^t,\r^t_{\B} \ra + \frac{1}{2}(\z_{\B} - \x_{\B}^t)\trans [\H^t]_{\B\B} (\z_{\B} - \x_{\B}^t)$. By the optimality of $\x_{\B}^{t+1}$, we have: $\mathcal{W}(\x^{t+1}_{\B})\leq \mathcal{W}(\x^t_{\B})=0$, leading to:
\beq \label{eq:opt:X}
\la \x_{\B}^{t+1} - \x^t_{\B}, \r^t_{\B}\ra + \frac{1}{2}(\x^{t+1}_{\B} - \x_{\B}^t)\trans [\H^t]_{\B\B} (\x^{t+1}_{\B} - \x_{\B}^t ) \leq 0.
\eeq
\noi We derive the following inequalities:
\beq
 \frac{1}{2} \Slow^t \|\x^{t+1}_{\B} - \x_{\B}^t\|_2^2  \overset{\step{172}}{\leq} \frac{1}{2} (\x^{t+1}_{\B} - \x_{\B}^t) \trans [\H^t]_{\B\B} (\x^{t+1}_{\B} - \x_{\B}^t) \overset{\step{173}}{\leq} -\la \x_{\B}^{t+1} - \x^t_{\B}, \r^t_{\B}\ra  \overset{\step{174}}{=} \|\x_{\B}^{t+1} - \x^t_{\B}\| \|\r^t_{\B}\|, \nn
\eeq
\noi where step \step{172} uses $\Slow^t\I_k \preceq [\H^t]_{\B\B} \preceq \Sup^t \I_k$; step \step{173} uses (\ref{eq:opt:X}); step \step{174} uses the Cauchy-Schwarz inequality. Dividing both sides by $\|\x^{t+1}_{\B} - \x_{\B}^t\|$, we have: $\E_{\xi^t} [\|\x^{t+1}_{\B} - \x_{\B}^t\|] \leq  \E_{\xi^t} [\frac{2}{\Slow^t}\|\r^t_{\B}\|]$. Using the result in Lemma \ref{lemma:bcd:simple:basic}, we have:
\beq \label{eq:xt1:xt}
Z_n^k\E_{\xi^t} [\|\x^{t+1} - \x^t\|] \leq  Z_n^k \E_{\xi^t} [\tfrac{2}{\Slow^t}\|\r^t\|].
\eeq

\noi \textbf{(b)} For notation convenience, we define:
\beq \label{eq:def:B1:B2}
\B_1 \triangleq \{i\,|\,\x_i^{t+1}\neq 0,\,i \in \B\},\,\text{and}\,\B_2 \triangleq \{i\,|\,\x_i^{t+1}=0,\,i \in \B\}.\nn
\eeq
\noi The solution $\x_{\B}^{t+1} \in \mathbb{R}^k$ is a local minimimizer for Problem (\ref{eq:dec:subp}) if and only if $[\nabla \mathcal{W}(\x_{\B}^{t+1})]_{\B_1}=0$. Using the optimality condition for $\x^{t+1}_{\B}$, we have:
\beq \label{eq:optiaml:condition:x}
\zero  &=& [\nabla \mathcal{W}(\z)]_{\x_{\B}^{t+1}} = \r^t_{\B_1} +  [[\H^t_{\B\B}]_{\B_1\B}] (\x^{t+1}_{\B}-\x^{t}_{\B}) \nn\\
&\overset{\step{172}}{=} & \r^t_{\B_1} +    [[\H_{\B\B}^t]_{\B_1\B_1} ] (\x^{t+1}_{\B_1}-\x^{t}_{\B_1}) +  [[\H_{\B\B}^t]_{\B_1\B_2} ] (\x^{t+1}_{\B_2}-\x^{t}_{\B_2}),
\eeq
\noi where step \step{172} uses $\B = \B_1 \cup \B_2$. We derive the following equalities:
\beq
&&\E_{\xi^t}[ \la [[\H^t]_{\B\B}]( \x_{\B}^{t+1} - \x_{\B}^{t}), \x_{\B}^{t+1} \ra ] \nn\\
&=&    \E_{\xi^t}[ \begin{bmatrix}
\x^{t+1}_{\B_1} - \x^{t}_{\B_1}   \\
  \x^{t+1}_{\B_2} - \x^{t}_{\B_2} \\
\end{bmatrix} \trans  \begin{bmatrix}
[\H^t_{\B\B}]_{\B_1\B_1}   & [\H^t_{\B\B}]_{\B_1\B_2}   \\
[\H^t_{\B\B}]_{\B_2\B_1}  & [\H^t_{\B\B}]_{\B_2\B_2}  \\
\end{bmatrix}  \begin{bmatrix}
\x^{t+1}_{\B_1}  \\
\x^{t+1}_{\B_2} \\
\end{bmatrix}   ]\nn\\
&\overset{\step{172}}{=}&  \E_{\xi^t}[\,\la [[\H^t_{\B\B}]_{\B_1\B_1}] (\x^{t+1}_{\B_1}-\x^{t}_{\B_1}   ) , \x^{t+1}_{\B_1} \ra +  \la [[\H^t_{\B\B}]_{\B_1\B_2}] (\x^{t+1}_{\B_2} - \x^{t}_{\B_2}), \x^{t+1}_{\B_1}  \ra] \nn\\
&\overset{\step{173}}{=}&  \E_{\xi^t}[\, \la - \r^t_{\B_1} - [[\H_{\B\B}^t]_{\B_1\B_2} ] (\x^{t+1}_{\B_2} -\x^{t}_{\B_2}) , \x_{\B_1}^{t+1}\ra + \la [[\H^t_{\B\B}]_{\B_1\B_2}]  (\x^{t+1}_{\B_2} - \x^{t}_{\B_2}), \x^{t+1}_{\B_1} \ra   + 0 + 0 ] \nn\\
&\overset{}{=}&  \E_{\xi^t}[\, \la - \r^t_{\B_1} , \x_{\B_1}^{t+1}\ra ] \nn\\
&\overset{\step{174}}{=}& \E_{\xi^t}[- \la \r_{\B_1}^t, \x^{t+1}_{\B_1}\ra -  \la \r_{\B_2}^t, \x^{t+1}_{\B_2}\ra  ] = \E[ - \la \r_{\B}^t, \x^{t+1}_{\B}\ra ] \nn\\
&\overset{\step{175}}{=}& -Z_n^k  \la \r^t, \x^{t+1}\ra, \label{eq:Bound:222}
\eeq
\noi where step \step{172} uses the fact that $[\x^{t+1}]_{\B_2} = \mathbf{0}$; step \ding{173} uses the optimality condition as in (\ref{eq:optiaml:condition:x}); step \step{174} uses $\B=[\B_1;\B_2]$ and the fact that $[\x^{t+1}]_{\B_2} = \mathbf{0}$; step \step{175} uses Lemma \ref{lemma:bcd:simple:basic} with $ Z_n^k =\frac{k}{n}$.

\noi \textbf{(c)} We derive the following equalities:
\beq
&& \E_{\xi^t}[ \frac{1}{2} \| \x^{t+1} - \bar{\x}\|_{\H_{\ast}^t}^2] - \E_{\xi^t}[\frac{1}{2} \| \x^{t} - \bar{\x}\|_{\H_{\ast}^t}^2]  +  \E_{\xi^t}[ \frac{1}{2} \| \x^{t+1} - \x^{t}\|_{\H_{\ast}^t}^2]  \nn\\
&\overset{\step{172}}{=}&  \E_{\xi^t}[ \la \H_{\ast}^t( \x^{t+1} - \x^{t}), \x^{t+1}-\bar{\x} \ra ]\nn\\
&\overset{\step{173}}{=}&  \E_{\xi^t}[ \la  [[\H^t]_{\B\B}]( \x_{\B}^{t+1} - \x_{\B}^{t}), \x_{\B}^{t+1} -\bar{\x}_{\B} \ra ]\nn\\
&=& \underbrace{\E_{\xi^t}[ \la  [[\H^t]_{\B\B}]( \x_{\B}^{t} - \x_{\B}^{t+1}), \bar{\x}_{\B} \ra ]}_{\Gamma_1} + \underbrace{\E_{\xi^t}[ \la [ [\H^t]_{\B\B} ]( \x_{\B}^{t+1} - \x_{\B}^{t}), \x_{\B}^{t+1} \ra ]}_{\Gamma_2},  \label{eq:KKKK}
\eeq
\noi where step \step{172} uses the Pythagoras relation; step \step{173} uses $[\H_{\ast}^t]_{\B\B}=[\H^t]_{\B\B}$ and $\x_{\Bc}^{t+1}-\x_{\Bc}^t=\zero$.

\noi We first bound the term $\Gamma_1$ in (\ref{eq:KKKK}) using the following inequalities:
\beq
\label{gamma:2222}
\Gamma_1 &=&\E_{\xi^t}[ \la   [ [\H^t]_{\B\B} ] ( \x^{t}_{\B} - \x^{t+1}_{\B}), \bar{\x}_{\B}\ra] \nn\\
&\overset{\step{172}}{\leq} & \E_{\xi^t}[  \|  [[\H^t]_{\B\B}] ( \x^{t}_{\B} - \x^{t+1}_{\B}) \|  \cdot \|\bar{\x}_{\B}\|]\nn\\
&\overset{\step{173}}{\leq} &\E_{\xi^t}[  \Sup^t \|\x^{t}_{\B} - \x^{t+1}_{\B}\|  \cdot \|\bar{\x}_{\B}\|]\nn\\
&\overset{\step{174}}{=} &   Z_n^k  \Sup^t \|\x^{t}- \x^{t+1}\|  \cdot \|\bar{\x}\| \nn\\
&\overset{\step{175}}{\leq} &  2  Z_n^k    \|\r^t\| \|\bar{\x}\| \Sup^t/\Slow^t = 2  Z_n^k \|\r^t\| \|\bar{\x}\|\kappa^t,
\eeq
\noi where step \step{172} uses the Cauchy-Schwarz inequality; step \step{173} uses $\Slow^t\I_k \preceq [\H^t]_{\B\B} \preceq \Sup^t \I_k$; step \step{174} uses Lemma \ref{lemma:bcd:simple:basic} with $ Z_n^k =\frac{k}{n}$; step \step{175} uses Inequality (\ref{eq:xt1:xt}).

\noi We now bound the term $\Gamma_2$ in (\ref{eq:KKKK}) using the following inequalities:
\beq\label{gamma:1111}
\Gamma_2 &\overset{\step{172}}{\leq}& -Z_n^k  \la \r^t, \x^{t+1} \ra \nn\\
&=& Z_n^k  \la \r^t, \bar{\x}-\x^t\ra  +  Z_n^k  \la \r^t, \x^t - \x^{t+1}\ra  +  Z_n^k  \la \r^t,-\bar{\x}\ra  \nn\\
&\overset{\step{173}}{\leq}&Z_n^k  (\Upsilon^t -\tfrac{V_s}{2}\| \bar{\x}-\x^t \|_2^2)   +  Z_n^k  \|\r^t\| \|\x^t - \x^{t+1}\|  +  Z_n^k  \|\r^t\| \|\bar{\x}\|  \nn\\
&\overset{\step{174}}{\leq}&Z_n^k (\Upsilon^t-\tfrac{V_s}{2}\| \bar{\x}-\x^t \|_2^2)+\tfrac{2  Z_n^k }{\Slow^t} \|\r^t\|_2^2+Z_n^k\|\r^t\|\|\bar{\x}\|,
\eeq
\noi where step \step{172} uses Equality (\ref{eq:Bound:222}); step \step{173} uses Lemma \ref{lemma:bound:negative:l2} that $\la \r^t, \bar{\x} - \x^t \ra\leq \Upsilon^t  -\frac{V_s}{2}\|\x^t -\bar{\x}\|_2^2$, and the Cauchy-Schwarz inequality; step \step{174} uses (\ref{eq:xt1:xt}).

\noi In view of (\ref{eq:KKKK}), (\ref{gamma:2222}), and (\ref{gamma:1111}), we have:
\beq
&&\E_{\xi^t}[ \tfrac{1}{2} \| \x^{t+1} - \bar{\x}\|_{\H_{\ast}^t}^2 -  \tfrac{1}{2} \| \x^{t} - \bar{\x}\|_{\H_{\ast}^t}^2  ]  \nn\\
&\leq &  - \E_{\xi^t}[ \tfrac{1}{2} \| \x^{t+1} - \x^{t}\|_{\H_{\ast}^t}^2]  + \Gamma_1 + \Gamma_2 \nn\\
&\overset{\step{172}}{\leq} &  0 +   Z_n^k  (\Upsilon^t-\tfrac{V_s}{2}\| \bar{\x}-\x^t \|_2^2 )   +  \tfrac{2  Z_n^k }{ \Slow^t } \|\r^t\|_2^2+(1+2\kappa^t)Z_n^k\|\r^t\| \|\bar{\x}\|, \nn
\eeq
\noi where step \step{172} uses $- \E_{\xi^t}[ \tfrac{1}{2} \| \x^{t+1} - \x^{t}\|_{\H_{\ast}^t}^2] \leq 0$.

\end{proof}

\subsection{Proof of Lemma \ref{lemma:condition:constant}}
\label{app:lemma:condition:constant}

\begin{proof}
We initially establish the subsequent inequality:
\beq \label{eq:abcd}
\frac{a+b}{c+d} \leq \max( \frac{a}{c},\frac{b}{d}), \forall a\geq 0, b\geq 0, c>0, d>0.
\eeq
We consider two cases: \bfit{(i)} $\frac{a}{c}\geq \frac{b}{d}$. we have: $b\leq \frac{ad}{c}$, leading to $\frac{a+b}{c+d}\leq \frac{a+\frac{ad}{c}}{c+d} = \frac{a}{c}\cdot \frac{c+d}{c+d}=\frac{a}{c}$. \bfit{(ii)} $\frac{a}{c}<\frac{b}{d}$. We have: $a \leq \frac{bc}{d}$, resulting in $\frac{a+b}{c+d}\leq \frac{\frac{bc}{d}+b}{c+d}=\frac{b}{d}\cdot \frac{c+d}{c+d}=\frac{b}{d}$. Therefore, Inequality (\ref{eq:abcd}) holds.

Using the definition of $\Sup^t$ and $\Slow^t$, we have:
\beq
\frac{\Sup^t}{\Slow^t} = \frac{ \frac{\Aup+\theta_1}{\mu^t}+ \Vup+\theta_2 }{  \frac{\Alow+\theta_1}{\mu^t} + \Vlow+\theta_2 } \overset{\step{172}}{\leq}\max(\frac{\Aup+\theta_1}{\Alow+\theta_1},\frac{\Vup+\theta_2}{ \Vlow+\theta_2}) \overset{\step{173}}{\leq}  1 + \epsilon,\nn
\eeq
\noi where step \step{172} uses Inequality (\ref{eq:abcd}); step \step{173} uses the fact that $\frac{\Aup+\theta_1}{\Alow+\theta_1} \leq 1+\epsilon$ if $\theta_1\geq \frac{\Aup - \Alow (1+\epsilon)}{\epsilon}$, and $\frac{\Vup+\theta_2}{ \Vlow+\theta_2} \leq 1+\epsilon$ if $\theta_2 \geq \frac{\Vup - \Vlow (1+\epsilon)}{\epsilon}$.

 \end{proof}

\subsection{Proof of Theorem \ref{theorem:DEC:stepsize:constant}} \label{app:theorem:DEC:stepsize:constant}

\begin{proof}

We consider constant stepsizes with $\mu^t = \bar{\mu}$ for all $t\geq 1$.

We define: $\H_{\ast}^t \triangleq \UBt \UBt\trans \H \UBt  \UBt\trans \in \mathbb{R}^{n\times n}$, $\H \triangleq (\A\trans \A +\theta_1 \I_n )/\bar{\mu} + \tilde{\M} + \theta_2 \I_n \in \mathbb{R}^{n\times n}$, $\Slow \triangleq \frac{\Alow+\theta_1}{\bar{\mu}} + \Vlow+\theta_2$, and $\Sup \triangleq \frac{\Aup+\theta_1}{\bar{\mu}} + \Vup+\theta_2$, $\Delta_{\x}^t \triangleq \E_{\xi^{t}} [\|\x^t - \bar{\x}\|_2^2]$, $\Delta_{F}^t \triangleq \E_{\xi^{t}} [(\min_{i=1}^t F(\x^i) )  -   F(\bar{\x})]$.

\noi First, using \textbf{Part (c)} in Lemma \ref{lemma:first:bound}, we have: $\|\r^t\|\leq L_F$.

\noi Second, using \textbf{Part (b)} of Lemma \ref{lemma:bound:negative:l2}, we have the upper bound of $\Upsilon^t$ that: $\forall t,\,\Upsilon^t \leq \tfrac{1}{2} \bar{\mu}L_h^2 - \Delta_{F}^t$.

\noi Third, with $\B^t$ and $\B^{t+1}$ randomly and uniformly chosen, for any $\z\in\mathbb{R}^n$, the following holds:
\beq
\E_{\B^t} [\|\z\|_{\H_{\ast}^t}^2] &=& \E_{\B^t} [\z\trans [ \H_{\ast}^t]\z] = \E_{\B^t} [\z\trans  \UBt \UBt\trans \H \UBt  \UBt\trans ]\z \nn\\
&=& \E_{\B^{t+1}} [\z\trans [ \UBtt \UBtt\trans \H \UBtt  \UBtt\trans ]\z ] \nn\\
&= & \E_{\B^{t+1}} [\|\z\|_{\H_{\ast}^{t+1}}^2] .\label{eq:Ht:Ht1}
\eeq


\noi \textbf{(a)} Building upon our prior discussions, we derive the following inequalities:
\beq \label{eq:based:based}
&& \E_{\xi^{t+1}}[  \tfrac{1}{2}\| \x^{t+1} - \bar{\x}\|_{\H_{\ast}^{t+1}}^2 ]      \nn\\
&\overset{\step{172}}{=}& \E_{\xi^{t}}[ \tfrac{1}{2}\| \x^{t+1} - \bar{\x}\|_{\H_{\ast}^t}^2]    \nn\\
&\overset{\step{173}}{\leq}&  \E_{\xi^{t}}[\tfrac{1}{2}\| \x^{t} - \bar{\x}\|_{\H_{\ast}^t}^2] - Z_n^k \tfrac{V_s}{2}\| \x^{t} - \bar{\x} \|_2^2 +  Z_n^k \{ \tfrac{2}{\Slow^t}\|\r^t\|_2^2 + \Upsilon^t +    ( 1+    2 \kappa^t )       \|\r^t\| \|\bar{\x}\|\} \nn\\
&\overset{\step{174}}{\leq}& \E_{\xi^{t}}[(1-\tfrac{V_s}{ \Sup})\cdot\tfrac{1}{2}\|\x^{t}-\bar{\x}\|_{\H_{\ast}^t}^2]  +  Z_n^k \{  \tfrac{2}{\Slow^t} (L_F)^2 + \tfrac{\bar{\mu}L_h^2}{2} -  \Delta_{F}^t +    ( 3 + 2 \epsilon) L_F \|\bar{\x}\| \} ,
\eeq
\noi where step \step{172} uses Equality (\ref{eq:Ht:Ht1}) with $\z=\x^{t+1} - \bar{\x}$, leading to $\E_{\xi^t}[\| \x^{t+1} - \bar{\x}\|_{\H_{\ast}^t}^2 ] = \E_{\xi^{t+1}}[\| \x^{t+1} - \bar{\x}\|_{\H_{\ast}^{t+1}}^2 ]$; step \step{173} uses the inequality in \textbf{Part (b)} of Lemma (\ref{lemma:bcd:bound}); step \step{174} uses $Z_n^k\frac{V_s }{2}\|\x^t-\bar{\x}\|_2^2 \geq  \E_{\xi^t}[\frac{V_s } {2 \Sup } \|\x^t -  \bar{\x}\|_{\H_{\ast}^t}^2]$.

\noi \textbf{(b)} Based on (\ref{eq:based:based}), we apply Lemma \ref{lemma:two:non:negative:sequences} with the following definitions:
\beq
\gamma \triangleq 1 - \tfrac{V_s}{\Sup},\,\Phi^{t} \triangleq \E_{\xi^{t}}[ \tfrac{1}{2}\| \x^{t} - \bar{\x}\|_{\H_{\ast}^t}^2], \,\Lambda^t\triangleq  Z_n^k \{  \tfrac{2}{\Slow^t} (L_F)^2+\tfrac{\bar{\mu}}{2}L_h^2 - \Delta_{F}^t +  ( 3 + 2 \epsilon)L_F \|\bar{\x}\| \}.
\eeq
\noi This results in the subsequent inequality for any integer $T\geq 1$:
\beq
&&\E_{\xi^{T+1}}[ \tfrac{1}{2}\| \x^{T+1} - \bar{\x}\|_{\H_{\ast}^t}^2] \nn\\
&\leq&    \gamma^T  \E_{\xi^{1}}[ \frac{1}{2} \| \x^{1} - \bar{\x}\|_{\H_{\ast}^1}^2 ] + \tfrac{ Z_n^k }{ 1 - \gamma} \max_{t=1}^T  \{ \tfrac{2}{\Slow} (L_F)^2 + \tfrac{\bar{\mu}}{2}L_h^2 - \Delta_{F}^t +   ( 3 + 2 \epsilon)      L_F \|\bar{\x}\| \} \label{BCD:fix:stepsize}
\eeq
\noi We further derive the following inequalities:
\beq
&&\E_{\xi^{T+1}} [\tfrac{1}{2}\|\x^{T+1} - \bar{\x}\|_2^2] \nn\\
&\overset{\step{172}}{\leq}& \tfrac{1}{\Slow  Z_n^k } \E[\tfrac{1}{2}\|\x^{T+1} - \bar{\x}\|_{\H_{\ast}^{T+1} }^2] \nn\\
&\overset{\step{173}}{\leq}& \tfrac{\gamma^T}{\Slow  Z_n^k }  \cdot \E_{\xi^{1}}[ \tfrac{1}{2} \| \x^{1} - \bar{\x}\|_{\H_{\ast}^1}^2 ] + \tfrac{1}{\Slow (1 - \gamma) }   \max_{t=1}^T  \{\tfrac{2 (L_F)^2 }{\Slow^t}  +  \tfrac{1}{2}\bar{\mu}L_h^2 -   \Delta_{F}^t  + ( 3 + 2 \epsilon)      L_F \|\bar{\x}\| \}    \nn\\
&\overset{\step{174}}{\leq}& \tfrac{\gamma^T   \Sup}{\Slow   }   \cdot       \E_{\xi^1}[\tfrac{1}{2}\| \x^{1} - \bar{\x}\|_{2}^2 ]   + \tfrac{1}{\Slow  (1-\gamma) } \cdot   \max_{t=1}^T  \{\tfrac{  2(L_F)^2 }{\Slow^t}  +  \tfrac{1}{2}\bar{\mu}L_h^2 -   \Delta_{F}^t +  ( 3 + 2 \epsilon)      L_F \|\bar{\x}\| \}      \nn\\
&\overset{\step{175}}{\leq}&\gamma^T \cdot (1+\epsilon) \cdot \E_{\xi^1}[\tfrac{1}{2}\|\x^1 - \bar{\x} \|_2^2] + \tfrac{1}{V_s} \{ 2(L_F)^2 \cdot  \tfrac{\bar{\mu}}{\Aup+\theta_1} +  \tfrac{1}{2}\bar{\mu}L_h^2 -   \Delta_{F}^T   + ( 3 + 2 \epsilon) L_F \|\bar{\x}\| \} \nn\\
&\overset{\step{176}}{=}& \tfrac{1}{V_s} \{ K_3 \gamma^T +  D_3 \bar{\mu} - \Delta_{F}^T + ( 3 + 2 \epsilon) L_F \|\bar{\x}\|\} ,  \label{eq:BCD:final:inequality}
\eeq
\noi where step \step{172} uses $\Slow  Z_n^k \E_{\xi^{T+1}}[ \|\x^{T+1} - \bar{\x}\|_2^2] \leq  \E_{\xi^{T+1}}[ \|\x^{T+1} - \bar{\x}\|_{ \H_{\ast}^{T+1}}^2]$; step \step{173} uses Inequality (\ref{BCD:fix:stepsize}); step \step{174} uses $\E_{\xi^{1}}[ \| \x^{1} - \bar{\x}\|_{\H_{\ast}^1}^2 ] \leq  Z_n^k  \Sup  \E_{\xi^{1}}[ \| \x^{1} - \bar{\x}\|_{2}^2 ] $; step \step{175} uses $\frac{\Sup}{\Slow}=\kappa \leq 1+\epsilon$, $\gamma \triangleq 1 - \frac{V_s}{\Sup}$, and $\Sup \triangleq \frac{\Alow+\theta_1}{\bar{\mu}} + \Vlow+\theta_2 \geq \frac{\Aup+\theta_1}{\bar{\mu}}$; step \step{176} uses $K_3\triangleq \E_{\xi^1}[(1+\epsilon)\cdot\tfrac{V_s}{2}\|\x^1 - \bar{\x} \|_2^2]  = \tfrac{V_s (1+\epsilon) }{2} \Delta_{\x}^1$ and $D_3\triangleq \frac{2(L_F)^2 }{  \theta_1 + \Aup  } + \frac{L_h^2}{2}$.

\noi We now focus on (\ref{eq:BCD:final:inequality}). Using the fact that $\E_{\xi^{T+1}} [\frac{1}{2}\|\x^{T+1} - \bar{\x}\|_2^2] \geq 0$, we obtain: $\Delta_{F}^T  \leq  K_3 \gamma^T    + D_3 \bar{\mu} + ( 3 + 2 \epsilon) L_F \|\bar{\x}\|$.

\noi Using the fact that $\Delta_{F}^T \geq 0$, we obtain: $\Delta_{\x}^{T+1}   \leq \frac{2}{V_s} \left( K_3 \gamma^T +  D_3 \bar{\mu}  + ( 3 + 2 \epsilon) L_F \|\bar{\x}\|\right) $.

\end{proof}

\subsection{Proof of Theorem \ref{theorem:DEC:stepsize:diminishing}} \label{app:theorem:DEC:stepsize:diminishing}

To finish the proof of this theorem, we first provide the following useful lemma.

\begin{lemma} \label{lemma:BCD:Ht1:Ht}
Assume $\mu^t = \frac{\eta}{t+t_0}$ with $\eta =  \frac{\Aup + \theta_1}{V_s}$. We have: $\E_{\xi^{t+1}}[ \frac{1}{2} \| \x^{t+1} - \bar{\x}\|_{\H_{\ast}^{t+1}}^2]  \leq \E_{\xi^t}[ \frac{1}{2} \| \x^{t+1} - \bar{\x}\|_{\H_{\ast}^{t}}^2] + \tfrac{V_s Z_n^k }{2} \|\x^{t+1} - \bar{\x} \|_2^2 $.

 \begin{proof}

We denote $\H_{\ast}^t \triangleq \UBt \UBt \trans \H^t \UBt \UBt\trans \in \mathbb{R}^{n\times n}$, where $\H^t\triangleq (\A\trans \A + \theta_1 \I_n)/ {\mu^t} + \tilde{\M} + \theta_2 \I_n \in \mathbb{R}^{n\times n}$.

\noi We have the following inequalities for all $\z \triangleq \x^{t+1} - \bar{\x} \in \mathbb{R}^n$:
\beq
&&\E_{\xi^{t+1}}[ \|\z \|_{\H_{\ast}^{t+1}}^2]  -  \E_{\xi^t}[ \| \z \|_{\H_{\ast}^{t}}^2] \nn\\
&\overset{\step{172}}{=}& \E_{\xi^{t+1}}[  \z \trans  ( \UBtt \UBtt\trans \H^{t+1} \UBtt \UBtt \trans ) \z  ]  -  \E_{\xi^t}[  \z \trans  ( \UBt \UBt\trans \H^t \UBt \UBt\trans ) \z  ]  \nn\\
&\overset{\step{173}}{=}&\E_{\xi^{t}}[  \z \trans  ( \UBt \UBt \trans \H^{t+1} \UBt \UBt \trans ) \z  ]  -  \E_{\xi^t}[  \z \trans  ( \UBt \UBt \trans \H^{t} \UBt \UBt \trans ) \z  ]  \nn\\
&\overset{}{= }& \E_{\xi^{t}}[  \z \trans  ( \UBt \UBt\trans [\H^{t+1}-\H^{t}] \UBt \UBt \trans ) \z  ]  \nn\\
&\overset{\step{174}}{=}& (\tfrac{1}{\mu^{t+1}} - \tfrac{1}{\mu^{t}}) \E_{\xi^{t}}[  \z \trans  ( \UBt [ \A\trans \A + \theta_1 \I_n]_{\B^t\B^t} \UBt \trans ) \z  ]  \nn\\
&\overset{\step{175}}{\leq }& \tfrac{1}{\eta} \E_{\xi^t}[  \z \trans \UBt \left(    \Aup \I_k + \theta_1 \I_k   \right)  \UBt\trans \z ] = \tfrac{\Aup + \theta_1}{\eta} \E_{\B^t}[  \z \trans \UBt \UBt\trans \z ] \nn\\
&\overset{\step{176}}{ = }& Z_n^k \tfrac{\Aup + \theta_1}{\eta} \|\z\|_2^2 \nn\\
&\overset{\step{177}}{ = }& Z_n^k V_s \|\z\|_2^2, \nn
\eeq
\noi where step \step{172} uses the definition of $\H_{\ast}^t$; step \step{173} uses the fact that both $\B^t$ and $\B^{t+1}$ are choosen randomly and uniformly; step \step{174} uses the choice $\mu^t=\frac{\eta}{t+t_0}$ that $\frac{1}{\mu^{t+1}} - \frac{1}{\mu^t}  = \frac{1}{\eta} \cdot \left(   (t+t_0+1) - (t + t_0)  \right) =   \frac{1}{\eta} $; step \step{175} uses $[\A\trans \A]_{\B^t\B^t} \preceq \Aup \I_k$; step \step{176} uses Lemma \ref{lemma:bcd:simple:basic}; step \step{177} uses the choice $\eta=\tfrac{\Vup+\theta_1}{V_s}$.

\end{proof}

\end{lemma}

We now prove the proof of this theorem.

\begin{proof}

We consider diminishing stepsizes with $\mu^t =   \tfrac{\eta}{t+t_0}$ for all $t\geq 1$, where $\eta =  \frac{\Aup + \theta_1}{V_s}$.

We define: $\H_{\ast}^t \triangleq \UBt \UBt\trans \H^t \UBt  \UBt\trans \in \mathbb{R}^{n\times n}$, $\H^t \triangleq (\A\trans \A +\theta_1 \I_n )/\mu^t + \tilde{\M} + \theta_2 \I_n  \in \mathbb{R}^{n\times n}$, $\Slow \triangleq \frac{\Alow+\theta_1}{\mu^t} + \Vlow+\theta_2$, and $\Sup^t \triangleq \frac{\Aup+\theta_1}{\mu^t} + \Vup+\theta_2$, $\Delta_{\x}^t \triangleq \E_{\xi^{t}} [\|\x^t - \bar{\x}\|_2^2]$, $\Delta_{F}^t \triangleq \E_{\xi^{t}} [(\min_{i=1}^t F(\x^i) )  -   F(\bar{\x})]$, and $\Phi^t \triangleq \E_{\xi^t} [\tfrac{1}{2}\| \x^{t} - \bar{\x}\|_{\H_{\ast}^t}^2] -  Z_n^k  V_s\tfrac{1}{2}\|\x^t-\bar{\x}\|_2^2$.

\noi First, using \textbf{Part (c)} in Lemma \ref{lemma:first:bound}, we have: $\|\r^t\|\leq L'_F$.

\noi Second, using Lemma \ref{lemma:bound:negative:l2}, we have: $\sum_{t=1}^T \Upsilon^t \leq  C_{\Upsilon}\left(1 + \ln(T)\right) -  T  \Delta_{F}^T$ for any $T\geq 1$.

\noi Third, using the definition of $\Slow^t$, we have:
\beq \label{eq:bound:Slow:t}
\Slow^t \triangleq  \frac{\Alow+\theta_1}{\mu^t} + \Vlow+\theta_2  \geq \frac{\Alow+\theta_1}{\mu^t} + V_s = \frac{(\Alow+\theta_1)(t+t_0)}{ \eta } + V_s = V_s (t+t_0) + V_s \geq V_s (t+1).
\eeq

\noi Fourth, we establish the upper bound for $(-\Phi^{T+1}+\Phi^{1})$ using the following inequalities:
\beq \label{eq:bound:Phi:T1:1}
&& -\Phi^{T+1}+\Phi^{1} \nn\\
&=& - \{\E_{\xi^{T+1}}[\tfrac{1}{2}\| \x^{T+1} - \bar{\x}\|_{\H_{\ast}^{T+1}}^2]  -  Z_n^k\tfrac{V_s}{2}\|\x^{T+1} - \bar{\x} \|_2^2 \}\nn\\
&& + \{ \E_{\xi^{1}}[\tfrac{1}{2}\| \x^{1} - \bar{\x}\|_{\H_{\ast}^{1}}^2]  -  Z_n^k\tfrac{V_s}{2}\|\x^{1} - \bar{\x} \|_2^2 \} \nn\\
& \overset{\step{172}}{\leq}  &  -\Slow^{T+1} \E_{\xi^{T+1}}[ \tfrac{1}{2}\|  \x^{T+1} - \bar{\x} \|_2^2] +  Z_n^k \tfrac{V_s}{2}\|\x^{T+1} - \bar{\x} \|_2^2 + \Sup^1 \E_{\xi^{1}}[\tfrac{1}{2}\| \x^{1} - \bar{\x}\|_{2}^2 ] \nn\\
& \overset{\step{173}}{\leq}  & - Z_n^k [ V_s(T+2)- V_s]  \cdot \tfrac{1}{2}\|  \x^{T+1} - \bar{\x} \|_2^2 + Z_n^k \Sup^1 \tfrac{1}{2}\Delta_{\x}^1  \nn\\
& \overset{\step{174}}{\leq}  &  Z_n^k \{-   (T+1)  \tfrac{V_s}{2}\Delta_{\x}^{T+1} + K_4\},
\eeq
\noi where step \step{172} uses $\Slow^{t} \E_{\xi^{t}}[\tfrac{1}{2}\| \x^{t} - \bar{\x}\|_2^2]\leq \E_{\xi^{t}}[\tfrac{1}{2}\| \x^{t} - \bar{\x}\|_{\H_{\ast}^{t}}^2]\leq \Sup^{t} \E_{\xi^{t}}[\tfrac{1}{2}\| \x^{t} - \bar{\x}\|_2^2]$ for all $t\geq 1$; step \step{173} uses $\E_{\xi^{T+1}}[\tfrac{1}{2}\| \x^{T+1} - \bar{\x}\|_2^2]=Z_n^k \tfrac{1}{2}\| \x^{T+1} - \bar{\x}\|_2^2$ and $\Slow^{T+1}\geq  (T+2) \tfrac{V_s}{2}\| \x^{T+1} - \bar{\x}\|_2^2$; step \step{174} uses the definition of $K_4\triangleq \tfrac{1}{2}\Sup^1 \Delta_{\x}^1$.

\noi \textbf{(a)} Using the inequality in \textbf{Part (b)} in Lemma \ref{lemma:bcd:bound}, we have:
\beq \label{eq:xk1:xk:minus}
&& \E_{\xi^t}[ \tfrac{1}{2}\| \x^{t+1} - \bar{\x}\|_{\H_{\ast}^t}^2 ] -  \E_{\xi^t}[\tfrac{1}{2}\| \x^{t} - \bar{\x}\|_{\H_{\ast}^t}^2] + Z_n^k     \tfrac{V_s}{2}\| \x^{t} - \bar{\x} \|_2^2  \nn\\
&\leq&     Z_n^k  \Upsilon^t +  (2   / \Slow^t)   Z_n^k \|\r^t\|_2^2  +  ( 1+   2\kappa^t )     Z_n^k  \|\r^t\| \|\bar{\x}\| \nn\\
& \overset{\step{172}}{\leq}  &     Z_n^k  \Upsilon^t + \tfrac{2  Z_n^k  (L_F')^2 }{V_s t} +  ( 3+   2\epsilon)     Z_n^k  L_F' \|\bar{\x}\|,
\eeq
\noi where step \step{172} uses $\Slow^t \geq V_s (t+1) > V_s t$ as shown in Inequality (\ref{eq:bound:Slow:t}), $\|\r^t\|\leq L_F'$, and $\kappa^t\leq 1+\epsilon$.

Using the results in Lemma \ref{lemma:BCD:Ht1:Ht}, we have:
\beq \label{eq:dfdf}
\E_{\xi^{t+1}}[\tfrac{1}{2}\|\x^{t+1} - \bar{\x} \|_{\H_{\ast}^{t+1}}^2]  -  \E_{\xi^t}[\tfrac{1}{2}\| \x^{t+1} - \bar{\x} \|_{\H_{\ast}^{t}}^2] \leq Z_n^k  \tfrac{V_s}{2}\| \x^{t+1} - \bar{\x} \|_2^2.
\eeq

\noi We define $\Phi^t \triangleq \E_{\xi^t} [\tfrac{1}{2}\| \x^{t} - \bar{\x}\|_{\H_{\ast}^t}^2] -  Z_n^k   \tfrac{V_s}{2}\|\x^t-\bar{\x}\|_2^2$. Adding the two inequalities in (\ref{eq:xk1:xk:minus}) and (\ref{eq:dfdf}) together, we have:
\beq \label{eq:recursive}
 \Phi^{t+1} -  \Phi^t
\leq   Z_n^k  (    \tfrac{ 2  (L_F')^2 }{V_s} \cdot \tfrac{1 }{t} +  \Upsilon^t  +  ( 3 +  2\epsilon )   L_F' \|\bar{\x}\|  ).
\eeq

\noi \textbf{(b)} Let $T\geq 1$ be any integer. Summing Inequality (\ref{eq:recursive}) over $t=1,...,T$, we have:
\beq \label{eq:last:BCD:BCD}
 0 &\leq & \textstyle - \Phi^{T+1}  + \Phi^{1}   + Z_n^k \{(\sum_{t=1}^T\frac{1}{t})\cdot \frac{2 (L_F')^2 }{V_s} +  \sum_{t=1}^T\Upsilon^t + T (3+2\epsilon)L_F' \|\bar{\x}\|\} \nn\\
&\overset{\step{172}}{\leq }& - \Phi^{T+1}  + \Phi^{1}   + Z_n^k \{  (\ln(T)+1) (\tfrac{ 2  (L_F')^2 }{V_s} + C_{\Upsilon} )    -    T  \Delta_{F}^T   +  T ( 3 +  2\epsilon )     L_F' \|\bar{\x}\| \}   \nn\\
&\overset{\step{173}}{= }& - \Phi^{T+1}  + \Phi^{1}   + Z_n^k  \{    (\ln(T)+1) D_4 -  T  \Delta_{F}^T       + T ( 3 +  2\epsilon )     L_F' \|\bar{\x}\|  \} \nn\\
&\overset{\step{174}}{\leq}& Z_n^k\{  - (T+1)  \tfrac{V_s}{2}\Delta_{\x}^{T+1} + K_4 +      D_4  (1 + \ln(T) ) -   T  \Delta_{F}^T   + T ( 3 +  2\epsilon )     L_F' \|\bar{\x}\|  )\},
\eeq
\noi where step \step{172} uses $\sum_{t+1}^T\tfrac{1}{t}\leq \ln(T)+1$ and the upper bound for $\sum_{t=1}^T \Upsilon^t \leq  C_{\Upsilon}\left(1 + \ln(T)\right) -  T  \Delta_{F}^T$; step \step{173} uses the definition of $D_4 \triangleq \frac{2 (L_F')^2 }{V_s} + C_{\Upsilon}$; step \step{174} uses Inequality (\ref{eq:bound:Phi:T1:1}).

\noi We now focus on (\ref{eq:last:BCD:BCD}). Using the fact that $\Delta_{\x}^{T+1} \geq 0$, we obtain: $\Delta_{F}^T \leq \frac{K_4}{T} +  \frac{ D_4 \cdot \left(1 + \ln(T)\right) }{T}         + ( 3 +    2\epsilon )   L_F' \|\bar{\x}\| $.

\noi Using the fact that $\Delta_{F}^T \geq 0$, we obtain: $\Delta_{\x}^{T+1} \leq  \left(\frac{K_4}{T+1} +  \frac{ D_4 \cdot \left(1 + \ln(T)\right) }{T+1}         + ( 3 +    2\epsilon )   L_F' \|\bar{\x}\| \right) \frac{2}{V_s} $.

\end{proof}

\section{Experiments}
\label{app:sect:exp}

This section demonstrates the effectiveness and efficiency of Algorithm \ref{algo:main} on two nonsmooth sparsity constrained optimization tasks, namely the sparsity constrained $\ell_1$ regression and sparsity constrained $\ell_\infty$ regression. Given an arbitrary design matrix $\A \in \mathbb{R}^{m\times n}$ and an observation vector $\b \in \mathbb{R}^m$, we aim to solve the following optimization problems:
\beq
&&\min_{\x}\,\frac{\lambda}{2}\|\x\|_2^2+ \|\A\x-\b\|_{1},\,s.t.\,\|\x\|_0 \leq s,\nn\\
\text{and}\,\,&&\min_{\x}\,\frac{\lambda}{2}\|\x\|_2^2+ \|\A\x-\b\|_{\infty},\,s.t.\,\|\x\|_0 \leq s,\,\,\,\,\,\nn
\eeq

\noi where  $s$ and $\lambda$ are given parameters.

$\blacktriangleright$ \textbf{Datasets}. Following \cite{yuankdd2020}, we examine four types of datasets for the design matrix $\A \in \mathbb{R}^{m\times n}$. \textbf{(i)} `random-m-n': The matrix of size $m\times n$ is generated by sampling from a standard Gaussian distribution. \textbf{(ii)} `e2006-m-n': We select $m$ examples and $n$ dimensions from the original real-world dataset `e2006', available for download at: \url{https://www.csie.ntu.edu.tw/~cjlin/libsvmtools/datasets}. This dataset contains 16087 examples and 150360 dimensions. \textbf{(iii)} `random-m-n-C': We create a matrix $\mathcal{V}(\A)\in\mathbb{R}^{m\times n}$ to verify the robustness of the algorithms. Here, $\mathcal{V}(\A)$ is a noisy version of $\A\in\mathbb{R}^{m\times n}$, with $2\%$ of entries in $\A$ corrupted by scaling the original values by 100 times \cite{yuankdd2020}. \textbf{(iv)} `e2006-m-n-C': We employ the same method to generate corrupted real-world data as used in the `random-m-n-C' dataset. We generate the observation vector $\mathbf{b}$ in $\mathbb{R}^m$ as follows: a sparse signal $\bar{\mathbf{x}}$ in $\mathbb{R}^n$ is created by randomly selecting a support set of size $100$, with values sampled from a standard Gaussian distribution. The observation vector $\mathbf{b}$ is then computed as $\mathbf{b} = \A\bar{\mathbf{x}} + 10\times\randn(m,1)$.

$\blacktriangleright$ \textbf{Compared Methods}. We compare \textbf{SPGM-IHT} and \textbf{SPGM-BCD} with 5 state-of-the-art nonsmooth sparsity constrained optimization algorithms: \bfit{(i)} Projective Subgradient Descent (PSGD) \cite{liu2019one}, \bfit{(ii)} Alternating Direction Method of Mutipliers based on IHT (ADMM-IHT) \cite{HeY12}, \bfit{(iii)} Dual Iterative Hard Thresholding(DIHT)\cite{yuan2020jmlr}, \bfit{(iv)} Convex $\ell_1$ Approximation Method (CVX-$\ell_1$) \cite{candes2005decoding}, and \bfit{(v)} Nonconvex $\ell_p$ Approximation Method (NCVX-$\ell_p$) \cite{xu2012l}. For CVX-$\ell_1$ and NCVX-$\ell_p$, we use standard linearized ADMM to solve nonsmooth $\ell_1$ norm and $\ell_{1/2}$ norm regularized problems $\min_{\x}\,F(\x) + \sigma\|\x\|_{p}^p$ with $p\in\{1,\frac{1}{2}\}$, sweeping the regularization parameter $\sigma$ over a range or values ($\sigma = \{2^{-9}, 2^{-7},...,2^{9}\}$). We run these two algorithms for 10 parameters, selecting the solution that leads to the smallest objective after hard thresholding projection and re-optimization over the support set. We employ an efficient closed-form solver to compute the $\ell_p$ norm proximal operator \cite{xu2012l}.

$\blacktriangleright$ \textbf{Experimental Settings}. We update the smoothing parameter $\mu$ every $K=10$ iterations by halving it: $\mu \Leftarrow \mu \times \tfrac{1}{2}$. For \textbf{SPGM-BCD}, the random strategy ensures a strong optimality guarantee by maintaining the block-$k$ stationary condition. However, the greedy strategy often yields faster convergence in practice. Therefore, we combine both methods, selecting $8$ coordinates using the random strategy and $2$ coordinates using the greedy strategy \cite{yuankdd2020}. We keep a record of the relative changes of the objective function values by $d_t = |F(\x^t)-F(\x^{t+1})|/(1+|F(\x^t)|)$. We let \textbf{SPGM} run up to $T$ iterations and stop it at iteration $t<T$ if $\text{mean}([{d}_{t-\text{min}(t,\upsilon)+1},{d}_{t-min(t,\upsilon)+2},...,{d}_t]) \leq \epsilon$. We use the default value $(\theta,\epsilon,\upsilon,T)=(10^{-3},10^{-5},100,1000)$ for \textbf{SPGM}. All code was implemented in Matlab on an Intel 3.20GHz CPU with 8 GB RAM. We assess the quality of the solution by comparing the objective values across different methods. Recognizing that the optimal solution is expected to be sparse, we initialize the solutions for all methods as $10^{-3} \times \randn(n,1)$ and project them to feasible solutions. We vary $s=\{5,10,20,...,80,90\}$ for different datasets and present the average results based on 5 random initial points. 

\begin{figure}[!ht]
\centering
\begin{subfigure}{.24\textwidth}\centering\includegraphics[width=\linewidth]{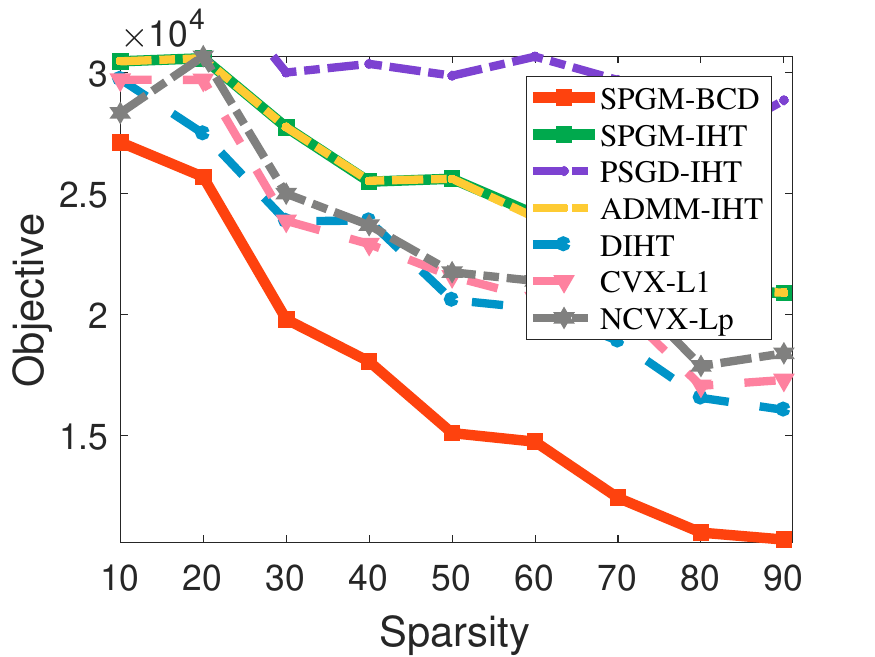}\caption{\scriptsize random-256-1024}\label{fig:sub1}\end{subfigure}
\begin{subfigure}{.24\textwidth}\centering\includegraphics[width=\linewidth]{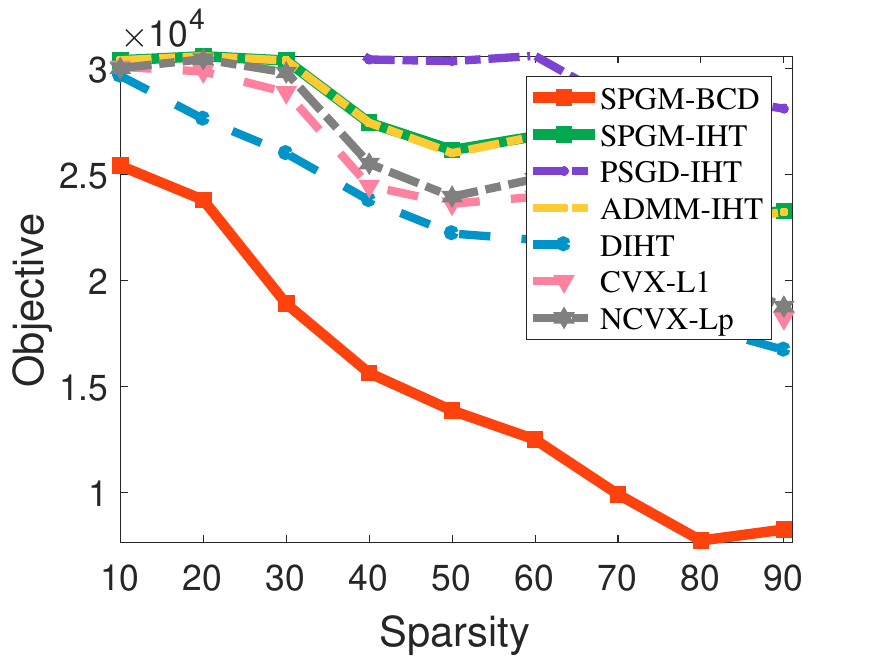}\caption{\scriptsize random-256-2048}\label{fig:sub2}\end{subfigure}
\begin{subfigure}{.24\textwidth}\centering\includegraphics[width=\linewidth]{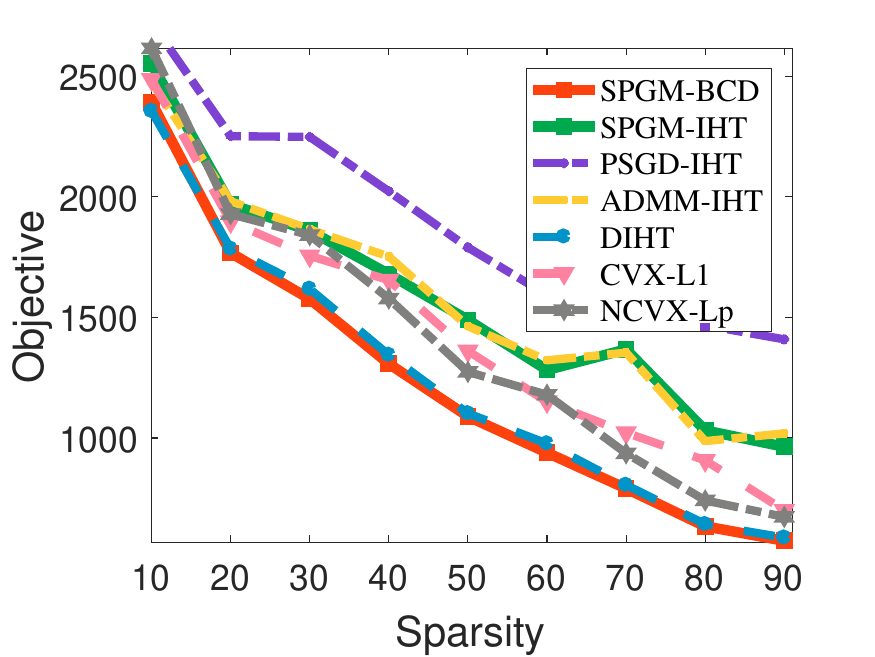}\caption{\scriptsize e2006-5000-1024}\label{fig:sub3}\end{subfigure}
\begin{subfigure}{.24\textwidth}\centering\includegraphics[width=\linewidth]{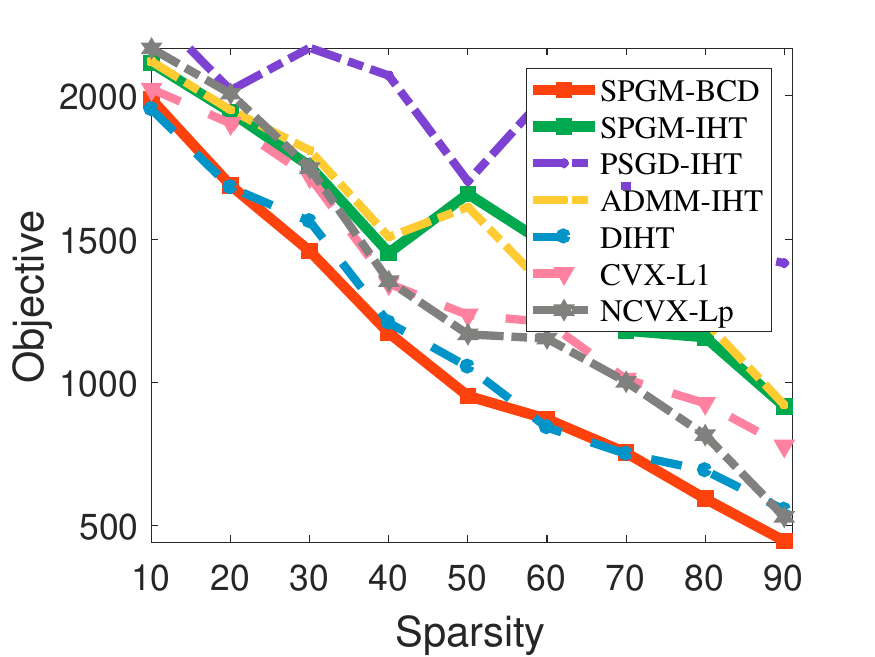}\caption{\scriptsize e2006-5000-2048}\label{fig:sub4}\end{subfigure}

\begin{subfigure}{.24\textwidth}\centering\includegraphics[width=\linewidth]{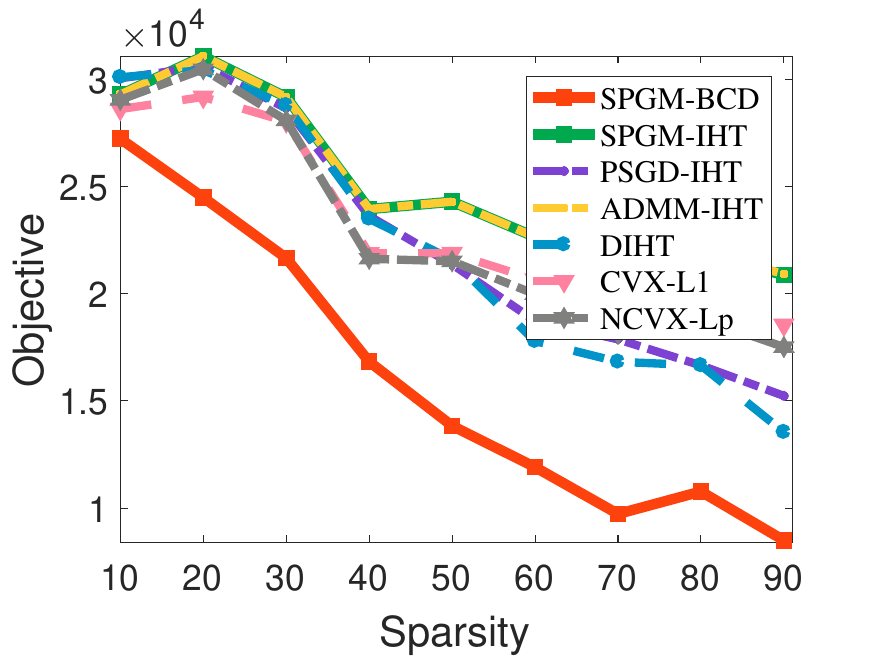}\caption{\scriptsize random-256-1024}\label{fig:sub1}\end{subfigure}
\begin{subfigure}{.24\textwidth}\centering\includegraphics[width=\linewidth]{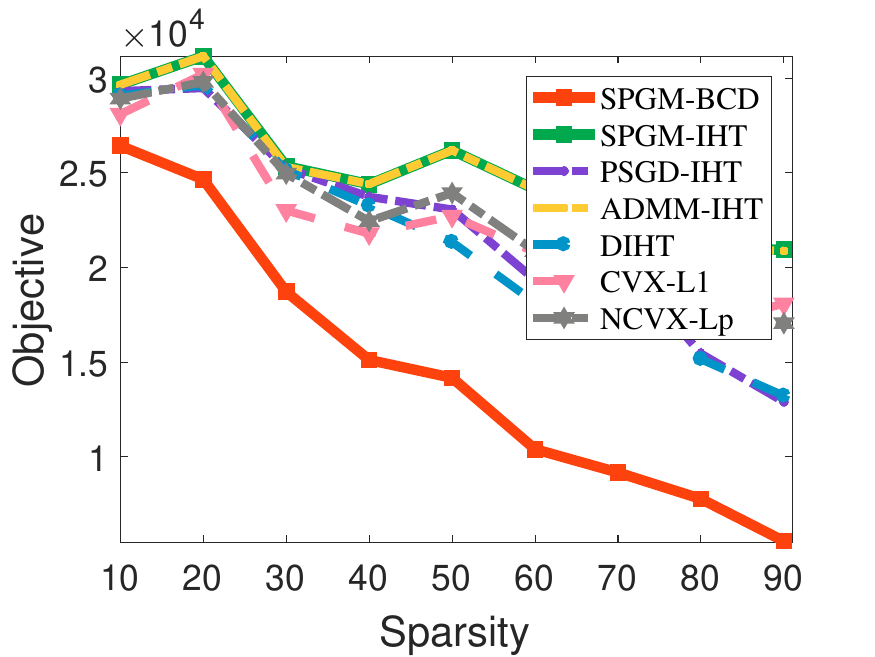}\caption{\scriptsize random-256-2048}\label{fig:sub2}\end{subfigure}
\begin{subfigure}{.24\textwidth}\centering\includegraphics[width=\linewidth]{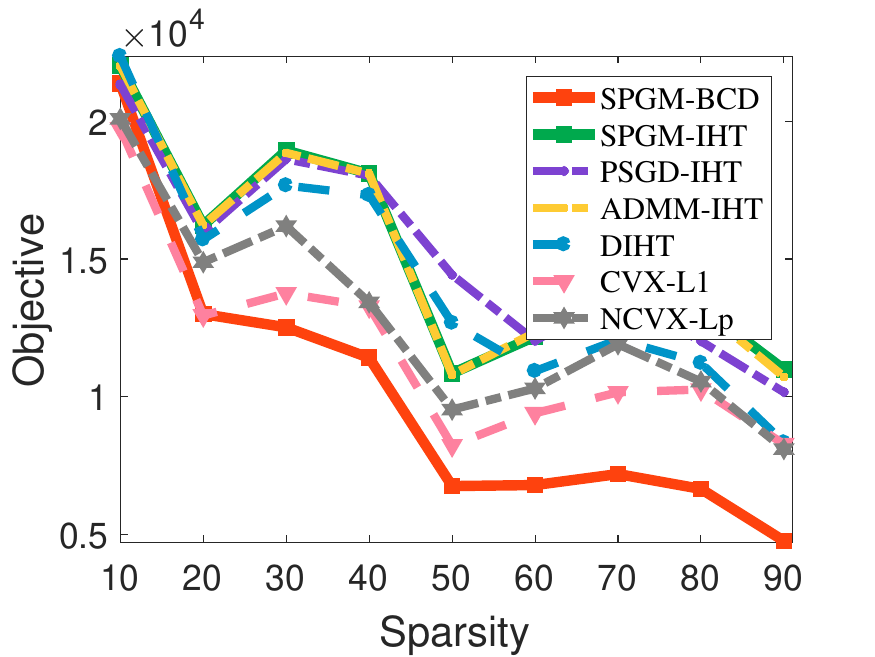}\caption{\scriptsize e2006-5000-1024}\label{fig:sub3}\end{subfigure}
\begin{subfigure}{.24\textwidth}\centering\includegraphics[width=\linewidth]{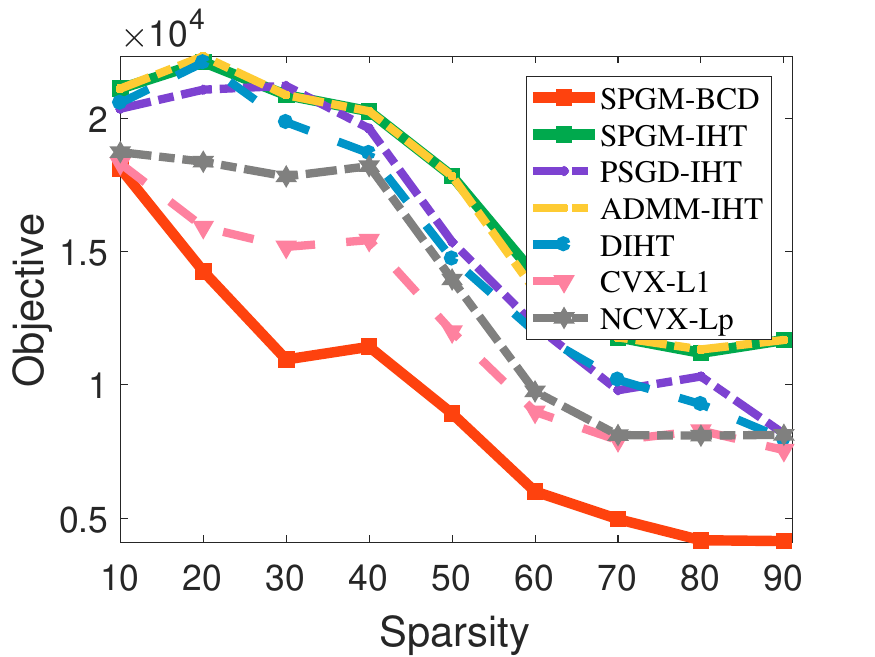}\caption{\scriptsize e2006-5000-2048}\label{fig:sub4}\end{subfigure}

\caption{Experimental results on sparsity constrained $\ell_1$ regression problems on different datasets with varying the sparsity of the solution.}
\label{exp:accuracy:1}
\end{figure}

\begin{figure}[!ht]
\centering
\begin{subfigure}{.24\textwidth}\centering\includegraphics[width=\linewidth]{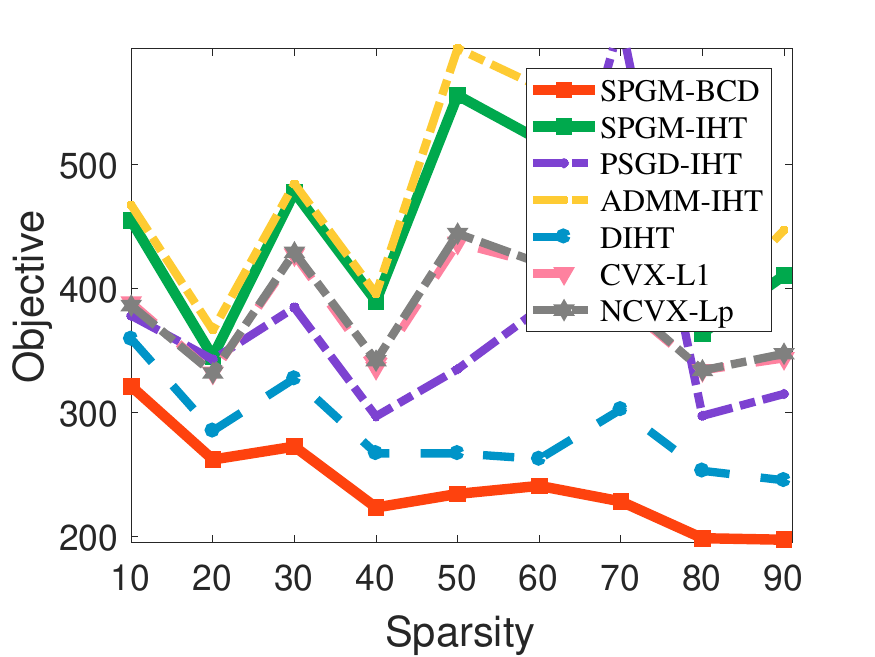}\caption{\scriptsize random-256-1024}\label{fig:sub1}\end{subfigure}
\begin{subfigure}{.24\textwidth}\centering\includegraphics[width=\linewidth]{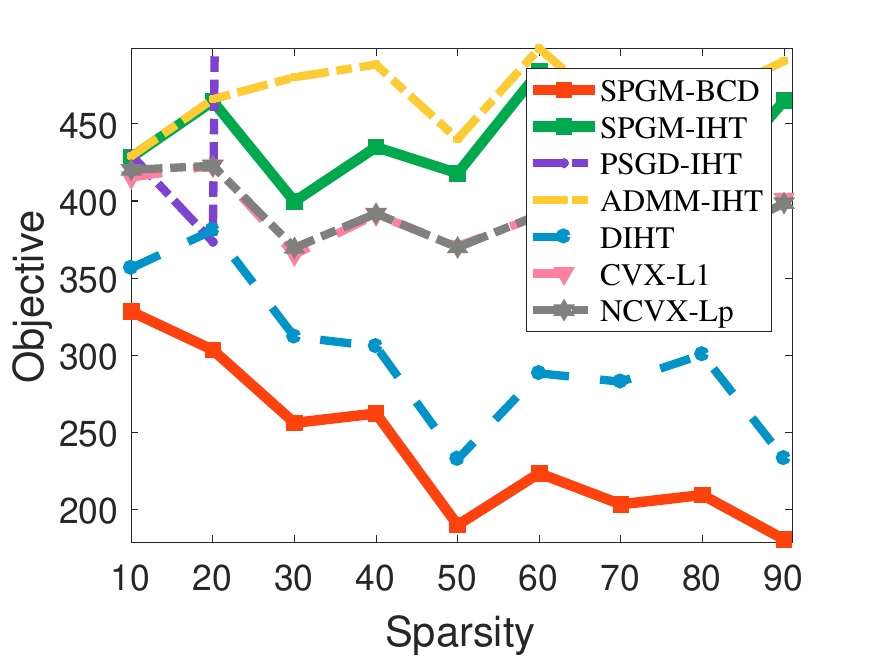}\caption{\scriptsize random-256-2048}\label{fig:sub2}\end{subfigure}
\begin{subfigure}{.24\textwidth}\centering\includegraphics[width=\linewidth]{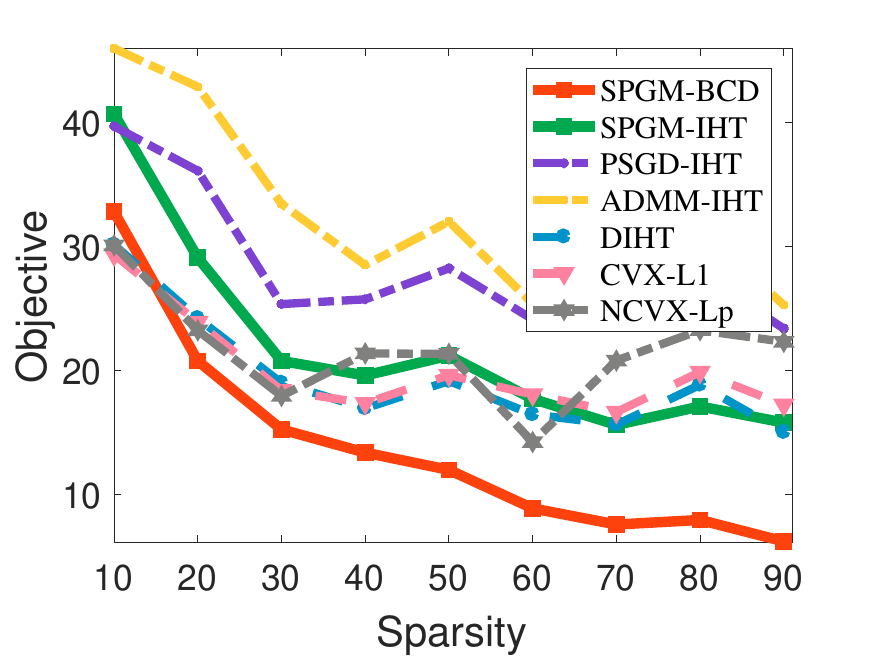}\caption{\scriptsize e2006-5000-1024}\label{fig:sub3}\end{subfigure}
\begin{subfigure}{.24\textwidth}\centering\includegraphics[width=\linewidth]{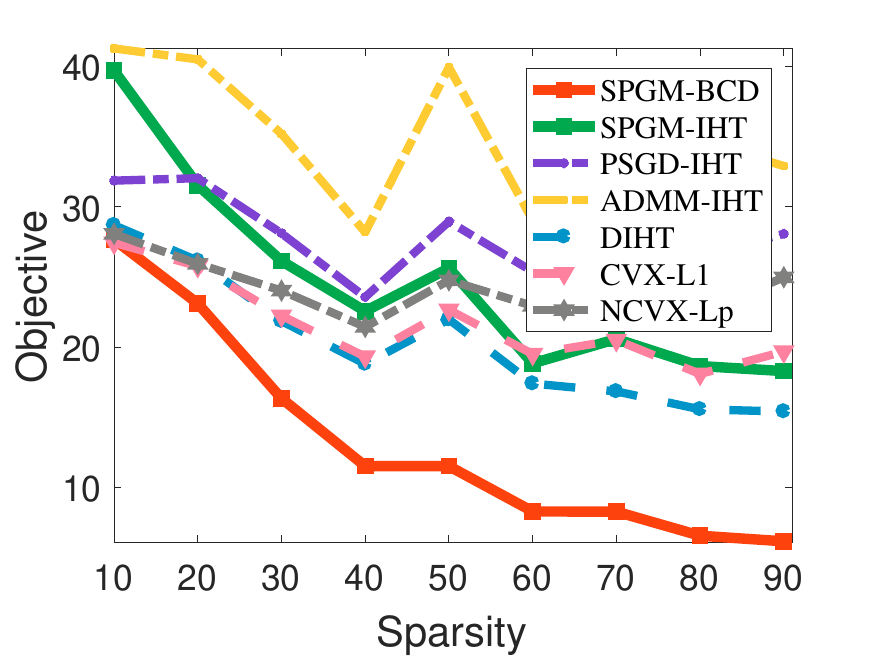}\caption{\scriptsize e2006-5000-2048}\label{fig:sub4}\end{subfigure}

\begin{subfigure}{.24\textwidth}\centering\includegraphics[width=\linewidth]{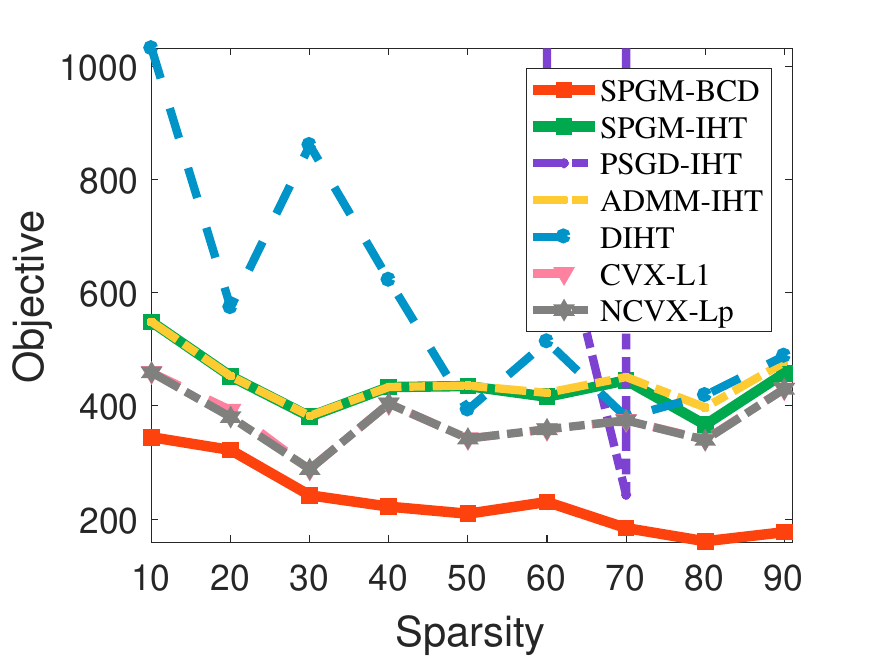}\caption{\scriptsize random-256-1024}\label{fig:sub1}\end{subfigure}
\begin{subfigure}{.24\textwidth}\centering\includegraphics[width=\linewidth]{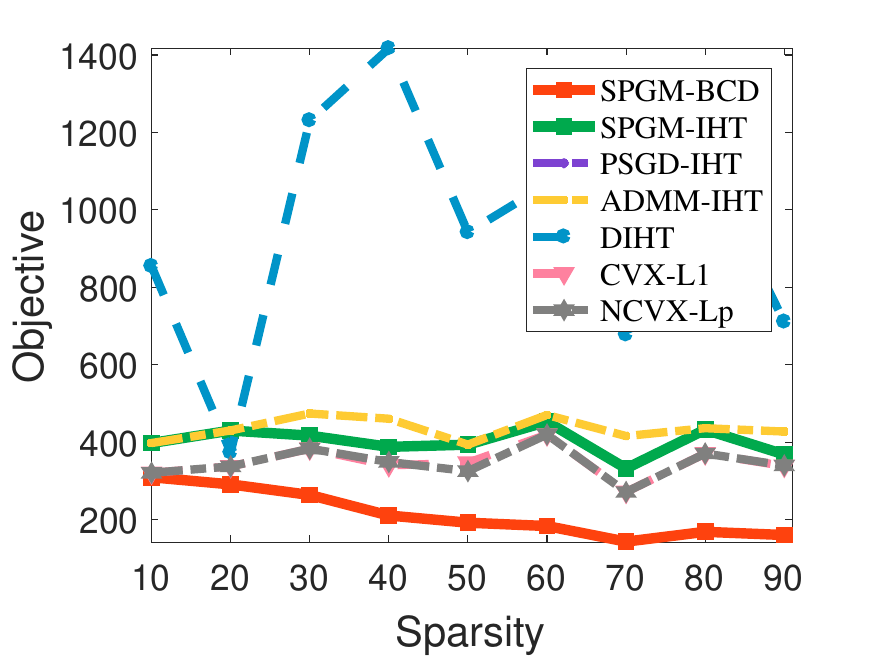}\caption{\scriptsize random-256-2048}\label{fig:sub2}\end{subfigure}
\begin{subfigure}{.24\textwidth}\centering\includegraphics[width=\linewidth]{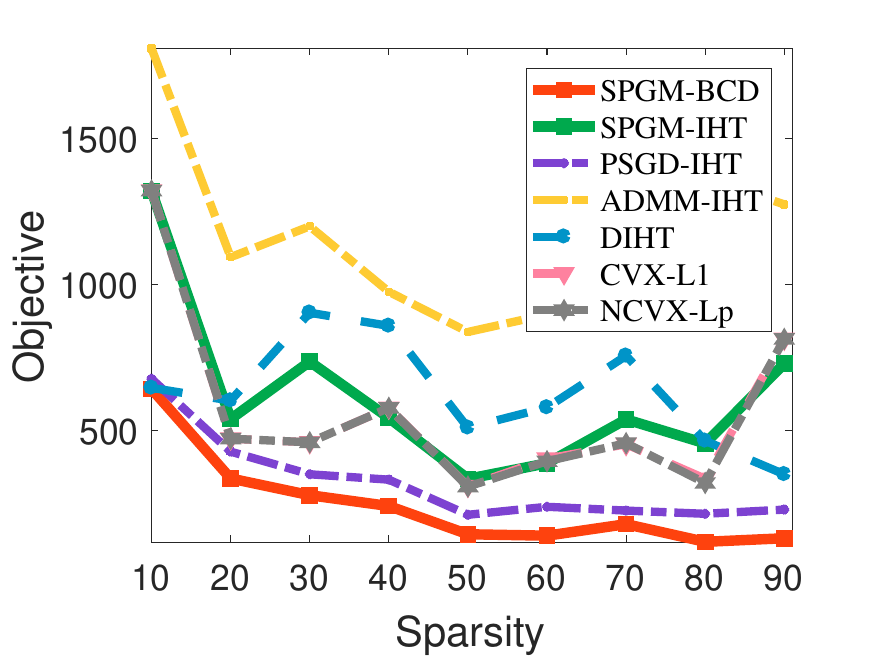}\caption{\scriptsize e2006-5000-1024}\label{fig:sub3}\end{subfigure}
\begin{subfigure}{.24\textwidth}\centering\includegraphics[width=\linewidth]{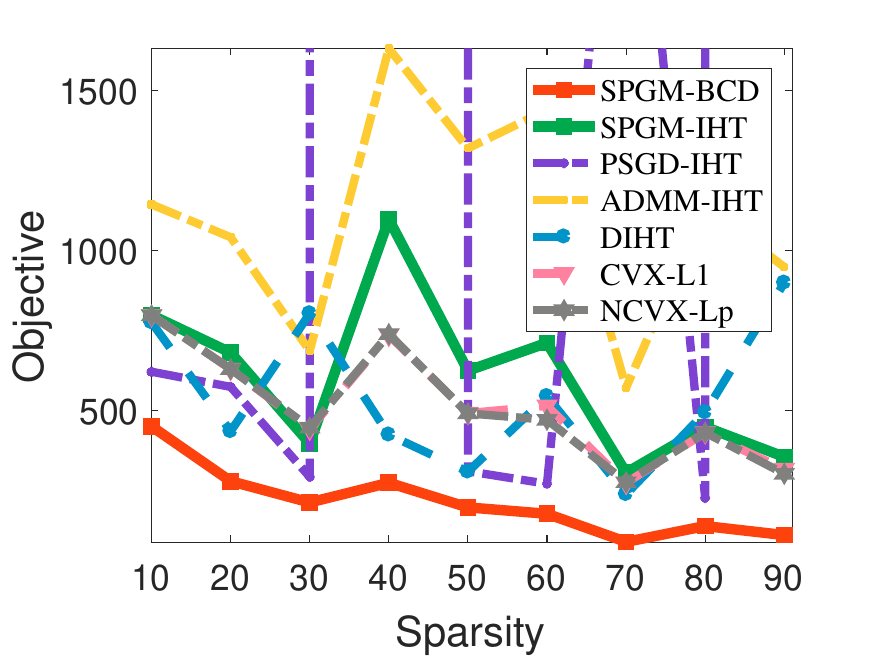}\caption{\scriptsize e2006-5000-2048}\label{fig:sub4}\end{subfigure}

\caption{Experimental results on sparsity constrained $\ell_{\infty}$ regression problems on different datasets with varying the sparsity of the solution.}
\label{exp:accuracy:2}
\end{figure}

$\blacktriangleright$ \textbf{Computational Effectiveness}. We demonstrate the computational effectiveness of \textbf{SPGM-IHT} and \textbf{SPGM-BCD} by comparing them to a set of methods (\{PSGD-IHT, ADMM-IHT, DIHT, CVX-$\ell_1$, NCVX-$\ell_p$\}). Several observations can be made from Figure \ref{exp:accuracy:1} and Figure \ref{exp:accuracy:2}. \bfit{(i)} DIHT achieves comparable results with \textbf{SPGM-BCD} on random-256-1024 and random-256-2048 in the $\ell_1$ regression. \bfit{(ii)} CVX-$\ell_1$ and NCVX-$\ell_p$ exhibit similar performance, generally outperforming others methods except \textbf{SPGM-BCD}. They achieve this by solving the relaxation problem ten times and fine-tuning the hyperparameter $\sigma$ to obtain $k$-sparsity solutions. \bfit{(iii)} PSGD-IHT generally yields worse results in our experiments. \bfit{(iv)} \textbf{SPGM-IHT} performs similarly to ADMM-IHT. \bfit{(v)} \textbf{SPGM-BCD} significantly outperforms most methods due to its ability to find stronger stationary points, which aligns with our theoretical results.

$\blacktriangleright$ \textbf{Computational Efficiency}. We present runtime comparisons for all the methods on various datasets for solving the sparsity constrained $\ell_1$ regression problem. Table \ref{tab:time:comparision} displays the average CPU times from three runs. \bfit{(i)} The convex and nonconvex relaxation methods are slightly slower than IHT-style methods because they need to run ten times to find the best regularization parameter. \bfit{(ii)} The computational efficiency of \textbf{SPGM-IHT} is comparable to that of other IHT-style methods since it is itself another IHT-style method. \bfit{(iii)} \textbf{SPGM-DEC} is slower than the other methods and typically takes about 20 seconds to converge in all instances while achieving better accuracy. \bfit{(iv)} Overall, the efficiency of both \textbf{SPGM-DEC} and \textbf{SPGM-IHT} is on par with existing methods. This is expected since they are block coordinate descent algorithms.

\begin{table} [!hb]
\begin{center}
\scalebox{0.85}{\begin{tabular}{|c|c|c|c|c|c|c|c|}
\hline
&{  PSGD-IHT} & {ADMM-IHT} &  { DIHT}  & { CVX-$\ell_1$ }&{ NCVX-$\ell_p$} &{ SPGM-IHT} &{ SPGM-BCD} \\
\hline
random-256-1024 &    $1\pm1$   &   $2\pm3$ &  $1\pm2$&       $4\pm1$          & $2\pm1$   &  $2\pm1$    & $14\pm3$  \\
random-256-2048 &    $1\pm1$   &   $2\pm1$ &  $3\pm2$&       $3\pm1$          & $2\pm1$   &  $2\pm1$    & $15\pm3$ \\
e2006-5000-1024 &   $4\pm1$   &   $2\pm1$ &  $2\pm1$&       $5\pm1$          & $4\pm1$   &  $2\pm1$    & $21\pm5$\\
e2006-5000-2048 &   $5\pm1$   &   $3\pm2$ &  $3\pm3$ &      $5\pm2$          & $4\pm1$   &  $2\pm1$    & $22\pm5$ \\
 \hline
\end{tabular}}
\caption{Comparisons of average times (in seconds) of all the methods on different datasets.}
\label{tab:time:comparision}
\end{center}
\end{table}

\end{document}